%% file: texte.tex
\documentclass[a4paper,11pt,reqno]{smfart}
\usepackage{amssymb,amsmath,mathrsfs,graphics,graphicx,mathptm,float}
\usepackage[T1]{fontenc}
\usepackage[english, francais]{babel}
\usepackage[all]{xy}

\theoremstyle{plain}
\newtheorem{thm}{Th\'eor\`eme}[section]
\newtheorem{pro}[thm]{Proposition}
\newtheorem{lem}[thm]{Lemme}
\newtheorem{cor}[thm]{Corollaire}

\theoremstyle{definition}
\newtheorem{defi}[thm]{D\'efinition}
\newtheorem{defis}[thm]{D\'efinitions}
\newtheorem{eg}[thm]{Exemple}
\newtheorem{egs}[thm]{Exemples}

\newtheorem{rem}[thm]{Remarque}

\def\og{\leavevmode\raise.3ex\hbox{$\scriptscriptstyle\langle\!\langle$~}}
\def\fg{\leavevmode\raise.3ex\hbox{~$\!\scriptscriptstyle\,\rangle\!\rangle$}}

\def\transp #1{\vphantom{#1}^{\mathrm t}\! {#1}}

\addtocounter{section}{0}             
\numberwithin{equation}{section}       

\begin{document}
\title[Quelques propri\'et\'es des transformations birationnelles du plan projectif complexe]{Quelques propri\'et\'es des transformations birationnelles du plan projectif complexe,\\
{\sl une histoire pour S.}
}

\date{}
\author{Julie D\'eserti}

\address{Institut de Math\'ematiques de Jussieu, UMR $7586$ du CNRS, Universit\'e Paris $7,$ Projet G\'eom\'etrie et Dynamique, Site Chevaleret, Case $7012,$ $75205$ Paris Cedex 13, France. 
Membre de l'ANR BLAN$06$-$3$\_$137237$}
\email{deserti@math.jussieu.fr}
\urladdr{http://people.math.jussieu.fr/~deserti/}

\maketitle{}

\begin{altabstract}
\selectlanguage{english}
We present some (unfortunately not all) known properties on the Cremona 
group; when it's possible we mentioned links with the most 
known group of polynomial automorphisms of the affine plane.
The mentioned properties are essentially algebraic properties: 
generators, relations, finite subgroups, subgroups of finite type,
automorphisms of the Cremona group, Tits alternative... but also 
dynamical properties: classification
of birational maps, centralizer, dynamic of an Heisenberg
subgroup... We deal with the construction
of entropic automorphisms on rational surfaces.
\end{altabstract}

\selectlanguage{french}

\begin{abstract}
On pr\'esente certaines (malheureusement pas toutes) 
propri\'et\'es connues du groupe de Cremona en faisant, lorsque 
c'est possible, un parall\`ele avec le groupe des 
automorphismes polynomiaux de $\mathbb{C}^2.$ Les 
propri\'et\'es abord\'ees seront essentiellement de nature alg\'ebrique:
th\'eor\`eme de g\'en\'eration, sous-groupes finis, sous-groupes de type 
fini, description du groupe d'automorphismes 
du groupe de Cremona, $\ldots$
mais aussi de nature dynamique: classification des transformations
birationnelles, centralisateur, dynamique d'un sous-groupe de 
Heisenberg $\ldots$  On \'evoque un peu les probl\`emes de 
construction d'automorphismes de type entropique sur les surfaces rationnelles.
\end{abstract}
\bigskip\bigskip

\begin{center}
\begin{minipage}{110mm}
\begin{small}
\textsl{Klein, in his 1872 Erlanger Programme, defined geometry as the 
study of those properties of figures that remain invariant under a 
particular group of transformations. ""Invariance" and "group" 
are the unifying concepts in Klein's Erlanger Programme. Groups 
of transformations had been used in geometry for many years, but 
Klein's originality consisted in reversing the roles, in making the 
group the primary object of interest and letting it operate on various 
geometries, looking for invariants" (\cite{Gr}).}
\end{small}
\end{minipage}
\end{center}

\bigskip

\bigskip

\section{Introduction}

\bigskip

\noindent D\`es la fin 
du XIX$^{\text{\`eme}}$ si\`ecle, de nombreux math\'ematiciens s'int\'eressent au
groupe des transformations birationnelles du plan 
projectif complexe appel\'e aussi groupe de Cremona: N\oe ther donne un 
th\'eor\`eme de g\'en\'eration; Bertini, 
Kantor et Wiman tentent de d\'ecrire ses sous-groupes finis 
et Castelnuovo de d\'eterminer les transformations birationnelles 
fixant une courbe de genre strictement plus grand que $1.$ Toutes ces questions 
ont \'et\'e reprises et approfondies au cours du XX$^{\text{\`eme}}$ si\`ecle; on en aborde certaines
aux \S \ref{generation} et \S \ref{gpefini}. On peut d\`es lors justifier le fait
que le groupe des transformations birationnelles de $\mathbb{P}^2(\mathbb{C})$
se distingue de celui des transformations birationnelles de $\mathbb{P}^n
(\mathbb{C})$ avec $n\geq 3$ en expliquant pourquoi il n'y a pas d'analogue 
au th\'eor\`eme de N\oe ther en dimension sup\'erieure.
Le \S \ref{decompzariski} est consacr\'e
\`a la d\'emonstration du 
fait suivant: tout \'el\'ement du groupe de Cremona s'\'ecrit comme une compos\'ee
de transformations \'el\'ementaires, les \'eclatements. Au \S \ref{autbir} on 
rappelle que la description du groupe d'automorphismes 
d'un groupe donn\'e est un probl\`eme classique; on donne aussi, entre autres, une d\'emarche qui a conduit \`a la description du groupe
des automorphismes du groupe de Cremona. 

\noindent Apr\`es avoir \'evoqu\'e certains aspects alg\'ebriques du groupe de Cremona 
on aborde des propri\'et\'es dynamiques comme la classification des 
transformations birationnelles \`a conjugaison pr\`es, la caract\'erisation
du centralisateur de certaines transformations de Cremona... on \'etablit
un parall\`ele entre le groupe des transformations birationnelles du plan 
projectif complexe et celui des automorphismes polynomiaux du plan 
(\S \ref{dyn}). Le \og programme de Zimmer\fg\, a motiv\'e l'\'etude 
des repr\'esentations de certains r\'eseaux dans le groupe de Cremona; cette 
probl\'ematique conduit aussi naturellement \`a la description des sous-groupes 
nilpotents (\S \ref{typefini}). Bien que le groupe de Cremona 
ne se plonge pas dans un groupe lin\'eaire (\S \ref{premspas}), il satisfait,
comme les groupes lin\'eaires, l'alternative de Tits; cet \'enonc\'e
s'obtient \`a partir des propri\'et\'es des sous-groupes de type fini
du groupe de Cremona (\S \ref{alttits}).
Au \S\ref{autrat} on s'int\'eresse au probl\`eme suivant: \'etant donn\'e une transformation
birationnelle de $\mathbb{P}^2(\mathbb{C})$ existe-t-il une suite finie d'\'eclatements
$\pi\colon\mathcal{X}\to\mathbb{P}^2(\mathbb{C})$ telle que l'application induite
$f_\mathcal{X}=\pi^{-1}f\pi$ soit un automorphisme de $\mathbb{P}^2(\mathbb{C})?$ 
On finit par mentionner un probl\`eme ouvert: le groupe de Cremona est-il simple ? On donne des r\'ef\'erences o\`u cette question a \'et\'e abord\'ee sans \^etre r\'esolue; on \'evoque aussi les r\'esultats obtenus pour le groupe des automorphismes polynomiaux du plan affine.
\bigskip

\section{Premiers pas dans le groupe de Cremona}\label{premspas}

\bigskip

\noindent Une \textsl{transformation rationnelle} du plan projectif complexe dans 
lui-m\^eme est une transformation de la forme
\begin{align*}
&\mathbb{P}^2(\mathbb{C})\dashrightarrow\mathbb{P}^2(\mathbb{C}),&&
(x:y:z)\mapsto(f_0(x,y,z):f_1(x,y,z):f_2(x,y,z)),
\end{align*}
\noindent les $f_i$ d\'esignant des polyn\^omes homog\`enes de m\^eme degr\'e.

\medskip

\noindent Une \textsl{transformation birationnelle} $f=(f_0:f_1:f_2)$ de 
$\mathbb{P}^2(\mathbb{C})$ dans $\mathbb{P}^2(\mathbb{C})$ est 
une transformation rationnelle qui admet un inverse lui-m\^eme rationnel.

\noindent Le \textsl{degr\'e} de $f,$ not\'e $\deg f,$ 
est le degr\'e des~$f_i;$ c'est aussi
si $L$ d\'esigne une droite g\'en\'erique du plan projectif complexe $\deg 
f^{-1}L.$

\medskip

\noindent Le \textsl{groupe de Cremona} est le groupe des transformations 
birationnelles de $\mathbb{P}^2(\mathbb{C})$ dans lui-m\^eme; on le note
$\mathrm{Bir}(\mathbb{P}^2).$ Ses \'el\'ements sont aussi
appel\'es \textsl{transformations de Cremona}.

\begin{egs}
\begin{itemize}
\item Tout automorphisme de $\mathbb{P}^2(\mathbb{C}),$ {\it i.e.}
tout \'el\'ement de $\mathrm{PGL}_3(\mathbb{C}),$
est une transformation birationnelle.

\item La transformation 
\begin{align*}
& \sigma\colon\mathbb{P}^2(\mathbb{C})\dashrightarrow\mathbb{P}^2(
\mathbb{C}), && (x:y:z)\mapsto(yz:xz:xy)
\end{align*}
\noindent est rationnelle. Dans la carte affine $z=1$ on a $\sigma=\left(\frac{1}{x},
\frac{1}{y}\right);$ on remarque en particulier que~$\sigma$ est une involution  birationnelle. Elle porte souvent le nom \textsl{d'involution
de Cremona}; comme on le verra elle joue un r\^ole particulier.

\item Un \textsl{automorphisme 
polynomial de $\mathbb{C}^2$} est une application bijective de la forme suivante
\begin{align*}
& \mathbb{C}^2\to\mathbb{C}^2, && (x,y)\mapsto(f_1(x,y),f_2(x,y)), && f_i\in\mathbb{C}[x,y].
\end{align*} 

\noindent Le groupe de Cremona contient le groupe $\mathrm{Aut}[\mathbb{C}^2]$ 
des automorphismes polynomiaux du plan complexe; en particulier le prolongement 
d'un  \textsl{automorphisme \'el\'ementaire}
\begin{align*}
(\alpha x+P(y),\beta y+\gamma), && \alpha,\hspace{0.1cm}\beta\in\mathbb{C}^*,
\hspace{0.1cm}\gamma\in\mathbb{C},\hspace{0.1cm} P\in\mathbb{C}[y],\hspace{0.1cm}\deg P\geq 2,
\end{align*}
\noindent au plan projectif complexe est une transformation birationnelle. Il en est
de m\^eme pour le prolongement d'une \textsl{application de H\'enon
g\'en\'eralis\'ee}
\begin{align*}
& (y,P(y)-\alpha x), &&\alpha\in\mathbb{C}^*,
&& P\in\mathbb{C}[y], && \deg P\geq 2.
\end{align*}
\end{itemize}
\end{egs}

\begin{defis}
Soit $f=(f_0:f_1:f_2)$ dans $\mathrm{Bir}(\mathbb{P}^2).$
Le \textsl{lieu d'ind\'etermination} de $f,$ ou encore l'ensemble des \textsl{points 
\'eclat\'es} par $f,$ est le lieu d'annulation des $f_i.$ On le d\'esigne par~$\mathrm{Ind}(f).$

\noindent Le \textsl{lieu exceptionnel} de $f,$ ou encore l'ensemble des \textsl{courbes 
contract\'ees} par $f,$ est le lieu des z\'eros du d\'eterminant jacobien 
de $f.$ On le note $\mathrm{Exc}(f).$
\end{defis}

\medskip

\begin{egs}\label{egsptdind}
\begin{itemize}
\item Si $f$ appartient \`a $\mathrm{Aut}(\mathbb{P}^2(\mathbb{C}))=\mathrm{PGL}_3
(\mathbb{C})$ on a $\mathrm{Ind}(f)=\mathrm{Exc}(f)=~\emptyset.$

\item Rappelons que l'involution de Cremona s'\'ecrit
\begin{align*}
& \sigma\colon\mathbb{P}^2(\mathbb{C})\dashrightarrow\mathbb{P}^2(
\mathbb{C}), && (x:y:z)\mapsto(yz:xz:xy);
\end{align*}
\noindent on remarque que 
\begin{align*}
&\mathrm{Ind}(\sigma)=\{(1:0:0),\,(0:1:0),\,
(0:0:1)\} &&\text{et} &&\mathrm{Exc}(\sigma)=\{x=0\}\cup\{y=~0\}\cup\{z=~0\}.
\end{align*}
 
\item Si $\rho$ est l'involution donn\'ee par
\begin{align*}
& \rho\colon\mathbb{P}^2(\mathbb{C})\dashrightarrow\mathbb{P}^2(
\mathbb{C}), && (x:y:z)\mapsto(xy:z^2:yz),
\end{align*}
\noindent on a $\mathrm{Ind}(\rho)=\{(1:0:0),\,(0:1:0)\}$ et $\mathrm{Exc}
(\rho)=\{y=0\}\cup\{z=0\}.$

\item Enfin l'involution $\tau$ d\'efinie par 
\begin{align*}
& \tau\colon\mathbb{P}^2(\mathbb{C})\dashrightarrow\mathbb{P}^2(
\mathbb{C}), && (x:y:z)\mapsto(x^2:xy:y^2-xz)
\end{align*}
\noindent satisfait $\mathrm{Ind}(\tau)=\{(0:0:1)\}$ et $\mathrm{Exc}
(\tau)=\{x=0\}.$
\end{itemize}
\end{egs}

\noindent \`A une transformation de Cremona $f$ on peut associer son graphe $\Gamma_f$ qui est une
sous-vari\'et\'e irr\'eductible de $\mathbb{P}^2(\mathbb{C})\times\mathbb{P}^2
(\mathbb{C})$ (\emph{voir} \cite{Fisch}).
Soit $p_1$ (resp. $p_2$) la projection de $\Gamma_f$ sur le premier (resp. 
second) facteur:
\begin{align*}
\xymatrix{&\Gamma_f\ar[dl]_{p_1}\ar[dr]^{p_2} &\\
\mathbb{P}^2(\mathbb{C})\ar@{-->}[rr]^f & & \mathbb{P}^2(\mathbb{C}) }
\end{align*}
\bigskip

\noindent Le lieu d'ind\'etermination de $f$ est l'ensemble (fini) des points o\`u
$p_1$ n'est pas localement inversible. Notons $\mathcal{E}(p_2)$ le lieu des points
o\`u $p_2$ n'est pas une application finie; $\mathrm{Exc}(f)$ 
co\"{\i}ncide avec~$p_1(\mathcal{E}(p_2)).$

\bigskip

\noindent Dans ce qui suit on verra que $\mathrm{Bir}(\mathbb{P}^2)$
poss\`ede de nombreuses propri\'et\'es du groupe lin\'eaire, n\'eanmoins on a 
l'\'enonc\'e suivant.

\begin{pro}[\cite{CD}]
{\sl Le groupe de Cremona ne se plonge pas dans $\mathrm{GL}_n
(\Bbbk)$ o\`u $\Bbbk$ d\'esigne un corps de caract\'eristique nulle.}
\end{pro}

\begin{proof}[{\sl D\'emonstration}]
Commen\c{c}ons par rappeler le r\'esultat suivant d\^u \`a Birkhoff: 
soient~$\Bbbk$ un corps de caract\'eristique nulle et $A,$ $B,$ $C$ trois
\'el\'e\-ments de~$\mathrm{GL}_n(\Bbbk)$ tels que $C$ soit d'ordre $p$ premier 
et~$[A,B]=~C,$ $[A,C]=[B,C]=\mathrm{id};$ alors $p\leq n$ (\emph{voir} \cite{Bi}).

\medskip

\noindent Supposons qu'il existe un morphisme injectif $\iota$ de
$\mathrm{Bir}(\mathbb{P}^2)$ dans
$\mathrm{GL}_n(\Bbbk).$ Pour tout nombre premier $p$
consid\'erons, dans la carte affine $z=1,$ le groupe
\begin{align*}
\langle\left(\exp\left(-\frac{2\mathrm{i}\pi}{p}\right)x,y\right),(x,xy),
\left(x,\exp\left(\frac{2\mathrm{i}\pi}{p}\right)y\right)\rangle.
\end{align*}
\noindent Les images par $\iota$ des trois g\'en\'erateurs satisfont le
lemme de Birkhoff donc $p\leq n;$ ceci \'etant valable
pour tout premier $p,$ on obtient le r\'esultat annonc\'e.
\end{proof}

\bigskip

\subsection{Le sous-groupe des automorphismes polynomiaux du plan}

\bigskip

\noindent Introduisons deux sous-groupes de $\mathrm{Aut}[\mathbb{C}^2]$
qui jouent un r\^ole particulier
\smallskip

\begin{itemize}
\item le \textsl{groupe affine} 
\begin{align*}
\mathtt{A}=\{(\alpha_1x+\beta_1y+\gamma_1,\alpha_2x+\beta_2y+\gamma_2)\,\vert\, \alpha_i,\hspace{0.1cm} \beta_i,\hspace{0.1cm} \gamma_i\in\mathbb{C},\hspace{0.1cm}\alpha_1\beta_2
-\alpha_2\beta_1\not=0\}
\end{align*}

\item et le \textsl{groupe \'el\'ementaire}
\begin{align*}
\mathtt{E}=\{(\alpha x+P(y),\beta y+\gamma)\,\vert\, \alpha,\hspace{0.1cm}\beta\in\mathbb{C}^*,
\hspace{0.1cm}\gamma\in\mathbb{C},\hspace{0.1cm} P\in\mathbb{C}[y],\hspace{0.1cm}\deg P\geq 2\}.
\end{align*}
\end{itemize} 

\smallskip

\noindent En $1942$ Jung a d\'emontr\'e le th\'eor\`eme de structure suivant. 

\begin{thm}[\cite{Ju}]\label{jung}
{\sl Le groupe $\mathrm{Aut}[\mathbb{C}^2]$ est un produit amalgam\'e
\begin{align*}
& \mathrm{Aut}[\mathbb{C}^2]=\mathtt{A}\ast_\mathtt{S}\mathtt{E}, && \mathtt{S}=
\mathtt{A}\cap \mathtt{E}.
\end{align*}}
\end{thm}

\noindent Dit autrement tout \'el\'ement $\phi$ de $\mathrm{Aut}[\mathbb{C}^2]$ n'appartenant
pas \`a $\mathtt{S}$ est de la forme 
\begin{align*}
& (a_1)e_1\ldots a_n(e_n), && a_i\in\mathtt{A}\setminus\mathtt{E}, && 
e_i\in\mathtt{E}\setminus\mathtt{A};
\end{align*}
\noindent de plus cette \'ecriture est unique modulo les relations suivantes
\begin{align*}
& a_ie_i=(a_is)(s^{-1}e_i), && e_{i-1}a_i=(e_{i-1}s')(s'^{-1}a_i), && s,\hspace{0.1cm} s'\in
\mathtt{S}.
\end{align*}

\bigskip

\noindent Il y a de nombreuses d\'emonstrations du th\'eor\`eme
\ref{jung}. On
renvoie \`a \cite{La3} pour un historique d\'etaill\'e et pour une preuve
originale et g\'eom\'etrique de cet \'enonc\'e dont on va juste
\'evoquer l'id\'ee. Soit 
\begin{align*}
\widetilde{f}\colon(x,y)\mapsto(\widetilde{f}_1(x,y),\widetilde{f}_2(x,y))
\end{align*}
\noindent un automorphisme polynomial de $\mathbb{C}^2$ de degr\'e $n;$ on peut prolonger $\widetilde{f}$ 
en une transformation de Cremona
\begin{align*}
& f\colon\mathbb{P}^2(\mathbb{C})\dashrightarrow\mathbb{P}^2(\mathbb{C}), &&
(x:y:z)\mapsto\left(z^n\widetilde{f}_1\left(\frac{x}{z},\frac{y}{z}\right):z^n
\widetilde{f}_2\left(\frac{x}{z},\frac{y}{z}\right):z^n\right).
\end{align*}
\noindent La d\'emonstration de Lamy s'effectue par r\'ecurrence sur le nombre
de points d'ind\'etermination de~$f:$ il s'agit de montrer qu'il existe
$\varphi\colon\mathbb{P}^2(\mathbb{C})\dashrightarrow\mathbb{P}^2(\mathbb{C}),$
prolongement d'un automorphisme polynomial de $\mathbb{C}^2,$ 
tel que $\#\,\mathrm{Ind}(f\varphi^{-1})<\#\,\mathrm{Ind}(f).$

\bigskip

\noindent Puisque $\mathrm{Aut}[\mathbb{C}^2]$ est le produit
amalgam\'e de $\mathtt{A}$ et $\mathtt{E}$ le long de $\mathtt{A}\cap
\mathtt{E}$ la th\'eorie de Bass-Serre assure que 
$\mathrm{Aut}[\mathbb{C}^2]$ agit de mani\`ere non triviale 
par translation \`a gauche sur un arbre $\mathcal{T}$ (\emph{voir}~\cite{Se}). L'arbre $\mathcal{T}$ est d\'efini comme suit: 
l'ensemble des sommets est l'union disjointe des classes \`a 
gauche $(\mathrm{Aut}[\mathbb{C}^2])/\mathtt{A}$ et $(\mathrm{Aut}[\mathbb{C}^2])/\mathtt{E}$ et celui 
des ar\^etes l'union disjointe des classes \`a gau\-che~$(\mathrm{Aut}[\mathbb{C}^2])/
\mathtt{S}.$ Pour tout automorphisme polynomial $f$ l'ar\^ete 
$f\mathtt{S}$ relie les sommets $f\mathtt{A}$ et~$f\mathtt{E}.$
Voici quelques sommets de $\mathcal{T}$  

\begin{figure}[H]
\begin{center}
\input{arbre.pstex_t}
\end{center}
\end{figure}

\noindent o\`u $a,$ $\widetilde{a}$ (resp. $e,$ $\widetilde{e}$)
d\'esignent des \'el\'ements de $\mathtt{A}$ (resp. $\mathtt{E}$).

\noindent On r\'ecup\`ere ainsi une repr\'esentation fid\`ele\footnote{\hspace{0.1cm} 
Pour certains produits amalgam\'es l'action induite n'est pas
fid\`ele: on peut v\'erifier que la matrice 
$-\mathrm{id}$ agit trivialement sur l'arbre associ\'e au produit amalgam\'e $\mathrm{SL}_2(\mathbb{Z})=
\mathbb{Z}/4\mathbb{Z}\ast_{\mathbb{Z}/2\mathbb{Z}}
\mathbb{Z}/6\mathbb{Z}.$} de $\mathrm{Aut}[\mathbb{C}^2]$
dans le groupe des isom\'etries de l'arbre~$\mathcal{T}.$
En \'etudiant cette action Lamy a d\'emontr\'e de 
nombreuses propri\'et\'es pour~$\mathrm{Aut}[
\mathbb{C}^2]$ (\emph{voir} \cite{La}), on en pr\'ecisera certaines dans la 
suite.

\bigskip

\subsection{Le groupe de de Jonqui\`eres}

\noindent Le \textsl{groupe de de Jonqui\`eres} est le groupe des 
transformations birationnelles du plan projectif complexe pr\'eservant un pinceau 
de courbes rationnelles; on le note $\mathrm{dJ}.$ Comme deux pinceaux de courbes rationnelles sont birationnellement conjugu\'es, $\mathrm{dJ}$ ne d\'epend pas, \`a conjugaison pr\`es, du pinceau choisi. Dit autrement on peut supposer \`a
conjugaison birationnelle pr\`es que $\mathrm{dJ}$ est, dans une carte affine $(x,y)$ de 
$\mathbb{P}^2(\mathbb{C}),$ le groupe maximal des transformations birationnelles
laissant la fibration $y=$ cte invariante. Une transformation $f$ de $\mathrm{dJ}$
permute les fibres de la fibration, $f$ induit un automorphisme de la base
$\mathbb{P}^1(\mathbb{C})$ soit un \'el\'ement de $\mathrm{PGL}_2(
\mathbb{C});$ lorsque $f$ pr\'eserve les fibres, $f$ agit comme une homographie dans celles-ci. Le 
groupe de de Jonqui\`eres s'identifie donc au produit semi-direct
$\mathrm{PGL}_2(\mathbb{C}(y))\rtimes\mathrm{PGL}_2(\mathbb{C}).$ 

\bigskip

\subsection{R\'eseaux homalo\"idaux}

\bigskip

\noindent On mentionne le lien entre 
r\'eseaux homalo\"idaux et transformations birationnelles, ces r\'eseaux ayant
conduit l'\'ecole de g\'eom\'etrie italienne \`a l'\'etude des transformations
birationnelles.

\begin{defis}
\noindent Soient $\mathrm{S}$ une surface et $p$ un point de $\mathrm{S}.$ Il existe une surface 
$\widetilde{\mathrm{S}}$ et un morphisme~$\pi\colon
\widetilde{\mathrm{S}}\to \mathrm{S},$ uniques \`a isomorphisme pr\`es, tels que
\smallskip
\begin{itemize}
\item $\pi_{\vert\pi^{-1}(\mathrm{S}\setminus\{p\})}\colon\pi^{-1}
(\mathrm{S}\setminus\{p\})\to\mathrm{S}\setminus\{p\}$ soit un isomorphisme;

\item $E:=\pi^{-1}(p)$ soit isomorphe \`a $\mathbb{P}^1(\mathbb{C}).$
\end{itemize}
\smallskip

\noindent On dit que $\pi$ est l'\textsl{\'eclatement}
de $\mathrm{S}$ en $p$ et $E$ la \textsl{courbe exceptionnelle} ou 
\textsl{diviseur exceptionnel}; l'inverse $\pi^{-1}$ de 
$\pi$ \textsl{contracte} $E.$

\noindent Consid\'erons une courbe irr\'eductible $\mathcal{C}$ sur $\mathrm{S}$ passant
par $p$ avec multiplicit\'e $\nu.$ L'adh\'erence de~$\pi^{-1}(\mathcal{C}\setminus
\{p\})$ dans $\widetilde{\mathrm{S}}$ est une courbe irr\'eductible $\widetilde{\mathcal{C}}$ 
appel\'ee \textsl{transform\'ee stricte} de 
$\mathcal{C}$ dans~$\widetilde{\mathrm{S}}.$ Soient $D$ et $D'$ deux diviseurs sur $\mathrm{S};$
nous avons 
\begin{equation}\label{blow}
\pi^*\mathcal{C}=\widetilde{\mathcal{C}}+\nu E,\hspace{5mm}(\pi^*D,\pi^*D')=(D,D'), 
\hspace{5mm} (E,\pi^*D)=0
\hspace{5mm}\text{et}\hspace{5mm}E^2=-1.
\end{equation}

\noindent Soient $\mathrm{S}$ une surface et $p$ un point de $\mathrm{S}.$ Le diviseur exceptionnel
obtenu en \'eclatant $p$ est appel\'e \textsl{premier voisinage infinit\'esimal}
de $p$ et les points de $E$ sont dits
\textsl{infiniment proches} de $p.$ Le
\textsl{$i$-\`eme voisinage infinit\'esimal} de $p$ est l'ensemble
des points contenus dans le premier
voisinage d'un certain point du $(i-1)$-\`eme voisinage infinit\'esimal de $p.$
Lorsqu'on souhaite distinguer les points de $\mathrm{S}$ des points infiniment proches on
appelle \textsl{points propres} les points de $\mathrm{S}.$
\end{defis}

\noindent Soit $f=(f_0:f_1:f_2)$ une transformation birationnelle du plan 
projectif dans lui-m\^eme. Le \textsl{r\'eseau homalo\"idal} associ\'e \`a $f$ est 
le syst\`eme de courbes~$\mathscr{R}_f$ d\'efini par
\begin{align*}
&\alpha_0f_0+\alpha_1f_1+\alpha_2f_2,&& (\alpha_0:\alpha_1:
\alpha_2)\in\mathbb{P}^2(\mathbb{C}).
\end{align*}
\noindent Les points base $p_i$ du r\'eseau $\mathscr{R}_f$ sont les points
par lesquels passent toutes les courbes du r\'eseau; on dit aussi que les $p_i$ sont les \textsl{points 
base} de $f.$ Ils peuvent \^etre dans $\mathbb{P}^2(\mathbb{C})$ ou
infiniment proches
de $\mathbb{P}^2(\mathbb{C});$ si $p_j$ n'est pas
propre, il appartient \`a un~$k$-i\`eme voisinage infinit\'esimal de $p_\ell,$ $\ell\not=j.$ Les points
base propres de $f$ sont les points d'ind\'etermination de~$f.$ La multiplicit\'e de $f$ en $p_i$ est la multiplicit\'e d'une
courbe g\'en\'erique de~$\mathscr{R}_f$ en $p_i,$ {\it i.e.} l'ordre en~$p_i$ d'un \'el\'ement
 g\'en\'erique de~$\mathscr{R}_f.$

\noindent On a vu qu'\`a toute transformation birationnelle on peut associer un r\'eseau
 homalo\"idal, r\'eciproquement on a l'\'enonc\'e suivant.

\begin{thm}[\cite{RS}]
{\sl Un r\'eseau homalo\"idal d\'efinit une infinit\'e de transformations birationnelles, chacune
pouvant \^etre obtenue \`a partir d'une autre via composition \`a gauche par un automorphisme de $\mathbb{P}^2(\mathbb{C}).$}
\end{thm}

\begin{egs}\hspace{1mm}
\begin{itemize}
\item Consid\'erons la transformation birationnnelle $\sigma.$
Les points d'ind\'etermination de $\sigma$ sont $(1:0:0),$ $(0:1:0)$ et $(0:0:1).$ Les
transformations de Cremona associ\'ees au r\'eseau homalo\"{\i}dal constitu\'e
des coniques passant par les points $(1:0:0),$ $(0:1:0)$ et~$(0:0:1)$ sont les 
transformations de la forme $A\sigma$ avec $A$ dans $\mathrm{Aut}(\mathbb{P}^2(
\mathbb{C})).$\smallskip

\item Soit $\mathcal{S}$ le r\'eseau homalo\"{\i}dal form\'e des coniques passant par
$(1:0:0),$ $(0:1:0)$ et tangentes \`a la droite d'\'equation $z=0.$ Les transformations
birationnelles associ\'ees \`a $\mathcal{S}$ sont du type~$A\rho$ avec $A$ automorphisme 
de $\mathbb{P}^2(\mathbb{C}).$
\end{itemize}
\end{egs}

\noindent Soient $f$ une transformation birationnelle et $\mathscr{R}_f$ le r\'eseau
homalo\"idal associ\'e \`a $f.$ Les courbes de~$\mathscr{R}_f$ satisfont les \'equations
suivantes (\cite{AC})
\begin{align*}
&\sum_{i=1}^q\mu_i^2=n^2-1,&&\sum_{i=1}^q\mu_i=3n-3
\end{align*}
\noindent o\`u $\mu_i$ d\'esigne la multiplicit\'e aux points base (qui rappelons-le ne sont pas 
n\'ecessairement propres),~$q$ le nombre de points base et $n$ le degr\'e de $f.$ La
premi\`ere \'equation traduit le fait que deux courbes g\'en\'eriques du r\'eseau
homalo\"{\i}dal se coupent en les points base et un unique autre point; la seconde
exprime que les courbes du r\'eseau sont rationnelles.

\bigskip

\section{Le th\'eor\`eme de factorisation de Zariski}\label{decompzariski}

\bigskip

\noindent L'involution $\sigma$ se d\'ecompose en deux suites d'\'eclatements 

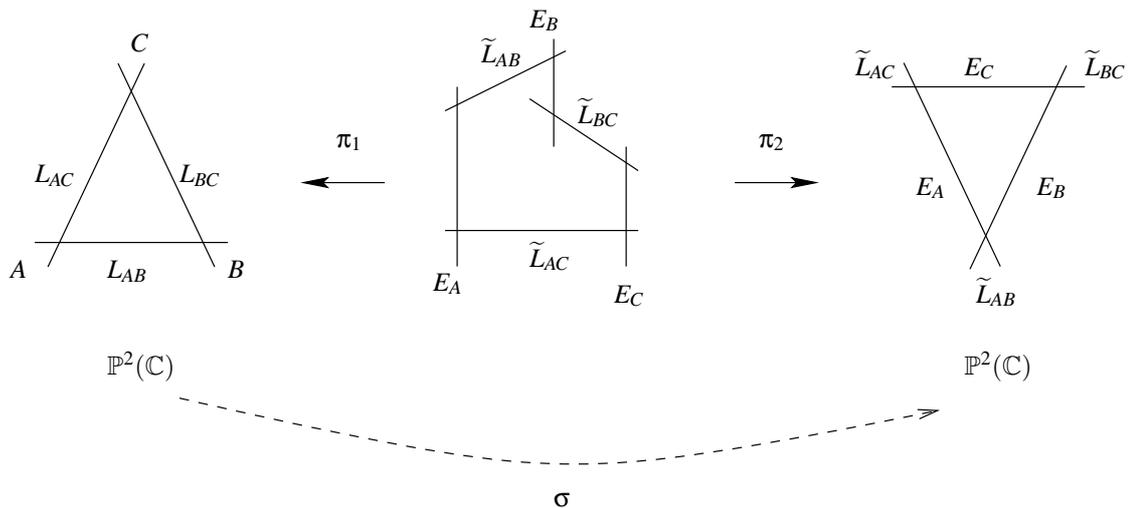
\begin{figure}[H]
\begin{center}
\input{decomp.pstex_t}
\end{center}
\end{figure}

\noindent avec
\begin{align*}
 & A=(1:0:0), && B=(0:1:0), && C=(0:0:1),
\end{align*}
\noindent $E_A$ (resp. $E_B,$ resp. $E_C$) le diviseur exceptionnel obtenu en \'eclatant
$A$ (resp. $B,$ resp. $C$) et~$\widetilde{L}_{AB}$ (resp. $\widetilde{L}_{AC},$ 
resp. $\widetilde{L}_{BC}$) la transform\'ee stricte de $L_{AB}$ (resp. 
$L_{AC},$ resp. $L_{BC}$).

\noindent Ceci n'est pas propre \`a $\sigma,$ c'est un fait g\'en\'eral: 
toute transformation de Cremona s'\'ecrit au moyen 
d'\'eclatements comme l'assure le th\'eor\`eme suivant d\^u \`a 
Zariski; ce r\'esultat est en fait valable pour toute
transformation birationnelle d'une surface projective lisse dans une 
autre. Avant de l'\'enoncer rappelons que si $\mathrm{X}$ est une vari\'et\'e 
irr\'eductible et $\mathrm{Y}$ une vari\'et\'e, une 
\textsl{transformation rationnelle}
$f\colon \mathrm{X}\dashrightarrow \mathrm{Y}$ est un 
morphisme d'un ouvert $\mathcal{U}$ de $\mathrm{X}$ dans $\mathrm{Y}$ qui n'est pas la
restriction d'un morphisme $\widetilde{\mathcal{U}}\to \mathrm{X}$ avec 
$\mathcal{U}\subsetneq\widetilde{\mathcal{U}}.$

\begin{thm}[Zariski, 1944]\label{Zariski}
{\sl Soient $\mathrm{S},$ $\widetilde{\mathrm{S}}$ deux surfaces projectives lisses et $f\colon
\mathrm{S}\dashrightarrow\widetilde{\mathrm{S}}$ une transformation birationnelle. Il existe une 
surface projective lisse $\mathrm{S}'$ et deux suites d'\'eclatements
\begin{align*}
& \pi_1\colon \mathrm{S}'\to \mathrm{S}, &&\pi_2\colon \mathrm{S}'\to \widetilde{\mathrm{S}}
\end{align*}
\noindent telles que 
$f=\pi_2\pi^{-1}_1$
\begin{align*}
\xymatrix{& \mathrm{S}'\ar[dl]_{\pi_1}\ar[dr]^{\pi_2} &\\
\mathrm{S}\ar@{-->}[rr]_f & & \widetilde{\mathrm{S}} }
\end{align*}}
\end{thm}

\noindent On va d\'emontrer le th\'eor\`eme \ref{Zariski} en suivant \cite{Be2}. 
La d\'emonstration se fait en deux temps
\begin{itemize}
\item tout d'abord on montre qu'une transformation rationnelle d'une surface
dans $\mathbb{P}^n(\mathbb{C})$ s'\'ecrit~$\phi\pi^{-1}$ o\`u $\pi$ d\'esigne une suite
d'\'eclatements et $\phi$ un morphisme
(th\'eor\`eme \ref{bebe1});

\item puis on \'etablit qu'un morphisme entre deux surfaces est la
compos\'ee d'un isomorphisme et d'une suite d'\'eclatements (th\'eor\`eme \ref{bebe2}).
\end{itemize}

\vspace{0.3cm} 

\noindent Si $\mathrm{V}$ d\'esigne une vari\'et\'e lisse on note $\mathrm{Pic}(\mathrm{V})$
le \textsl{groupe de Picard de $\mathrm{V},$} {\it 
i.e.} le groupe des classes d'isomorphisme de fibr\'es en droites sur $\mathrm{V}.$

\begin{rem}\label{syslin}
Soit $\mathrm{S}$ une surface. \`A une transformation rationnelle $f\colon
\mathrm{S}\dashrightarrow~\mathbb{P}^n(\mathbb{C})$ et un syst\`eme
d'hyperplans $H$ dans $\mathbb{P}^n(\mathbb{C})$ on peut associer le syst\`eme 
lin\'eaire $f^*\vert H\vert.$ En effet la transformation $f$ est un morphisme d'un ouvert maximal $U
$ de $\mathrm{S}$ dans $\mathbb{P}^n(\mathbb{C}).$ L'ensemble $\mathrm{S}\setminus U$ est fini donc 
$\mathrm{Pic}(\mathrm{S})$ et 
$\mathrm{Pic}(\mathrm{S}\setminus U)$ sont isomorphes; par cons\'equent on peut parler d'image inverse par
$f$ d'un syst\`eme lin\'eaire $\mathcal{P}$ que l'on note~$f^*\mathcal{P}.$ 

\noindent R\'eciproquement soit $\mathcal{P}$ un syst\`eme lin\'eaire de $\mathrm{S}$ sans
composante fixe; notons $\check{\mathcal{P}}$ l'espace projectif dual de $\mathcal{P}.$
On peut d\'efinir une transformation rationnelle $f\colon \mathrm{S}
\dashrightarrow\check{\mathcal{P}}:$ \`a un point $p$ de~$\mathrm{S}$ on
associe l'hyperplan de $\mathcal{P}$ constitu\'e des diviseurs passant par $p;$ notons
que $f$ est d\'efini en $p$ si et seulement si $p$ n'est pas un point base 
de~$\mathcal{P}.$

\noindent Il y a donc une bijection entre 
\smallskip
\begin{itemize}
\item les transformations rationnelles $f$ de $\mathrm{S}$
dans $\mathbb{P}^n(\mathbb{C})$ telles que $f(\mathrm{S})$ ne soit
pas contenu dans un hyperplan

\hspace{-1.3cm} et 

\item les syst\`emes lin\'eaires sur $\mathrm{S}$ sans composante fixe et de dimension $n.$
\end{itemize}
\end{rem}

\begin{thm}[\cite{Be2}]\label{bebe1}
{\sl Soit $f\colon \mathrm{S}\dashrightarrow\mathbb{P}^n(\mathbb{C})$
une transformation rationnelle d'une surface~$\mathrm{S}$ dans l'espace projectif $\mathbb{P}^n
(\mathbb{C}).$ Il existe une surface $\mathrm{S}',$ une suite d'\'eclatements $\pi
\colon \mathrm{S}'\to \mathrm{S}$ et un morphisme $\phi\colon
\mathrm{S}'\to\mathbb{P}^n(\mathbb{C})$ tels que $f=\phi\pi^{-1}$
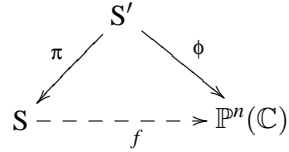
\begin{align*}
\xymatrix{&\mathrm{S}'\ar[dl]_{\pi}\ar[dr]^{\phi} &\\
\mathrm{S}\ar@{-->}[rr]_f & & \mathbb{P}^n(\mathbb{C}) }
\end{align*}}
\end{thm}

\begin{proof}[{\sl D\'emonstration}]
On peut se ramener au cas o\`u $f(\mathrm{S})$ n'est pas contenu dans un hyperplan 
de~$\mathbb{P}^n(\mathbb{C}).$ \`A $f$ correspond donc un syst\`eme lin\'eaire $\mathcal{P}
\subset\vert D\vert$ (si $D$ est un diviseur sur une surface~$\mathrm{S}$ on d\'esigne par $\vert D\vert$ l'ensemble
des diviseurs effectifs sur $S$ lin\'eairement \'equivalents \`a~$D$) de dimension $n$
sur $S$ sans composante fixe (Remarque \ref{syslin}). Si $\mathcal{P}$ n'a pas de point
 base, $f$ est un morphisme et le th\'eor\`eme est d\'emontr\'e.

\noindent Supposons que $\mathcal{P}$ ait un point base $p.$ Soit $\pi_1
\colon \mathrm{S}_1\to \mathrm{S}$ l'\'eclatement de $\mathrm{S}$ au point $p;$ la courbe
exceptionnelle $E_1$ appartient \`a une composante fixe de $\pi_1^*\mathcal{P}\subset\vert
\pi_1^*D\vert$ avec une certaine multiplici\-t\'e~$\nu_1\geq 1;$ autrement dit le syst\`eme 
\begin{align*}
\mathcal{P}_1=\pi_1^*\mathcal{P}-\nu_1E_1\subset\vert D_1\vert=\vert\pi_1^*D-\nu_1 E_1\vert
\end{align*}
\noindent n'a pas de composante fixe. \`A $\mathcal{P}_1$ correspond une transformation
rationnelle $f_1\colon \mathrm{S}_1\dashrightarrow\mathbb{P}^n
(\mathbb{C})$ qui co\"incide avec $f\pi_1.$ Si $f_1$ est un morphisme le th\'eor\`eme est
d\'emontr\'e; sinon on it\`ere le proc\'ed\'e. On obtient alors une suite d'\'eclatements
$(\pi_q\colon \mathrm{S}_q\to \mathrm{S}_{q-1})_q$ et une suite de
syst\`emes lin\'eai\-res~$(\mathcal{P}_q\subset\vert D_q\vert=\vert\pi_q^*D_{q-1}-\nu_qE_q\vert)_q$ sur
$\mathrm{S}_q$ sans composante fixe. \`A partir de (\ref{blow}) on obtient
 \begin{align*}
D_q^2=D_{q-1}^2-\nu_q^2<D_{q-1}^2.
\end{align*}
\noindent Comme $\mathcal{P}_q$ n'a pas de composante fixe on a $D_q^2\geq 0$ ce
qui assure que la suite d'\'eclatements est finie: le proc\'ed\'e conduit n\'ecessairement,
en un nombre fini d'\'etapes, \`a un syst\`eme lin\'eaire $\mathcal{P}_q$ sans point base et
d\'efinissant un morphisme $f\colon \mathrm{S}_q\to\mathbb{P}^n
(\mathbb{C}).$
\end{proof}

\begin{eg}
Soit $\mathrm{S}$ une surface quadrique lisse de $\mathbb{P}^3(\mathbb{C});$ une telle surface
est donn\'ee par
\begin{align*}
\sum q_{ij}x_ix_j=0
\end{align*}
\noindent o\`u $(q_{ij})$ d\'esigne une matrice sym\'etrique non d\'eg\'en\'er\'ee. Les formes quadratiques
sym\'etriques non d\'eg\'en\'er\'ees sur $\mathbb{C}^4$ sont isomorphes; il en 
r\'esulte que les surfaces quadriques lisses de~$\mathbb{P}^3(\mathbb{C})$ sont
projectivement isomorphes. Consid\'erons l'application de Segr\'e
\begin{align*}
&\xi\colon\mathbb{P}^1(\mathbb{C})\times\mathbb{P}^1
(\mathbb{C})\to\mathbb{P}^3(\mathbb{C}), && ((x:y),\,(u:v)) \mapsto(xu:xv:yu:
yv).
\end{align*}
\noindent L'application $\xi$ est un plongement dont l'image est la quadrique lisse
\begin{align*}
X_0X_3-X_1X_2=0.
\end{align*}
\noindent Ainsi toute surface quadrique $\mathrm{S}$ de $\mathbb{P}^3(\mathbb{C})$ est isomorphe \`a $
\mathbb{P}^1(\mathbb{C})\times\mathbb{P}^1(\mathbb{C}).$

\noindent Soit $p$ un point de $\mathrm{S};$ l'ensemble des droites de $\mathbb{P}^3(\mathbb{C})$ passant par $p$ s'identifie \`a
l'espace pro\-jectif~$\mathbb{P}^2(\mathbb{C}).$ \`A tout point $m$ de $S$ on peut
associer la droite passant par $m$ et $p;$ ceci induit une application rationnelle $f$ de
$\mathrm{S}$ dans $\mathbb{P}^2(\mathbb{C}),$ application rationnelle non d\'efinie en $p$ qui 
s'\'etend en un mor\-phisme~$\widetilde{f}\colon\widetilde{\mathrm{S}}\to
\mathbb{P}^2(\mathbb{C}).$ L'image par $\widetilde{f}^{-1}$ d'un point de $\mathbb{P}^2
(\mathbb{C}),$ correspondant \`a une droite~$L$ passant par $p,$ est constitu\'ee
\begin{itemize}
\item ou bien du second point d'intersection de $L$ avec $\mathrm{S};$

\item ou bien, si $L$ est contenue dans $\mathrm{S},$ de tous les points de $L.$
\end{itemize}
Ainsi $\widetilde{f}$ est un morphisme qui contracte les deux g\'en\'eratrices de $\mathrm{S}$
passant par $p.$ On a donc
\begin{align*}
\xymatrix{&\widetilde{\mathrm{S}}\ar[dl]_{\pi}\ar[dr]^{\widetilde{f}} &\\
\mathrm{S}\ar@{-->}[rr]_f & & \mathbb{P}^2(\mathbb{C}) }
\end{align*}
\noindent o\`u $\widetilde{f}$ est, \`a isomorphisme pr\`es, l'\'eclatement de deux points distincts et 
$\pi$ l'\'eclatement de $\mathrm{S}$ en~$p.$
\end{eg}

\begin{lem}[\cite{Be2}]\label{bobo0}
{\sl Soit $f\colon \mathrm{S}\to \mathrm{S}'$ un morphisme birationnel d'une surface irr\'eductible \'eventuellement singuli\`ere $\mathrm{S}$ dans une surface lisse $\mathrm{S}'.$
Supposons que la transformation rationnelle~$f^{-1}$ ne soit pas d\'efinie en un point $p$
de $\mathrm{S}';$ alors $f^{-1}(p)$ est une courbe sur~$\mathrm{S}.$}
\end{lem}

\noindent \`A partir de cet \'enonc\'e on obtient le r\'esultat qui suit.

\begin{lem}[\cite{Be2}]\label{bobo}
{\sl Soit $f\colon \mathrm{S}\dashrightarrow \mathrm{S}'$ une transformation
birationnelle entre deux surfaces~$\mathrm{S}$ et~$\mathrm{S}'.$ S'il existe un point $p$ de $\mathrm{S}'$ en lequel
$f^{-1}$ n'est pas d\'efinie il existe une courbe~$\mathcal{C}$ sur~$\mathrm{S}$ telle que
$f(\mathcal{C})=p.$}
\end{lem}

\begin{proof}[{\sl D\'emonstration}]
L'application $f$ induit un morphisme $\phi$ d'un ouvert $U$ de $\mathrm{S}$ dans $\mathrm{S}'.$ 
Notons~$\Gamma_\phi\subset U\times \mathrm{S}'$ le graphe de $\phi.$ L'adh\'erence $\overline{
\Gamma_\phi}$ de $\Gamma_\phi$ dans $\mathrm{S}\times \mathrm{S}'$ est une surface irr\'eductible,
\'eventuellement singuli\`ere. On d\'esigne par $\pi_1$ (resp. $\pi_2$) la projection de 
$\overline{\Gamma_\phi}$ sur $\mathrm{S}$ (resp. $\mathrm{S}'$); ces deux projections sont birationnelles
et $f=\pi_2\pi_1^{-1}$
\begin{align*}
\xymatrix{&\overline{\Gamma_\phi}\ar[dl]_{\pi_1}\ar[dr]^{\pi_2} &\\
\mathrm{S}\ar@{-->}[rr]_f & & \mathrm{S}' }
\end{align*}
\noindent Le fait que $f^{-1}$ ne soit pas d\'efini en $p$ implique que $\pi_2^{-1}$ aussi; le lemme
\ref{bobo0} assure l'existence d'une courbe $\widetilde{\mathcal{C}}$ sur 
$\overline{\Gamma_\phi}$ telle que $\pi_2(\widetilde{\mathcal{C}})=p.$ L'image de 
$\widetilde{\mathcal{C}}$ par $\pi_1$ est une courbe $\mathcal{C}$ satisfaisant~$f(\mathcal{C})=p.$
\end{proof}

\begin{pro}[\cite{La3}]\label{bobo2}
{\sl Soit $f\colon \mathrm{S}\to \mathrm{S}'$ un morphisme birationnel entre deux
surfaces $\mathrm{S}$ et~$\mathrm{S}'.$ Supposons que la transformation rationnelle $f^{-1}$ ne soit pas 
d\'efinie en un point $p$ de~$\mathrm{S}';$ alors~$f$ s'\'ecrit $\pi\phi$ o\`u $\pi\colon
\widetilde{\mathrm{S}}\to \mathrm{S}'$ d\'esigne l'\'eclatement de $\mathrm{S}'$ en $p$ et $\phi$ un
morphisme birationnel de~$\mathrm{S}$ dans~$\widetilde{\mathrm{S}}$
\begin{align*}
\xymatrix{&\widetilde{\mathrm{S}}\ar[dr]^{\pi} &\\
\mathrm{S}\ar[ru]^\phi\ar[rr]_f & & \mathrm{S}' }
\end{align*}}
\end{pro}

\begin{proof}[{\sl D\'emonstration}]
Soit $\pi\colon\widetilde{\mathrm{S}}\to \mathrm{S}'$ l'\'eclatement de $\mathrm{S}'$ au
point $p;$ posons $\phi:=\pi^{-1} f.$ Supposons que $\phi$ ne soit pas un morphisme;
soit $m$ un point de $\mathrm{S}$ o\`u $\phi$ n'est pas d\'efini. On a
\begin{itemize}
\item $f(m)=p;$

\item $f$ n'est pas localement inversible en $m;$

\item il existe une courbe sur $\widetilde{\mathrm{S}}$ contract\'ee sur $m$ par $\phi^{-1};$ cette
courbe ne peut \^etre que le diviseur exceptionnel $E$ associ\'e \`a $\pi.$ 
\end{itemize}

\noindent Consid\'erons $q_1,$ $q_2$ deux points distincts de $E$ o\`u $\phi^{-1}$ est
bien d\'efini et $\mathcal{C}_1,$ $\mathcal{C}_2$ deux germes de courbes lisses
transverses \`a $E$ en $q_1$ et $q_2$ respectivement. Alors $\pi(\mathcal{C}_1)$ et 
$\pi(\mathcal{C}_2)$ sont deux germes de courbes lisses transverses en $p=f(m)$ qui sont
image par $f$ de deux germes de courbes en $m$ 

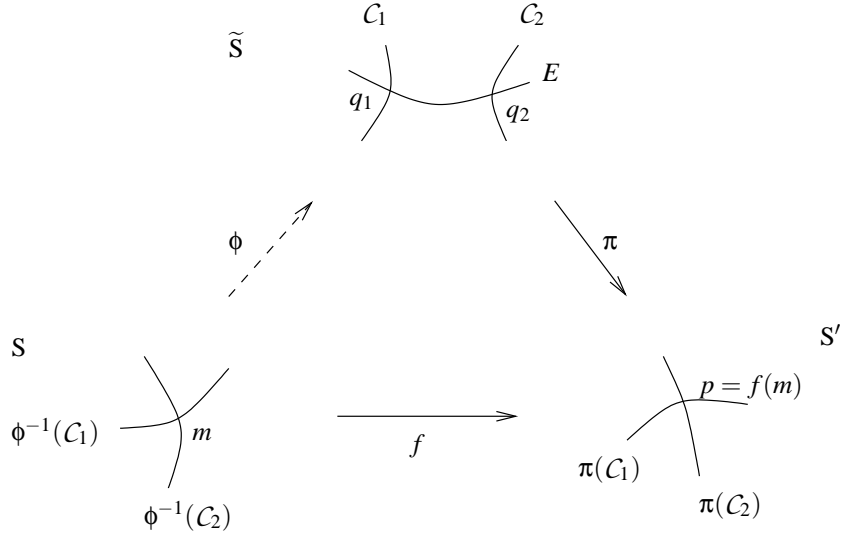
\begin{figure}[H]
\begin{center}
\input{zariski.pstex_t}
\end{center} 
\end{figure}

\noindent La diff\'erentielle de $f$ en $m$ est donc de rang $2$ ce qui contredit le fait que
 $f$ n'est pas localement inversible en $m.$
\end{proof}

\noindent L'\'enonc\'e qui pr\'ec\`ede nous permet de terminer la seconde \'etape.

\begin{thm}[\cite{Be2}]\label{bebe2}
{\sl Soit $f\colon \mathrm{S}\to \mathrm{S}_0$ un morphisme birationnel entre les
surfaces $\mathrm{S}$ et~$\mathrm{S}_0.$ Il existe une suite d'\'eclatements $(\pi_k\colon
\mathrm{S}_k\to \mathrm{S}_{k-1})_{k=1..n}$ et un isomorphisme $\phi\colon
\mathrm{S}\stackrel{\sim}{\longrightarrow}~\mathrm{S}_n$ tels que 
\begin{align*}
& f=\pi_1\ldots\pi_n\phi.
\end{align*}}
\end{thm}

\begin{proof}[{\sl D\'emonstration}]
Si $f$ est un isomorphisme la preuve est termin\'ee; sinon $f^{-1}$ n'est pas d\'efini en un
certain point $p$ de $\mathrm{S}_0.$ La proposition \ref{bobo2} assure que $f$ s'\'ecrit comme la
compos\'ee de l'\'eclatement $\pi_1$ de $\mathrm{S}_0$ en $p$ et du morphisme birationnel $f_1
\colon \mathrm{S}\to \mathrm{S}_1.$ Si $f_1$ n'est pas un isomorphisme on
it\`ere le proc\'ed\'e; reste \`a montrer qu'it\'erer un nombre fini de fois suffit. Raisonnons
par l'absurde: supposons qu'on ait une suite infinie d'\'eclatements $(\pi_k
\colon \mathrm{S}_k\to \mathrm{S}_{k-1})_k$ et de morphismes birationnels $(f_k
\colon \mathrm{S}\to~\mathrm{S}_k)_k$ tels que 
\begin{align*}
& f_{k-1}=\pi_k f_k, &&\forall\hspace{0.1cm} k\geq 1.
\end{align*}
\noindent Notons $\nu(f_k)$ le nombre de courbes irr\'eductibles contract\'ees par $f_k.$ Puisque
$f_{k-1}=~\pi_k f_k$ toute courbe contract\'ee par $f_k$ est donc contract\'ee par 
$f_{k-1}.$ De plus il existe au moins une courbe irr\'eductible $\mathcal{C}$ sur $\mathrm{S}_k$
telle que $f_k(\mathcal{C})$ soit le diviseur exceptionnel de $\pi_k$ d'o\`u $\nu(f_k)<
\nu(f_{k-1});$ pour $k$ suffisamment grand on a donc $\nu(f_k)<0$ ce qui est impossible. 
\end{proof}

\bigskip

\section{G\'en\'eration du groupe de Cremona}\label{generation}

\bigskip

\bigskip

\subsection{\'Enonc\'es de N\oe ther et Iskovskikh}

\bigskip

\noindent En $1871$ N\oe ther donne un r\'esultat de 
g\'en\'eration pour le groupe de Cremona.

\begin{thm}[N\oe ther, \cite{No, No2, No3}]\label{nono}
{\sl Le groupe de Cremona est engendr\'e par l'involution de 
Cremona $\sigma$ et $\mathrm{Aut}(\mathbb{P}^2(\mathbb{C})).$}
\end{thm}

\begin{egs}
\begin{itemize}
\item La d\'ecomposition de $\rho$ (\emph{voir} Exemples \ref{egsptdind})
est la suivante 
\begin{align*}
(z-y:y-x:y)\sigma(y+z:z:x)\sigma(x+z:y-z:z).
\end{align*}

\item La transformation $\tau$ (d\'efinie dans Exemples \ref{egsptdind})
s'\'ecrit aussi $\ell_1\sigma\ell_2\sigma \ell_3\sigma\ell_4\sigma\ell_5$
o\`u 

\begin{small}
\begin{align*}
& \ell_1=(y-x:2y-x:z-y+x), && \ell_2=(x+z:x:y), && \ell_3=(-y:x+z-3y:x), \\
& \ell_4=(x+z:x:y), && \ell_5=(y-x:-2x+z:2x-y). &&
\end{align*}
\end{small}

\item Un calcul montre que 
\begin{align*}
&\mathbb{P}^2(\mathbb{C})\dashrightarrow\mathbb{P}^2(\mathbb{C}),&& (x:y:z)
\mapsto(xz(x+y):yz(x+y):xy^2)
\end{align*}
\noindent se d\'ecompose comme suit
\begin{align*}
(2y-x:-y:-z)\sigma(y+z:z:x)\sigma(z:-2x-y:x+y)\sigma.
\end{align*}
\end{itemize}
\end{egs}

\noindent La d\'emonstration du
th\'eor\`eme \ref{nono} donn\'ee par N\oe ther \'etait incompl\`ete. 
Il ne consid\'erait que les transformations birationnelles quadratiques
ayant trois points d'ind\'etermination propres
alors que certaines en ont moins ({\it cf.} les transformations $\rho$
et $\tau$ de Exemples \ref{egsptdind}). Il faudra attendre $1901$ et 
Castelnuovo pour une preuve compl\`ete (\cite{Cas}). Au cours
du XX$^{\text{\`eme}}$ si\`ecle il y a eu de nombreuses d\'emonstrations
du th\'eor\`eme \ref{nono}, on renvoie \`a \cite{AC} pour un historique
d\'etaill\'e.

\noindent Dans les ann\'ees $80$ Gizatullin et 
Iskovskikh ont donn\'e une pr\'esentation de $\mathrm{Bir}
(\mathbb{P}^2)$ par g\'en\'erateurs et relations (\cite{Gi, Is});
\'enon\c{c}ons le r\'esultat d'Iskovskikh qu'il pr\'esente dans $\mathbb{P}^1
(\mathbb{C})\times\mathbb{P}^1(\mathbb{C})$ qui est birationnellement 
isomorphe \`a $\mathbb{P}^2(\mathbb{C}).$

\begin{thm}\label{isk}
{\sl Le groupe des transformations birationnelles de $\mathbb{P}^1(\mathbb{C})
\times\mathbb{P}^1(\mathbb{C})$ est engendr\'e par $\mathrm{Aut}(\mathbb{P}^1(
\mathbb{C})\times\mathbb{P}^1(\mathbb{C}))$ et le groupe $\mathrm{dJ}$ de de
Jonqui\`eres\footnote{\hspace{0.1cm} Le groupe de de Jonqui\`eres est birationnellement 
conjugu\'e au sous-groupe de $\mathrm{Bir}(\mathbb{P}^1(\mathbb{C})\times\mathbb{P}^1(\mathbb{C}))$ qui 
pr\'eserve la premi\`ere projection $p\colon\mathbb{P}^1(\mathbb{C})\times\mathbb{P}^1(\mathbb{C})
\to\mathbb{P}^1(\mathbb{C}).$}.

\noindent De plus les relations dans $\mathrm{Bir}(\mathbb{P}^1(\mathbb{C})\times\mathbb{P}^1(\mathbb{C}))$ sont 
les relations internes au groupe $\mathrm{dJ},$ au grou\-pe~$\mathrm{Aut}(\mathbb{P}^1(
\mathbb{C})\times\mathbb{P}^1(\mathbb{C}))$ auxquelles s'ajoutent la relation 
\begin{align*}
&(\eta e)^3=\left(\frac{1}{x},\frac{1}{y}\right)&& \text{o\`u} &&
\eta\colon(x,y)\mapsto(y,x) && \text{et} &&
e\colon(x,y)\mapsto\left(x,\frac{x} {y}\right).
\end{align*}}
\end{thm}

\medskip

\begin{rem} 
Il n'y a pas d'analogue au th\'eor\`eme
de N\oe ther en dimension sup\'erieure: il faut une infinit\'e non 
d\'enombrable de transformations de degr\'e strictement sup\'erieur \`a 
$1$ pour engendrer le groupe $\mathrm{Bir}(\mathbb{P}^n)$ des 
transformations birationnelles de $\mathbb{P}^n(\mathbb{C})$ dans 
lui-m\^eme (\emph{voir}~\cite{Pan}). L'id\'ee de la d\'emonstration est la suivante. Pan construit
pour toute cubique plane lisse une transformation birationnelle de 
$\mathbb{P}^n(\mathbb{C}),$ avec $n\geq 3,$ qui contracte
cette cubique. Or l'ensemble des classes d'isomorphismes des cubiques
planes lisses est une famille \`a un param\`etre; donc l'ensemble des types 
birationnels\footnote{\hspace{0.1cm} Soient $X$ et $X'$ deux sous-vari\'et\'es de 
$\mathbb{P}^n(\mathbb{C});$ on dit que $X$ et $X'$ ont m\^eme \textsl{type birationnel} si $X$ et $X'$ sont birationnellement \'equivalentes.} des composantes du lieu exceptionnel 
des \'el\'ements de 
$\mathrm{Bir}(\mathbb{P}^n)$ est infini.
Supposons que~$\mathrm{Bir}(\mathbb{P}^n)$ soit engendr\'e par un nombre 
fini de transformations $f_1,$ $\ldots,$ $f_r$ et $\mathrm{Aut}(\mathbb{P}^n
(\mathbb{C})).$ Les types birationnels des composantes du lieu exceptionnel de tout 
\'el\'ement de $\mathrm{Bir}(\mathbb{P}^n)$ s'obtiennent alors \`a partir de ceux des $f_i;$ 
l'ensemble des types birationnels des composantes du lieu exceptionnel des 
transformations de $\mathrm{Bir}(\mathbb{P}^n)$ est donc fini: contradiction.\smallskip

\noindent Par contre on dispose d'un \'enonc\'e semblable au th\'eor\`eme de N\oe ther
pour les transformations birationnelles r\'eelles qui sont des diff\'eomorphismes
de $\mathbb{P}^2(\mathbb{R})$ (\emph{voir} \cite{RV}); on a un \'enonc\'e comparable pour la sph\`ere $\mathbb{S}^2$ (\emph{voir} \cite{KM}).  
\end{rem}

\begin{rem}
Wright a reformul\'e le th\'eor\`eme de N\oe ther en terme de produit amalgam\'e: le 
groupe de Cremona est le produit amalgam\'e de  
$\mathrm{Aut}(\mathbb{P}^1(\mathbb{C})\times\mathbb{P}^1(\mathbb{C})),$
$\mathrm{dJ}$ et~$\mathrm{Aut}(\mathbb{P}^2(\mathbb{C}))$
suivant leur intersection (\cite{Wr2}).
\end{rem}

\bigskip

\subsection{Transformations birationnelles quadratiques}

\bigskip

\noindent Le th\'eor\`eme de N\oe ther conduit naturellement \`a s'int\'eresser
aux transformations birationnelles de degr\'e $2;$ notons $\mathrm{Bir}_2$
l'ensemble qu'elles forment. Le grou\-pe~$\mathrm{PGL}_3(\mathbb{C})\times\mathrm{PGL}_3
(\mathbb{C})$ agit sur $\mathrm{Bir}_2$
\begin{align*}
&\mathrm{PGL}_3(\mathbb{C})\times\mathrm{Bir}_2\times\mathrm{PGL}_3
(\mathbb{C})\to\mathrm{Bir}_2, && (A,f,B)\mapsto AfB^{-1}.
\end{align*}
\noindent L'orbite d'un \'el\'ement $f$ de $\mathrm{Bir}_2$ sous cette action est 
not\'ee $\mathcal{O}(f).$

\begin{thm}[\cite{CD}] \label{descr}
{\sl Si
\begin{align*}
& \Sigma^0=\mathcal{O}(x(x:y:z)), &&\Sigma^1=\mathcal{O}(\tau), 
&&\Sigma^2=\mathcal{O}(\rho), &&\Sigma^3=\mathcal{O}
(\sigma),
\end{align*}
\noindent alors
\begin{align*}
\mathrm{Bir}_2=\Sigma^0\cup\Sigma^1\cup\Sigma^2\cup\Sigma^3.
\end{align*}

\noindent On a 
\begin{align*}
& \dim\Sigma^0=10, && \dim\Sigma^1=12, &&\dim \Sigma^2=13, && \dim\Sigma^3=14
\end{align*}
\noindent ainsi que les conditions d'incidence suivantes
\begin{align*}
& \overline{\Sigma^0}=\Sigma^0, && \overline{\Sigma^1}=\Sigma^0\cup\Sigma^1, && 
\overline{\Sigma^2}=\Sigma^0\cup\Sigma^1\cup\Sigma^2, && \overline{\Sigma^3}=
\mathrm{Bir}_2=\Sigma^0\cup\Sigma^1\cup\Sigma^2\cup\Sigma^3.
\end{align*}}
\end{thm}

\begin{rem}
Les \'el\'ements de $\Sigma^i$ poss\`edent $i$ points d'ind\'etermination.
\end{rem}

\noindent D'apr\`es le th\'eor\`eme \ref{descr} il suffit de montrer que 
l'adh\'erence de $\Sigma^3$
est lisse le long de $\Sigma^0$ pour obtenir l'\'enonc\'e qui suit.

\begin{thm}[\cite{CD}] 
{\sl L'ensemble des transformations
birationnelles quadratiques est lisse dans l'ensemble des transformations rationnelles
quadratiques.}
\end{thm}

\begin{rem} 
La description, \`a composition \`a gauche et \`a droite pr\`es par un automorphisme de $\mathbb{P}^2(\mathbb{C}),$ 
des transformations birationnelles cubiques se fait aussi en \'etudiant la nature
du lieu des z\'eros du d\'eterminant jacobien; elle est sensiblement plus compliqu\'ee
($32$ mod\`eles). Alors que $\mathrm{Bir}_2$ est lisse et irr\'eductible, il n'en est 
pas de m\^eme pour l'ensemble des transformations birationnelles cubiques vu 
comme sous-ensemble de $\mathbb{P}^{29}(\mathbb{C})$ (le projectivis\'e de l'espace des polyn\^omes homog\`enes de degr\'e $3$ en $3$ variables s'identifie \`a $\mathbb{P}^{29}(\mathbb{C})$).
\end{rem}

\medskip

\noindent On peut se demander \`a quelle(s) condition(s) une transformation rationnelle
quadratique est birationnelle. Dans \cite{CD} deux crit\`eres sont \'etablis; l'un d'entre eux
est tr\`es simple \`a \'enoncer.

\begin{thm}[\cite{CD}]
{\sl Soit $Q$ une transformation rationnelle quadratique dont le d\'eterminant jacobien 
n'est pas identiquement nul. Supposons que 
$Q$ contracte le lieu des z\'eros $\mathcal{Z}$ du d\'eterminant jacobien de $Q;$
alors $\mathcal{Z}$ est l'union de trois
droites non concourantes et $Q$ est birationnelle.

\noindent De plus \`a composition \`a gauche et \`a droite pr\`es par un 
automorphisme de $\mathbb{P}^2(\mathbb{C})$ on a

\smallskip
\begin{itemize}

\item si $\mathcal{Z}$ est l'union de trois droites en position g\'en\'erale,
$Q$ est $\sigma;$

\item si $\mathcal{Z}$ est l'union d'une droite double et d'une droite
simple, $Q$ co\"incide $\rho;$

\item enfin si $\mathcal{Z}$ est une droite triple, $Q=\tau.$
\end{itemize}}
\end{thm}

\noindent Notons que ce crit\`ere ne se g\'en\'eralise pas en degr\'e 
sup\'erieur.

\begin{cor}[\cite{CD}]\label{corcrit}
{\sl Une transformation rationnelle quadratique du plan projectif dans
lui-m\^eme appartient \`a $\Sigma^3$ si et seulement si elle admet
trois points d'ind\'etermination.}
\end{cor}

\bigskip

\subsection{Transformations birationnelles et feuilletages}

\bigskip

\noindent Un \textsl{feuilletage} $\mathcal{F}$ de degr\'e $\nu$ sur 
$\mathbb{P}^2(\mathbb{C})$ est donn\'e par une $1$-forme
diff\'erentielle homog\`ene
\begin{align*}
\omega=F_0\mathrm{d}x+F_1\mathrm{d}y+F_2\mathrm{d}z,
\end{align*}
\noindent les $F_i$ d\'esignant des polyn\^omes homog\`enes de degr\'e
$\nu+1$ satisfaisant
\begin{align*}
& {\rm pgcd}(F_0,F_1,F_2)=1 && \text{et} && xF_0+yF_1+zF_2=0.
\end{align*}
\noindent Rappelons que le lieu singulier $\mathrm{Sing}(\mathcal{F})$ du
feuilletage $\mathcal{F}$ d\'efini par $\omega$ est le projectivis\'e
du lieu singulier de $\omega$
\begin{align*}
\mathrm{Sing}(\omega)=\{F_0=F_1=F_2=0\};
\end{align*}
\noindent l'identit\'e d'Euler $xF_0+yF_1+zF_2=0$ assure que cet
ensemble n'est jamais vide. Plus pr\'ecis\'ement on a une formule
de type Bezout
\begin{align*}
\#\,\mathrm{Sing}(\mathcal{F})=\nu^2+\nu+1,
\end{align*}
\noindent chaque point singulier \'etant compt\'e avec multiplicit\'e. \`A une 
transformation rationnel\-le~$f=(f_0:f_1:f_2)$ du plan projectif
on peut associer un feuilletage $\mathcal{F}(f)$ de 
$\mathbb{P}^2(\mathbb{C})$ d\'ecrit par la $1$-forme
\begin{align*}
(yf_2-zf_1)\mathrm{d}x+(zf_0-xf_2)\mathrm{d}y+(xf_1-yf_0)\mathrm{d}z.
\end{align*} 

\begin{eg}
Le feuilletage $\mathcal{F}(\sigma)$ associ\'e \`a $\sigma$ est
d\'efini par 
\begin{small}
\begin{align*}
& x(z^2-y^2)\mathrm{d}x+y(x^2-z^2)\mathrm{d}y+z(y^2-x^2)\mathrm{d}z=\frac{1}{2}(x^2,y^2,z^2)^*((z-y)\mathrm{d}x+(x-z)\mathrm{d}y+
(y-x)\mathrm{d}z).
\end{align*}
\end{small}
\noindent Le feuilletage $\mathcal{F}(\sigma)$ est de degr\'e $2,$ admet pour
int\'egrale premi\`ere $\frac{x^2-z^2}{y^2-z^2}$ et poss\`ede sept 
points singuliers qui sont les points fixes et les points d'ind\'etermination
de $\sigma.$ On a donc l'\'enonc\'e suivant.
\end{eg}

\begin{pro}[\cite{CD}]
{\sl Si $f$ est une transformation birationnelle quadratique g\'en\'erique,~$\mathcal{F}(f)$ poss\`ede sept points singuliers; quatre sont des points
fixes de $f,$ les trois autres des points d'ind\'etermination.}
\end{pro}

\noindent L'ensemble des applications rationnelles quadratiques 
s'identifie \`a $\mathbb{P}^{17}(\mathbb{C})$ et celui des feuilletages
de degr\'e $2$ sur $\mathbb{P}^2(\mathbb{C})$ \`a $\mathbb{P}^{14}
(\mathbb{C});$ l'application $f\mapsto \mathcal{F}(.)$ induit 
une application lin\'eaire qu'on note encore $\mathcal{F}(.)$ de 
$\mathbb{P}^{17}(\mathbb{C})$ dans $\mathbb{P}^{14}(\mathbb{C}).$

\noindent Soit $\mathcal{F}$ un feuilletage quadratique g\'en\'erique
sur le plan projectif complexe, {\it i.e.} dont les sept points singuliers
sont en position g\'en\'erale. Le th\'eor\`eme de division de de
Rham-Saito (\cite{Sa}) assure l'existence de 
polyn\^omes homog\`enes $f_0,$ $f_1$ et $f_2$ de degr\'e $2$
tels que  $\mathcal{F}$ soit d\'efini par
\begin{align*}
(yf_2-zf_1)\mathrm{d}x+(zf_0-xf_2)\mathrm{d}y+(xf_1-yf_0)\mathrm{d}z.
\end{align*}
\noindent Notons $f$ la transformation quadratique dont les composantes sont $f_0,$
$f_1$ et $f_2.$ L'ensem\-ble~$\mathrm{Sing}(\mathcal{F})$ est le lieu des
points $p_i$ pour lesquels il existe $\eta_i$ dans $\mathbb{C}$ tel
que $f(p_i)=\eta_ip_i.$ Choisissons trois points $p_1,$ $p_2$ et $p_3$
dans $\mathrm{Sing}(\mathcal{F}).$ Les $p_i$ n'\'etant pas align\'es
il existe une forme lin\'eaire~$\ell$ telle que~$\ell(p_i)=-\eta_i.$ 
L'application $\phi=f+\ell.\mathrm{id}$ satisfait $\mathcal{F}(f)=
\mathcal{F}(\phi)=\mathcal{F}$ et les points $p_i$ sont d'ind\'etermination
pour $f.$ En particulier $f$ est birationnelle (Corollaire \ref{corcrit}). La restriction 
de~$\mathcal{F}(.)$ \`a~$\mathrm{Bir}_2$ est donc dominante et par suite \`a fibre 
g\'en\'erique finie puisque la dimension de $\mathrm{Bir}_2$ 
co\"incide avec celle de l'ensemble des feuilletages quadratiques 
sur $\mathbb{P}^2(\mathbb{C}).$

\noindent Un feuilletage quadratique g\'en\'erique est d\'etermin\'e par la 
position de ses sept points singuliers et toute configuration g\'en\'erique est 
r\'ealis\'ee (\cite{CO, GMK}); on h\'erite d'une propri\'et\'e analogue pour les \'el\'ements
de~$\Sigma^3.$

\begin{pro}[\cite{CD}]
{\sl Une transformation birationnelle quadratique g\'en\'erique est d\'etermin\'ee
par la position de ses points d'ind\'etermination et de ses points fixes.}
\end{pro}

\noindent On en d\'eduit que la fibre g\'en\'erique de la restriction de 
l'application $\mathcal{F}(.)$ \`a $\mathrm{Bir}_2$ compte $35$ 
points exactement (choix de $3$ points d'ind\'etermination parmi $7$ points).

\bigskip

\section{Sous-groupes finis du groupe de Cremona}\label{gpefini}

\bigskip

\noindent L'\'etude des sous-groupes finis du groupe de Cremona a commenc\'e dans les ann\'ees $1870$ avec Bertini, Kantor et Wiman (\cite{Ber, K, Wi}). De nombreux auteurs
se sont depuis int\'eress\'es \`a ce sujet; citons par exemple 
\cite{BaBe, Be, BeBl, Bl, dF, DI}. En $2006$ Dolgachev et Iskovskikh ont am\'elior\'e
les r\'esultats de Kantor et Wiman et ont donn\'e une liste des sous-groupes finis de $\mathrm{Bir}(\mathbb{P}^2);$
avant de donner un th\'eor\`eme cl\'e sur lequel ils s'appuient, introduisons 
quelques notions. Commen\c{c}ons par les \textsl{surfaces $\mathrm{F}_n$ de Hirzebruch}. 
Posons $\mathrm{F}_0=\mathbb{P}^1(\mathbb{C})\times
\mathbb{P}^1(\mathbb{C});$ la surface $\mathrm{F}_1$ est la surface obtenue en
\'eclatant $\mathbb{P}^2(\mathbb{C})$ en $(1:0:0);$ cette surface est un
compactifi\'e de $\mathbb{C}^2$ naturellement muni d'une fibration
rationnelle correspondant aux droites $y=$ cte. Le diviseur \`a l'infini est constitu\'e
de deux courbes rationnelles s'intersectant transversalement en un point. On a 
\begin{itemize}
\item d'une part la transform\'ee stricte de la droite \`a l'infini dans $\mathbb{P}^2(
\mathbb{C}),$ c'est une fibre;

\item d'autre part le diviseur exceptionnel de l'\'eclatement qui est une section pour
la fibration.
\end{itemize}
Plus g\'en\'eralement $\mathrm{F}_n$ est, pour tout $n\geq
1,$ un compactifi\'e de $\mathbb{C}^2$ muni d'une fibration rationnelle telle que le
diviseur \`a l'infini soit constitu\'e de deux courbes rationnelles transverses, une
fibre~$f$ et une section $s_n$ d'auto-intersection $-n$
\begin{align*}
\mathrm{F}_n=\mathbb{P}_{\mathbb{P}^1(\mathbb{C})}(\mathcal{O}_{\mathbb{P}^1
(\mathbb{C})}\oplus\mathcal{O}_{\mathbb{P}^1(\mathbb{C})}(n)), && \forall
\hspace{0.1cm} n\geq 2.
\end{align*}

\noindent Pla\c{c}ons-nous sur 
$\mathrm{F}_n;$ on note $p$ l'intersection de $s_n$ et $f,$ $\pi_1$ 
l'\'eclatement de $\mathrm{F}_n$ en $p$ et $\pi_2$ la contraction de la
transform\'ee stricte $\widetilde{f}$ de $f.$ On passe de $\mathrm{F}_n$ 
\`a $\mathrm{F}_{n+1}$ via $\pi_2\pi_1^{-1}$

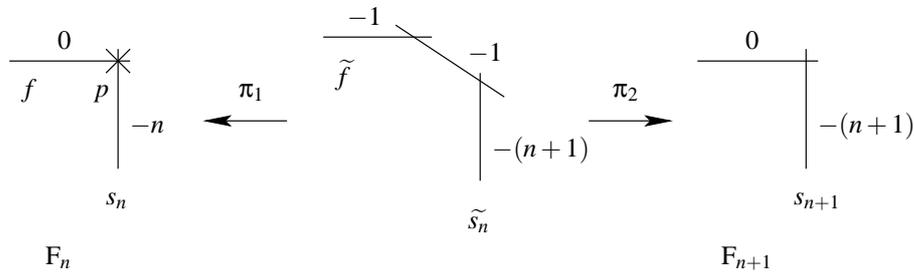
\begin{figure}[H]
\begin{center}
\input{el1.pstex_t}
\end{center}
\end{figure}

\noindent On peut aussi passer de la $(n+1)$-i\`eme \`a la $n$-i\`eme surface de
Hirzebruch via $\pi_2\pi_1^{-1}$ o\`u $\pi_1$ est l'\'eclatement de $\mathrm{F}_{n+1}$ 
en un point $p$ de la fibre $f$ qui n'appartient pas \`a $s_{n+1}$ et 
$\pi_2$ la contraction de la transform\'ee stricte $\widetilde{f}$ de $f$

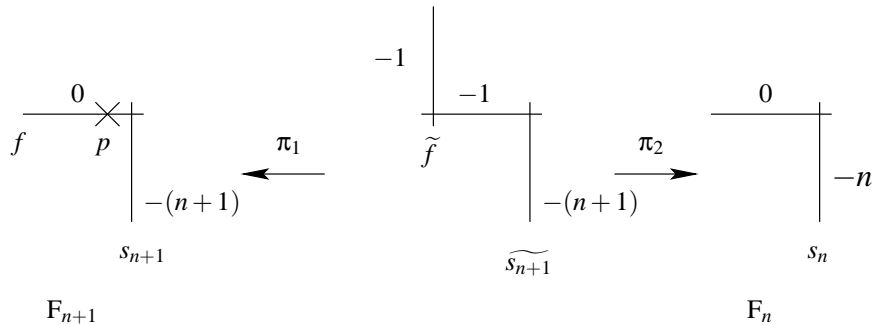
\begin{figure}[H]
\begin{center}
\input{el2.pstex_t}
\end{center}
\end{figure}

\noindent Soit $\mathrm{S}$ une surface projective lisse; un
\textsl{fibr\'e en coniques} $\eta\colon \mathrm{S}\to\mathbb{P}^1(
\mathbb{C})$ est un morphisme dont les fibres g\'en\'erales sont 
de genre $0$ et les fibres singuli\`eres l'union de deux droites.
Une surface munie d'un fibr\'e en coniques est isomorphe
\begin{itemize}
\item ou bien \`a $\mathrm{F}_n;$ 

\item ou bien \`a $\mathrm{F}_n$ \'eclat\'e en un nombre fini de points, 
tous appartenant \`a des fibres diff\'erentes. Le nombre d'\'eclatements 
co\"{i}ncide avec le nombre de fibres singuli\`eres du fibr\'e.
\end{itemize}

\noindent Une \textsl{surface de Del Pezzo} est une surface lisse $\mathrm{S}$ dont le 
fibr\'e anticanonique $-K_\mathrm{S}$ est ample. Le degr\'e d'une surface de Del Pezzo est 
$d=K_\mathrm{S}^2.$ D'apr\`es la formule de N\oe ther on a $1\leq d\leq 9.$ 

\noindent Si $d=9,$ alors $\mathrm{S}\simeq\mathbb{P}^2(\mathbb{C});$ si $d=8$ alors 
$\mathrm{S}\simeq\mathbb{P}^1(\mathbb{C})\times \mathbb{P}^1(\mathbb{C})$ ou $S\simeq
\mathrm{F}_1.$ Pour $d\leq 7$ la surface $S$ est isomorphe \`a $\mathbb{P}^2(\mathbb{C})$
\'eclat\'e en $n=9-d$ points $p_i,$ ces points satisfaisant les conditions suivantes
\begin{itemize}
\item il n'y en a pas $3$ align\'es;

\item il n'y en a pas $6$ sur une conique;

\item si de plus $n=8$ les $p_i$ ne sont pas sur une cubique plane dont le point singulier
serait un~des~$p_i.$
\end{itemize}

\begin{thm}[\cite{Ma, Is2}]\label{ma} {\sl Soit $\mathrm{G}$ un sous-groupe
fini du groupe de Cremona. Il existe une surface projective lisse $\mathrm{S}$ et une 
transformation birationnelle $\phi\colon \mathbb{P}^2(\mathbb{C})\dashrightarrow S$ telles que
$\phi\mathrm{G}\phi^{-1}$ soit un sous-groupe de~$\mathrm{Aut}
(\mathrm{S}).$ De plus on peut supposer que
\smallskip 
\begin{itemize}
\item ou bien $\mathrm{S}$ est une surface de Del Pezzo;

\item ou bien il existe un fibr\'e en coniques $\mathrm{S}\to\mathbb{P}^1(\mathbb{C})$ 
invariant par $\phi\mathrm{G}\phi^{-1}.$
\end{itemize}}
\end{thm}

\noindent Notons que l'alternative qui pr\'ec\`ede n'est pas exclusive: il 
y a des fibr\'es en coniques sur les surfaces de Del Pezzo.

\noindent Dolgachev et Iskovskikh donnent une caract\'erisation 
des couples $(\mathrm{G},\mathrm{S})$ satisfaisant l'une des \'eventualit\'es du 
th\'eor\`eme \ref{ma}. Ils utilisent ensuite la th\'eorie de Mori afin de 
pouvoir d\'eterminer quand deux couples sont 
birationnellement conjugu\'es. Notons que si le premier point avait \'et\'e 
partiellement r\'esolu par Wiman et Kantor, le second non. 
Il reste encore des questions ouver\-tes~(\cite{DI}, \S 9), comme par exemple d\'ecrire
les vari\'et\'es alg\'ebriques qui param\`etrent les classes de conjugaisons des
sous-groupes finis de $\mathrm{Bir}(\mathbb{P}^2).$
Blanc donne une r\'eponse \`a cette question pour
les sous-groupes ab\'eliens finis de $\mathrm{Bir}(\mathbb{P}^2)$ 
dont aucun \'el\'ement ne fixe une courbe de genre 
positif et aussi pour les sous-groupes cycliques
d'ordre fini et les \'el\'ements d'ordre fini de $\mathrm{Bir}(\mathbb{P}^2)$ 
(\emph{voir}~\cite{Bl, Bl4}).

\bigskip

\noindent Avant d'\'enoncer le r\'esultat de Bertini, repris par Bayle
et Beauville, introduisons trois types
d'involutions qui, comme on le verra, jouent un r\^ole tr\`es particulier.

\noindent Soient 
$p_1,$ $\ldots,$ $p_7$ sept points de $\mathbb{P}^2(\mathbb{C})$ en position 
g\'en\'erale. D\'esignons par $L$ le syst\`eme lin\'eaire de cubiques passant
par les $p_i;$ il est de dimension $2.$ Soit $p$ un point g\'en\'erique
de $\mathbb{P}^2(\mathbb{C});$ consid\'erons le pinceau $L_p$
constitu\'e des \'el\'ements de $L$ passant
par $p.$ Un pinceau de cubiques g\'en\'eriques ayant neuf points 
base, on d\'efinit par $i_G(p)$ le neuvi\`eme point base de $L_p.$ 
L'involution $i_G$ qui \`a~$p$ associe $i_G(p)$ ainsi construite est 
appel\'ee \textsl{involution de Geiser}. On peut v\'erifier qu'une telle 
involution est birationnelle de degr\'e $8;$ ses  
points fixes forment une courbe non hyperelliptique de genre~$3,$
de degr\'e $6$ avec des points double aux sept points choisis.

\noindent Soient $p_1,$ $\ldots,$ $p_8$ huit points de $\mathbb{P}^2(\mathbb{C})$ en 
position g\'en\'erale. Consid\'erons le pinceau de cubiques passant
par ces huit points; il a un neuvi\`eme point base qu'on notera
$p_9.$ Soit $p$ un point g\'en\'erique de $\mathbb{P}^2(\mathbb{C});$ il y a une 
unique cubique $\mathcal{C}(p)$ qui passe par les $p_i$ et $p.$ Sur $\mathcal{C}
(p)$ il y a une loi de groupe avec $p_9$ comme \'el\'ement neutre. On note $i_B$ l'involution qui \`a $p$ associe $-p;$ alors 
$i_B$ d\'efinit une involution birationnelle de $\mathbb{P}^2(\mathbb{C})$ de degr\'e~$17$ appel\'ee \textsl{involution de Bertini}. Les points fixes de $i_B$
forment une courbe non hyperelliptique de genre~$4,$ de degr\'e $9$
avec des points triples aux huit points choisis.

\noindent Pour finir introduisons les involutions de de Jonqui\`eres. Soit 
$\mathcal{C} $ une courbe de degr\'e $\nu\geq 2$ et soit~$p$ un point 
sur $\mathcal{C}$ de multiplicit\'e $\nu-2$ (si $\nu=2,$ le point $p$
n'appartient pas \`a $\mathcal{C}$). On suppose que $p$ est l'unique
point singulier de $\mathcal{C}.$ Au couple 
$(\mathcal{C},p)$ on va associer l'unique involution 
birationnelle~$i_{\mathrm{dJ}}$ qui fixe la courbe $\mathcal{C}$ et pr\'eserve les
droites passant par $p.$ Soit~$m$ un point g\'en\'erique de~$\mathbb{P}^2(
\mathbb{C});$ notons~$q_m,$ $r_m$ les deux points d'intersection, distincts de
$p,$ de la droite~$(pm)$ avec $\mathcal{C}.$ La transformation $i_{\mathrm{dJ}}$ associe au 
point $m$ le conjugu\'e harmonique de $m$ sur la droite $(pm)$
par rapport \`a~$q_m$ et~$r_m;$ autrement dit le point $i_{\mathrm{dJ}}(m)$ v\'erifie
la propri\'et\'e suivante: le birapport de $m,$ $i_{\mathrm{dJ}}(m),$~$q_m$ et~$r_m$ vaut~$-1.$ 
La transformation $i_{\mathrm{dJ}}$ est une \textsl{involution de de Jonqui\`eres}
de degr\'e $\nu$ centr\'ee en~$p$ et pr\'eservant~$\mathcal{C}.$ La normalis\'ee de la 
courbe de points fixes de $i_{\mathrm{dJ}}$ est une courbe hyperelliptique de 
genre $\nu-2$ d\`es que~$\nu\geq 3$ (par convention on dira qu'une courbe 
elliptique est hyperelliptique). 

\noindent L'\'enonc\'e suivant donne la classification des involutions birationnelles.

\begin{thm}[\cite{Ber, BaBe}]
{\sl Une involution de Cremona est \`a conjugaison birationnelle pr\`es de l'un 
des types suivants
\smallskip
\begin{itemize}
\item une involution projective;

\item une involution de de Jonqui\`eres $i_{\mathrm{dJ}}$ de degr\'e $\nu\geq 2;$

\item une involution de Bertini $i_B;$

\item une involution de Geiser $i_G.$
\end{itemize}}
\end{thm}

\noindent Dans \cite{BaBe} les auteurs montrent de plus que les classes 
de conjugaison des sous-groupes finis cycliques 
d'or\-dre~$2$ de $\mathrm{Bir}(\mathbb{P}^2)$ sont 
uniquement d\'etermin\'ees par le type birationnel des courbes de points 
fixes de genre positif. Plus pr\'ecis\'ement l'ensemble des classes
de conjugaison des sous-groupes d'ordre~$2$ de $\mathrm{Bir}
(\mathbb{P}^2)$ est une vari\'et\'e alg\'ebrique non 
connexe dont les composants connexes sont respectivement 
isomorphes 
\begin{itemize}
\item \`a l'espace des modules des courbes hyperelliptiques
de genre $g$ (involutions de de Jonqui\`eres);

\item \`a 
l'espace des modules des courbes canoniques de genre 
$3$ (involutions de Geiser);

\item \`a l'espace des modules des courbes canoniques de
genre $4$ avec une theta caract\'eristique nulle, isomorphes
\`a une intersection non singuli\`ere d'une surface
cubique et d'un c\^one quadratique dans $\mathbb{P}^3$
(involutions de Bertini).
\end{itemize}

\noindent Il n'est donc pas surprenant de constater que simultan\'ement 
\`a l'\'etude des sous-groupes finis de~$\mathrm{Bir}(\mathbb{P}^2)$ se 
d\'eveloppe la caract\'erisation des
transformations birationnelles qui fixent une courbe de genre 
donn\'e. Soit $\mathcal{C}$ une courbe irr\'eductible dans $\mathbb{P}^2
(\mathbb{C});$ le \textsl{groupe d'inertie} de $\mathcal{C},$ not\'e~$\mathrm{Ine}(\mathcal{C}),$ est le sous-groupe
des \'el\'ements de~$\mathrm{Bir}(\mathbb{P}^2)$ qui fixent $\mathcal{C}$ point par point.
Soit~$\mathcal{C}$ une courbe de~$\mathbb{P}^2(\mathbb{C})$ de genre $>1;$ 
dans les ann\'ees $1890$ Castelnuovo a montr\'e qu'un \'el\'ement de 
$\mathrm{Ine}(\mathcal{C})$ est ou bien une transformation de de Jonqui\`eres, 
ou bien une transformation de Cremona d'ordre $2,$ $3$ ou $4$ (\emph{voir} 
\cite{C}). Blanc, Pan et Vust pr\'ecisent l'\'enonc\'e de Castelnuovo.

\begin{thm}[\cite{BPV}]
{\sl Soit $\mathcal{C}$ une courbe irr\'eductible de $\mathbb{P}^2(\mathbb{C})$
de genre strictement sup\'erieur \`a $1.$ Un \'el\'ement $f$ de $\mathrm{Ine}
(\mathcal{C})$ est ou bien conjugu\'e \`a une transformation de de Jonqui\`eres,
ou bien d'ordre $2$ ou $3.$ Dans le premier cas si $f$ est d'ordre fini,
c'est une involution.}
\end{thm}

\noindent Pour le d\'emontrer ils suivent la m\^eme id\'ee que Castelnuovo; ils
construisent le syst\`eme lin\'eaire adjoint de $\mathcal{C}:$ soient $\pi\colon
Y\to\mathbb{P}^2(\mathbb{C})$ une r\'esolution plong\'ee des singularit\'es de $\mathcal{C}$
et~$\widetilde{\mathcal{C}}$ la transform\'ee stricte de $\mathcal{C}$ par 
$\pi^{-1}.$ Si le syst\`eme lin\'eaire $\vert\widetilde{C}+K_Y\vert$ n'est ni vide, ni r\'eduit
\`a un diviseur, $\pi_*\vert\widetilde{C}+K_Y\vert$ priv\'e de ses
\'eventuelles composantes fixes est le syst\`eme lin\'eaire adjoint. En it\'erant
cette construction, appel\'ee m\'ethode des syst\`emes lin\'eaires adjoints successifs, ils 
obtiennent que tout \'el\'ement $f$ de $\mathrm{Ine}(\mathcal{C})$ pr\'eserve 
une fibration $\mathcal{F},$ rationnelle ou elliptique. Si $\mathcal{F}$
est rationnelle,~$f$ est de de
Jonqui\`eres. Supposons que $\mathcal{F}$ soit elliptique. Puisque $g(
\mathcal{C})>1$ la restriction de~$f$ \`a une fibre g\'en\'erique est un automorphisme
avec au moins deux points fixes;  par suite $f$ est d'ordre $2,$ $3$ ou $4.$ 
En appliquant certains r\'esultats classiques sur les automorphismes des courbes
elliptiques, ils montrent que $f$ est d'ordre~$2$ ou~$3.$ 
 Finalement ils remarquent que ce r\'esultat ne se g\'en\'eralise
pas aux courbes de genre inf\'erieur ou \'egal \`a $1.$ Le cas des 
courbes de genre $1$ a \'et\'e trait\'e avec des techniques diff\'erentes dans~\cite{Pan2} 
et~\cite{Bl2}.

\begin{rem}
Signalons que Diller, Jackson et 
Sommese classifient les courbes invariantes par une
transformation birationnelle d'une surface projective complexe (\cite{DJS}).
Leur approche utilise des techniques de dynamique complexe; n'ayant pas encore
introduit le vocabulaire ad\'equat on abordera leurs r\'esultats au \S \ref{dyn}.
\end{rem}

\noindent Le nombre de classes de conjugaison
des transformations de Cremona d'ordre $2$ dans $\mathrm{Bir}
(\mathbb{P}^2)$ est infini (\cite{BaBe}). Montrons qu'il en est de m\^eme pour 
les \'el\'ements d'ordre $3$ et $5.$ 
Si deux transformations de Cremona $f$ et $g$ sont 
conjugu\'ees via $\phi$ alors $\phi$ envoie les courbes
non rationnelles fix\'ees par~$f$ sur les courbes non rationnelles
fix\'ees par $g.$ Soient~$f$ et $g$ deux transformations de  
Cremona d'ordre $3$ (resp.~$5$); les courbes de points
fixes de celles-ci peuvent \^etre n'importe quelle courbe
elliptique.  Comme le nombre de classes d'isomorphisme
de telles courbes est infini, le nombre de classes de conjugaison
dans~$\mathrm{Bir}(\mathbb{P}^2)$ des \'el\'ements de Cremona
d'ordre $3$ (resp.~$5$) aussi. 
\'Etant donn\'e un entier positif~$n$ on peut se demander combien 
vaut le nombre $\nu(n)$ de classes 
de conjugaison d'une transformation de Cremona d'ordre~$n$
dans $\mathrm{Bir}(\mathbb{P}^2).$ 
De Fernex r\'epond \`a cette question pour $n$ premier (\cite{dF});
on trouve dans \cite{Bl3} une r\'eponse pour tout~$n.$

\begin{thm}[\cite{Bl3}]\label{number1}
{\sl \textbf{\textit{1.}} Pour $n$ pair $\nu(n)$ est infini; il en est de m\^eme pour $n=3$
ou $n=5.$

\noindent\textbf{\textit{2.}} Soit $n$ un entier impair distinct de $3$ et $5;$ alors $\nu(n)$ est fini. 
De plus 
\begin{itemize}
\item $\nu(9)=3;$

\item $\nu(15)=9;$

\item $\nu(n)=1$ sinon.
\end{itemize}}
\end{thm}

\noindent Commen\c{c}ons par donner une id\'ee de la 
d\'emonstration de la premi\`ere assertion pour $n$
distinct de $2,$ $3$ et $5,$ \'eventualit\'es pr\'ec\'edemment
abord\'ees. 
Supposons $n\geq 2$ et consid\'erons un \'el\'ement~$P$ 
de~$\mathbb{C}[x^n]$ sans racine multiple. Blanc montre
 qu'il existe une transformation de Cremona~$f$
d'ordre~$2n$ telle que $f^n$ soit l'involution $(x,y)\mapsto(x,P(x)/y)$ 
qui fixe la courbe hyperellipti\-que~$y^2=P(x);$ le cas $n=2$ 
permet alors de conclure pour $n$ pair sup\'erieur ou \'egal \`a $4.$

\medskip

\noindent Pour d\'emontrer la seconde partie de l'\'enonc\'e Blanc
applique le th\'eor\`eme \ref{ma} au groupe fini engendr\'e par une 
transformation de Cremona d'ordre impair sup\'erieur ou \'egal \`a $7.$

\bigskip

\section{Groupe des automorphismes du groupe de Cremona}\label{autbir}

\bigskip

\noindent Il est classique d'\'etudier les automorphismes d'un groupe
$\mathrm{G}$ donn\'e. On peut, par exemple, d\'ecrire le
groupe d'automorphismes de $\mathrm{PGL}_3(\mathbb{C}).$

\begin{thm}[\cite{Die}]\label{dieudon}
{\sl Un automorphisme de $\mathrm{PGL}_3(\mathbb{C})$ s'obtient \`a partir
\smallskip
\begin{itemize}
\item des automorphismes int\'erieurs;

\item de la contragr\'ediente (l'involution $u\mapsto\transp u^{-1}$);

\item de l'action des automorphismes du corps $\mathbb{C}.$
\end{itemize}}
\end{thm}

\noindent Ce th\'eor\`eme est qualifi\'e de th\'eor\`eme g\'eom\'etrique. 
Consid\'erons un espace vectoriel $F$ sur un corps $\Bbbk$ et $\tau$ un 
automorphisme du corps $\Bbbk.$ Une \textsl{collin\'eation} (relative \`a 
$\tau$) est une application semi-lin\'eaire bijective $\ell\colon F\to F$ satisfaisant
\begin{align*}
&\ell(x+y)=\ell(x)+\ell(y), && \ell(x\alpha)=\ell(x)\tau(\alpha), &&\forall
\, x,\, y\in F,\hspace{1cm} \forall\,\alpha\in\Bbbk.
\end{align*} 
\noindent Une collin\'eation donne, par passage au quotient, une application bijective
$\overline{\ell}$ de $\mathbb{P}(F)$ dans lui-m\^eme; les
collin\'eations sont pr\'ecis\'ement les transformations bijectives respectant
les relations d'incidence. Le th\'eor\`eme fondamental de la g\'eom\'etrie 
projective s'\'enonce comme suit. 

\begin{thm}[\cite{Die}]
{\sl Soient $F$ un espace vectoriel de dimension $n$ sur un corps $\Bbbk$
et $f$ une application bijective de $\mathbb{P}(F)$ dans lui-m\^eme. 
Supposons que $f$ pr\'eserve l'alignement. Si $n\geq 3$ il existe une 
collin\'eation relative \`a un automorphisme du corps $\Bbbk$ telle que 
$f=\overline{\ell}.$}
\end{thm}

\noindent Dans le m\^eme esprit on peut d\'ecrire le groupe des automorphismes
du groupe des automorphismes de la droite affine complexe
\begin{align*}
& \mathrm{G}=\{\alpha z+\beta\hspace{0.1cm}\vert\hspace{0.1cm}\alpha\in\mathbb{C}^*,\hspace{0.1cm}\beta\in\mathbb{C}\}.
\end{align*}

\begin{thm} 
{\sl Tout automorphisme de $\mathrm{G}$ s'obtient \`a partir des automorphismes
int\'erieurs et de l'action d'un automorphisme du corps $\mathbb{C}.$} 
\end{thm}

\noindent En effet soit $\phi$ un \'el\'ement de $\mathrm{Aut}(\mathrm{G}).$ L'image d'un sous-groupe
ab\'elien maximal de $\mathrm{G}$ par $\phi$ est encore un sous-groupe 
ab\'elien maximal de $\mathrm{G}.$ Les sous-groupes ab\'eliens maximaux de
$\mathrm{G}$ sont de l'un des types suivants
\begin{align*}
& \mathrm{T}=\{z+\alpha\hspace{0.1cm}\vert\hspace{0.1cm}\alpha\in\mathbb{C}\}, &&
\mathrm{D}_{z_0}=\{\alpha(z-z_0)+z_0\hspace{0.1cm}\vert\hspace{0.1cm}\alpha\in\mathbb{C}^*\}.
\end{align*}

\noindent On constate que $\mathrm{T}$ n'a pas d'\'el\'ement de torsion alors que les $\mathrm{D}_{z_0}$
en ont; par suite $\mathrm{T}$ est invariant par $\phi$ et tout $\mathrm{D}_{z_0}$ est 
envoy\'e par $\phi$ sur un $\mathrm{D}_{z'_0}.$ \`A conjugaison pr\`es par une translation on
peut supposer que $\phi(\mathrm{D}_0)=\mathrm{D}_0.$ Autrement dit il existe un morphisme
multiplica\-tif~$\kappa_1\colon\mathbb{C}^*\to~\mathbb{C}^*$
tel que 
\begin{align*}
& \phi(\alpha z)=\kappa_1(\alpha)z, &&\forall\hspace{0.1cm}\alpha\in\mathbb{C}^* 
\end{align*}
\noindent et un morphisme additif $\kappa_2\colon\mathbb{C}\to\mathbb{C}$ tel que 
\begin{align*}
& \phi(z+\beta)=z+\kappa_2(\beta),&& \forall\hspace{0.1cm}\beta\in\mathbb{C}.
\end{align*}
\noindent De plus, on a d'une part 
\begin{align*}
&\phi(\alpha z+\alpha)=\phi(\alpha z)\phi(z+1)
=\kappa_1(\alpha)z+\kappa_1(\alpha)\kappa_2(1)
\end{align*}
\noindent et d'autre part 
\begin{align*}
& \phi(\alpha z+\alpha)=\phi(z+\alpha)\phi(\alpha z)
=\kappa_1(\alpha)z+\kappa_2(\alpha);
\end{align*}
\noindent d'o\`u l'\'egalit\'e $\kappa_2(\alpha)=\kappa_1(\alpha)\kappa_2(1).$ 
Puisque $\kappa_2(1)$ est non nul, $\kappa_1$ est un morphisme 
additif et multiplicatif: $\kappa_1$ est un automorphisme du corps 
$\mathbb{C}.$ Notons $\mathrm{Aut}(\mathbb{C},+,\cdot)$ le groupe des
automorphismes du corps $\mathbb{C};$ posons $\gamma=\kappa_2(1),$ on
a 
\begin{align*}
&\phi(\alpha z+\beta)=\kappa_1(\alpha)z+\gamma\kappa_1(\beta)
=\kappa_1(\alpha z+\kappa_1^{-1}(\gamma)\beta)
=\kappa_1(\kappa_1^{-1}(\gamma)z\circ(\alpha z+\beta)\circ\kappa_1(\gamma)z);
\end{align*}
\noindent {\it i.e.} modulo l'action d'un automorphisme de corps, d'une translation
et d'une homoth\'etie, $\phi$ est trivial.

\noindent En s'appuyant sur la structure de produit amalgam\'e
de $\mathrm{Aut}[\mathbb{C}^2]$ on montre que tout automorphisme
de $\mathrm{Aut}[\mathbb{C}^2]$ s'obtient \`a partir d'une
conjugaison int\'erieure et de l'action d'un automorphisme du
corps $\mathbb{C}$ (\emph{voir} \cite{De6}). Le groupe de Cremona ne 
poss\`ede pas de telle structure n\'eanmoins en \'etudiant les sous-groupes 
ab\'eliens maximaux non d\'enombrables du 
groupe de Cremona on peut d\'eterminer $\mathrm{Aut}(\mathrm{Bir}(\mathbb{P}^2)).$

\begin{thm}[\cite{De}]\label{autaut}
{\sl Le groupe des automorphismes 
ext\'erieurs du groupe de Cremona s'identifie au groupe des 
automorphismes du corps $\mathbb{C}.$ 

\noindent Autrement dit: soit $\phi$ dans $\mathrm{Aut}(\mathrm{Bir}(
\mathbb{P}^2));$ il existe $\kappa$ dans $\mathrm{Aut}(\mathbb{C},+,\cdot)$ 
et $\psi$ dans $\mathrm{Bir}(\mathbb{P}^2)$ tels que 
\begin{align*}
&\phi(f)=\kappa(\psi f\psi^{-1}), &&\forall\hspace{1mm} f\in\mathrm{Bir}(\mathbb{P}^2).
\end{align*}}
\end{thm}

\noindent L'id\'ee est la suivante. Consid\'erons un automorphisme $\phi$ de 
$\mathrm{Bir}(\mathbb{P}^2)$ et  un sous-groupe ab\'elien
ma\-ximal $\mathrm{G}$ de $\mathrm{Bir}(\mathbb{P}^2);$ alors $\phi(\mathrm{G})$ 
est un sous-groupe ab\'elien maximal de $\mathrm{Bir}(\mathbb{P}^2).$
On modifie~$\phi,$ uniquement \`a l'aide de conjugaisons int\'erieures et
d'action par automorphismes de corps, de sorte que 
\begin{align*}
&\phi(\mathrm{G}_1)=\mathrm{G}_1, &&\ldots, &&\phi(\mathrm{G}_n)=\mathrm{G}_n
\end{align*}
\noindent les $\mathrm{G}_i$ d\'esignant certains sous-groupes ab\'eliens maximaux de $\mathrm{Bir}
(\mathbb{P}^2)$ bien choisis. On en d\'eduit que 
$\mathrm{PGL}_3(\mathbb{C})$ est invariant point par point par $\phi$ puis que $\phi(\sigma)=
\sigma;$ le th\'eor\`eme de N\oe ther permet alors de conclure.

\medskip 

\noindent Avant de pr\'esenter un des points cl\'es de la d\'emonstration du th\'eor\`eme
\ref{autaut} introduisons quelques notations et d\'efinitions. Si $\mathrm{S}$
est une surface complexe compacte, un \textsl{feuilletage}
$\mathcal{F}$ sur $S$ est la donn\'ee d'une famille $(\chi_i)_i$ de champs de vecteurs
holomorphes \`a z\'eros isol\'es d\'efinis sur les ouverts~$\mathcal{U}_i$ d'un recouvrement de $\mathrm{S}.$
Les champs $\chi_i$ sont soumis \`a des conditions de compatibi\-lit\'e~:  il existe~$g_{ij}$
dans $\mathcal{O}^*(\mathcal{U}_i\cap\mathcal{U}_j)$ tel que $\chi_i$
co\"{\i}ncide avec $g_{ij}\chi_j$ sur $\mathcal{U}_i\cap~\mathcal{U}_j$.
Notons qu'un champ de vecteurs m\'eromorphe non trivial sur $\mathrm{S}$
d\'efinit un tel feuilletage.

\begin{lem}\label{amel}
{\sl Soit $\mathrm{G}$ un sous-groupe ab\'elien 
non d\'enombrable du groupe de Cremona.
Il existe un champ de vecteurs rationnel $\chi$ tel
que 
\begin{align*}
&f_*\chi=\chi, &&\forall\hspace{1mm}f\in\mathrm{G}.
\end{align*}

\noindent En particulier $\mathrm{G}$ pr\'eserve un feuilletage.}
\end{lem}

\begin{proof}[{\sl D\'emonstration}]
Puisque le groupe $\mathrm{G}$ n'est pas d\'enombrable, il existe un entier
$\nu$ tel que l'ensemble
\begin{align*}
\mathrm{G}_\nu=\{f\in\mathrm{G}\hspace{1mm}|\hspace{1mm}\deg f=\nu\}
\end{align*}
\noindent ne soit pas d\'enombrable. Par suite, l'adh\'erence de Zariski
$\overline{\mathrm{G}_\nu}$ de $\mathrm{G}_\nu$ dans  
\begin{align*}
\mathrm{Bir}(\mathbb{P}^2)_\nu=\{f\in\mathrm{Bir}(\mathbb{P}^2)\hspace{1mm}|\hspace{1mm}
\deg f\leq \nu\}
\end{align*}
\noindent est un ensemble alg\'ebrique de dimension
sup\'erieure ou \'egale \`a un; consid\'erons alors une courbe dans 
$\overline{\mathrm{G}_\nu},$ {\it i.e.} une application 
\begin{align*}
&\eta\colon\mathbb{D}\to\overline{\mathrm{G}_\nu}, && t\mapsto \eta(t).
\end{align*}
\noindent Notons que les
\'el\'ements de $\overline{\mathrm{G}_\nu}$ sont des applications
rationnelles qui commutent.
D\'efinissons alors le champ de vecteurs rationnel $\chi$
en tout point $m$ n'appartenant pas au lieu d'ind\'etermination
de $\eta(0)^{-1}$ par 
\begin{align*}
\chi(m)=\frac{\partial \eta(s)}{\partial s}\Big|_{s=0}
(\eta(0)^{-1}(m)).
\end{align*}
\noindent Soit $f$ dans $\overline{\mathrm{G}_\nu};$ en d\'erivant l'identit\'e $f\eta(s)
f^{-1}(m)=\eta(s)(m)$ par rapport \`a $s$ \`a $m$ fix\'e, on obtient~: $f_*\chi=\chi.$
Autrement dit $\chi$ est invariant par les \'el\'ements de 
$\overline{\mathrm{G}_\nu};$ on constate qu'en fait~$\chi$ est invariant 
par tout \'el\'ement de $\mathrm{G}.$
\end{proof}

\noindent Ainsi \`a chaque fois que l'on se donne un sous-groupe ab\'elien
non d\'enombrable $\mathrm{G}$ de $\mathrm{Bir}(\mathbb{P}^2)$ 
on h\'erite d'un feuilletage sur le plan projectif complexe invariant par $\mathrm{G}.$

\noindent Brunella, McQuillan et Mendes ont 
classifi\'e, \`a \'equivalence birationnelle pr\`es, les feuilletages 
holomorphes singuliers sur une surface complexe, compacte et projective 
(\cite{Br, McQ, Me}). Soit 
$\mathrm{S}$ une surface projective munie d'un feuilletage $\mathcal{F};$ notons 
$\mathrm{Bir}(\mathrm{S},\mathcal{F})$ (resp.~$\mathrm{Aut}(\mathrm{S},\mathcal{F})$) 
le groupe des transformations birationnelles (resp. holomorphes) 
laissant le feuilletage $\mathcal{F}$ invariant sur la surface $\mathrm{S}.$ G\'en\'eriquement
$\mathrm{Bir}(\mathrm{S},\mathcal{F})$ co\"incide avec $\mathrm{Aut}(\mathrm{S},\mathcal{F})$ et est fini. Dans \cite{CaFa} 
Cantat et Favre \'etudient les feuilletages qui ne v\'erifient pas cette 
philosophie. Avant d'\'enoncer la classification qu'ils obtiennent rappelons
ce qu'est une surface de Kummer (g\'en\'eralis\'ee). Soit~$\Lambda$ un r\'eseau de~$\mathbb{C}^2;$ il induit un tore complexe $T=
\mathbb{C}^2/\Lambda$ de dimension $2.$ Par exemple le produit d'une
courbe elliptique par elle-m\^eme est un tore 
complexe. Une application affine $f$ pr\'eservant~$\Lambda$ induit un 
automorphisme du tore $T.$ Supposons que la partie 
lin\'eaire de $f$ soit d'ordre infini, alors
\begin{itemize}
\item ou bien la partie lin\'eaire de $f$ est hyperbolique auquel cas
$f$ induit un automorphisme d'Anosov pr\'eservant deux feuilletages
lin\'eaires;

\item ou bien la partie lin\'eaire de $f$ est unipotente et $f$ pr\'eserve
une fibration elliptique.
\end{itemize}

\noindent Il arrive que le tore $T$ poss\`ede un groupe fini 
d'automorphismes normalis\'e par $f.$ Notons $\widetilde{T/\mathrm{G}}$ la 
d\'esingularis\'ee de~$T/\mathrm{G}.$  L'automorphisme $\widetilde{f}$ induit par $f$ 
sur $\widetilde{T/\mathrm{G}}$  pr\'eserve 
\begin{itemize}
\item une fibration elliptique lorsque la partie lin\'eaire de $f$ est unipotente ;

\item les transform\'es des feuilletages stables et instables 
de $f$ lorsque la partie lin\'eaire de $f$ est hyperbolique.  
\end{itemize}
On dit que $\widetilde{T/\mathrm{G}}$ est une \textsl{surface de Kummer} lorsque $\mathrm{G}=\{\mathrm{id},\hspace{0.1cm} (-x,-y)\};$ par analogie pour $\mathrm{G}$ 
quelconque $\widetilde{T/\mathrm{G}}$ est une \textsl{surface de Kummer 
g\'en\'eralis\'ee}.

\noindent Les $5$ classes de feuilletages obtenus par Cantat et 
Favre sont les suivantes 
\begin{itemize}
\item $\mathcal{F}$ est invariant par un champ de vecteurs holomorphe;

\item $\mathcal{F}$ est une fibration elliptique;

\item $\mathrm{S}$ est une surface de Kummer g\'en\'eralis\'ee et $\mathcal{F}$ est
la projection sur $\mathrm{S}$ du feuilletage stable ou instable d'un certain
automorphisme d'Anosov;

\item $\mathcal{F}$ est une fibration rationnelle;

\item il existe des entiers $p,$ $q,$ $r$ et $s$ tels qu'\`a rev\^etement
fini pr\`es on ait
\begin{align*}
& \mathrm{Bir}(S,\mathcal{F})=\{(x^py^q,x^ry^s),\hspace{0.1cm}(\alpha x,\beta y)\hspace{0.1cm} \vert\hspace{0.1cm}\alpha,\hspace{0.1cm}\beta\in\mathbb{C}^*\}.
\end{align*}
\end{itemize}

\noindent G\'en\'eriquement une application bim\'eromorphe qui laisse un 
feuilletage invariant pr\'eserve un pinceau de courbes rationnelles ou 
elliptiques; Cantat et Favre d\'ecrivent les transformations birationnelles
dont la dynamique n'est pas triviale et qui ne satisfont pas cette g\'en\'eralit\'e.

\noindent Cette \'etude conduit \`a la caract\'erisation des sous-groupes 
ab\'eliens maximaux non d\'enombrables du groupe de Cremona.

\begin{thm}[\cite{De}]
{\sl Si $\mathrm{G}$ est un sous-groupe ab\'elien maximal non 
d\'enombrable de~$\mathrm{Bir}(\mathbb{P}^2),$ il 
v\'erifie l'une des propri\'et\'es suivantes
\begin{itemize}
\item $\mathrm{G}$ poss\`ede des \'el\'ements de torsion;

\item $\mathrm{G}$ est conjugu\'e \`a un sous-groupe de $\mathrm{dJ}.$
\end{itemize}}
\end{thm}

\noindent Cet \'enonc\'e permet de d\'eterminer l'image par $\phi$ des
groupes ab\'eliens maximaux non d\'enombrables
\begin{align*}
&\mathrm{dJ_a}=\{(x+a(y),y)\hspace{1mm}\vert\hspace{1mm} a \in\mathbb{C}(y)\}, &&
\mathrm{T}=\{(x+\alpha,y+\beta)\hspace{1mm}\vert\hspace{1mm} \alpha,\hspace{1mm}\beta
\in\mathbb{C}\};
\end{align*}
\noindent modulo une conjugaison int\'erieure ils sont invariants par $\phi.$ Ensuite
on d\'emontre que, quitte \`a faire agir un automorphisme du corps $\mathbb{C}$
et un automorphisme int\'erieur, on a: 
$\mathrm{T}$ et 
\begin{align*}
\mathrm{D}=\{(\alpha x,\beta y)\,\vert\,\alpha,
\,\beta\in\mathbb{C}^*\}
\end{align*}
\noindent sont invariants point par point par $\phi$
et $\phi(\mathrm{dJ_a})=\mathrm{dJ_a};$ on en 
d\'eduit que $\mathrm{PGL}_3(\mathbb{C})$ et $\sigma$ le sont aussi
d'o\`u l'\'enonc\'e.

\bigskip

\section{Un peu de dynamique}\label{dyn}

\bigskip 

\subsection{Classification des transformations birationnelles}

\noindent Soit $f$ une transformation birationnelle sur une surface complexe
compacte $\mathrm{S}.$ La transformation $f$ est dite
\textsl{alg\'ebriquement stable} s'il n'existe pas de courbe~$\mathcal{C}$ dans $\mathrm{S}$ telle
que $f^k(\mathcal{C})$ appartienne \`a $\mathrm{Ind}(f)$ pour un
certain entier~$k\geq 0.$ Dit autrement une transformation est
alg\'ebriquement stable si la situation suivante n'arrive pas

\begin{figure}[H]
\begin{center}
\input{as.pstex_t}
\end{center}
\end{figure}

\noindent Les deux conditions suivantes sont 
\'equivalentes (\cite{DiFa}).
\smallskip
\begin{itemize}
\item Il n'existe pas de courbe $\mathcal{C}$ dans $\mathbb{P}^2(\mathbb{C})$ telle
que $f^k(\mathcal{C})$ appartienne \`a $\mathrm{Ind}(f)$ pour un
certain entier $k\geq 0;$

\item pour tout entier $n$ on a $\deg f^n=(\deg f)^n.$  
\end{itemize}
\smallskip

\noindent D'apr\`es Diller et Favre toute transformation birationnelle sur 
une surface complexe compacte est alg\'ebriquement stable modulo un 
changement de coordonn\'ees.

\begin{pro}[\cite{DiFa}]
{\sl Soit $f\colon \mathrm{S}\dashrightarrow \mathrm{S}$ une transformation birationnelle sur
une surface complexe compacte $\mathrm{S}.$ Il existe une suite d'\'eclatements 
$\pi\colon\widetilde{\mathrm{S}}\to\mathrm{S}$ telle que $\pi^{-1}f\pi\colon
\widetilde{\mathrm{S}}\dashrightarrow \widetilde{\mathrm{S}}$ soit alg\'ebriquement stable.}
\end{pro}

\noindent Donnons une id\'ee de la d\'emonstration de cet 
\'enonc\'e.
Supposons que $f$ ne soit pas alg\'ebriquement 
stable; il existe alors une cour\-be~$\mathcal{C}$ et un 
entier $k$ tels que
\begin{figure}[H]
\begin{center}
\input{as4.pstex_t}
\end{center}
\end{figure}

\noindent L'id\'ee mise en \oe uvre par Diller et 
Favre est la suivante: quitte \`a \'eclater les $p_i$
l'image de $\mathcal{C}$ est, pour $i=1,$ $\ldots,$ 
$k,$ une courbe.

\bigskip

\noindent Soient $f$ et $g$ deux \'el\'ements du groupe de Cremona: en g\'en\'eral
$\deg gfg^{-1}\not=\deg f,$ {\it i.e.} le degr\'e n'est pas un invariant birationnel.
Par contre il existe deux constantes positives $\alpha$ et $\beta$ telles que
\begin{align*}
&\alpha\deg f^n\leq\deg(gf^ng^{-1})\leq\beta\deg f^n, &&\forall\hspace{1mm}n\in\mathbb{N};
\end{align*}
\noindent le taux de croissance des degr\'es est un invariant birationnel. On introduit donc la notion 
suivante: soit $f$ dans $\mathrm{Bir}(\mathbb{P}^2),$ 
le \textsl{premier degr\'e dynamique} de $f$ est la quantit\'e
\begin{align*}
\lambda(f)=\liminf (\deg f^n)^{1/n}.
\end{align*}

\noindent Lorsque $f$ est une transformation de Cremona 
alg\'ebriquement stable on a $\lambda(f)=\deg f.$

\begin{thm}[\cite{DiFa}]\label{jefcha}
{\sl Une transformation $f$ de Cremona satisfait, \`a conjugaison birationnelle
pr\`es, une et une seule des propri\'et\'es suivantes
\smallskip
\begin{itemize}
\item la suite $(\deg f^n)_n$ est born\'ee, $f$ est
un automorphisme d'une surface rationnelle $\mathrm{S}$ et un it\'er\'e de $f$
appartient \`a $\mathrm{Aut}^0(\mathrm{S}),$ la composante neutre de $\mathrm{Aut}(\mathrm{S});$

\item $\deg f^n\sim n$ alors $f$ n'est
pas un automorphisme et $f$ pr\'eserve une unique fibration qui est rationnelle;

\item $\deg f^n\sim n^2$ auquel cas $f$ est un automorphisme pr\'eservant une 
unique fibration qui est elliptique;

\item $\deg f^n\sim\alpha^n,$ $\alpha>1.$
\end{itemize}
\smallskip
\noindent Dans les trois premiers cas on a $\lambda(f)=1,$ dans
le dernier $\lambda(f)>1.$ }
\end{thm}

\begin{defi}
Soit $f$ une transformation birationnelle du plan projectif complexe
dans lui-m\^eme. 
\smallskip
\begin{itemize}
\item si la suite $(\deg f^n)_n$ est born\'ee, on dit que $f$ est \textsl{virtuellement
isotope \`a l'identit\'e};

\item si $\deg f^n\sim n,$ on dit que $f$ est 
\textsl{un twist de de Jonqui\`eres};

\item si $\deg f^n\sim n^2,$ on dit que $f$ est 
\textsl{un twist d'Halphen};

\item si $\deg f^n\sim\alpha^n$ avec $\alpha>1,$ on parle de
transformation \textsl{entropique}. 
\end{itemize}
\end{defi}

\begin{egs}
\begin{itemize}
\item Toute transformation birationnelle d'ordre
fini, tout automorphisme du plan projectif complexe est virtuellement 
isotope \`a l'identit\'e.

\item La transformation $f=(xz:xy:z^2)$ est un twist de de Jonqui\`eres (dans la carte
affi\-ne~$z=1$ son it\'er\'e $n$-i\`eme s'\'ecrit $(x,x^ny)$).

\item Le prolongement d'un automorphisme de H\'enon
g\'en\'eralis\'e $(y,y^2-\delta x),$ avec $\delta\in\mathbb{C}^*,$ est entropique.

\item  Soit $f_M$ la transformation de la forme 
\begin{align*}
&f_M=(x^ay^b,x^cy^d), && M=\left[
\begin{array}{cc}
a & b\\
c & d
\end{array}
\right]\in\mathrm{SL}_2(\mathbb{Z});
\end{align*}
\noindent on a l'alternative
\begin{itemize}
\item ou bien $\vert\mathrm{tr}\,M\vert\leq 2$
et $\lambda(f_M)=1;$

\item ou bien $M$ a pour valeurs propres
$\lambda$ et $\lambda^{-1}$ avec $\lambda^{-1}<
1<\lambda$ et $\lambda(f_M)=\lambda.$
\end{itemize}
Par suite $f_M$ est virtuellement isotope \`a l'identit\'e (resp. un twist 
de de Jonqui\`eres, resp. un twist d'Halphen, resp. entropique) d\`es que
$\left[
\begin{array}{cc}
a & b\\
c & d
\end{array}
\right]$ est elliptique (resp. parabolique, resp. parabolique, resp. hyperbolique).

\noindent Supposons que $a,$ $b,$ $c,$ $d$ soient strictement 
positifs; $f_M^*$ n'est pas une isom\'etrie. Remarquons que 
$\deg f_M^2<(\deg f_M)^2;$ en particulier $f_M\in\mathrm{Bir}
(\mathbb{P}^2)$ n'est pas alg\'ebriquement stable. 
On peut voir $f_M$ comme un \'el\'ement de $\mathbb{P}^1
(\mathbb{C})\times\mathbb{P}^1(\mathbb{C})$ avec $(x,y)$
dans $\mathbb{P}^1(\mathbb{C})\times\mathbb{P}^1(\mathbb{C}).$ On~a 
\begin{align*}
&\mathrm{H}^2(\mathbb{P}^1(\mathbb{C})\times\mathbb{P}^1
(\mathbb{C}),\mathbb{Z})=\mathbb{Z}\mathrm{L}_x\oplus\mathbb{Z}
\mathrm{L}_y , &&\mathrm{L}_x=\{0\}\times\mathbb{P}^1(\mathbb{C}), 
\hspace{1mm}\mathrm{L}_y=\mathbb{P}^1(\mathbb{C})\times\{0\};
\end{align*}
\noindent d'o\`u
\begin{align*}
& f_M^*\mathrm{L}_x=a\mathrm{L}_x+b\mathrm{L}_y, &&f_M^*
\mathrm{L}_x=c\mathrm{L}_x+d\mathrm{L}_y.
\end{align*}
\noindent Il s'en suit que la matrice de $f_M^*$ co\"incide avec 
$\transp M.$ Finalement $f_M$ vue comme 
transformation birationnelle de $\mathbb{P}^1(\mathbb{C})
\times\mathbb{P}^1(\mathbb{C})$ est alg\'ebriquement 
stable et $\lambda(f_M)$ est la plus grande valeur propre de $M.$
\end{itemize}
\end{egs}

\begin{rem}
La notion de premier degr\'e dynamique et le th\'eor\`eme 
\ref{jefcha} se g\'en\'eralisent au cas des transformations 
bim\'eromorphes sur des surfaces complexes compactes k\"{a}hl\'eriennes.

\noindent Une transformation
rationnelle $f$ sur une surface complexe compacte k\"{a}hl\'erienne 
$\mathrm{S}$ d\'efinit un op\'erateur lin\'eaire $f^*$ sur les groupes de 
cohomologie de $\mathrm{S}$ pr\'eservant la d\'ecomposition de Hodge
\begin{align*}
\mathrm{H}^k(\mathrm{S},\mathbb{C})=\bigoplus_{i+j=k}\mathrm{H}^{i,j}(\mathrm{S},\mathbb{C}).
\end{align*}
\noindent On d\'efinit le premier degr\'e dynamique de $f$ par
\begin{align*}
\lambda(f)=\lim_{n\to +\infty}\|(f^n)^*\|^{1/n}
\end{align*}
\noindent o\`u $\|.\|$ d\'esigne une norme sur $\mathrm{End}(\mathrm{H}^{1,1}(\mathbb{R},\mathrm{S}));$
toujours d'apr\`es \cite{DiFa} cette quantit\'e est un invariant birationnel. Elle 
co\"incide sur $\mathbb{P}^2(\mathbb{C})$ avec la notion d\'ej\`a introduite.
\end{rem}

\noindent Si $f$ est une transformation bim\'eromorphe sur une 
surface k\"ahl\'erienne, 
$\lambda(f)$ est un entier alg\'ebrique dont les conjugu\'es
$\mu_i$ sont de module inf\'erieur \`a $1$ (\emph{voir} \cite{DiFa}). En particulier $\lambda
(f)$ est ou bien un nombre de Pisot (on a $\vert\mu_i\vert<1$
pour tout $i$), ou 
bien un nombre de Salem (les conjugu\'es de $\lambda$ 
sont $\lambda,$ $\lambda^{-1}$ et $\mu_i$ avec $\vert\mu_i\vert=1$).

\medskip

\noindent La d\'emonstration du th\'eor\`eme \ref{jefcha} se d\'ecompose comme suit
\begin{itemize}
\item $(\|f^{n*}\|)_n$ est born\'ee si et seulement si $f$ est conjugu\'e
\`a un automorphisme dont un it\'er\'e est isotope \`a l'identit\'e;

\item si $(\|f^{n*}\|)_n$ n'est pas born\'ee, $f$ est conjugu\'e
\begin{itemize}
\item ou bien \`a un automorphisme,

\item ou bien \`a une transformation birationnelle qui pr\'eserve une 
fibration rationnelle, de plus $\|f^{n*}\|\sim~n;$
\end{itemize}

\item si $f$ est un automorphisme et si $(\|f^{n*}\|)_n$ n'est pas born\'ee,
$f$ pr\'eserve une fibration ellipti\-que~(\cite{Gi2, Ca2, Ca4});

\item si $f$ est un automorphisme qui pr\'eserve une fibration elliptique
et si $(\|f^{n*}\|)_n$ est non born\'ee, alors $\|f^{n*}\|\sim n^2.$ 

\item si $f$ pr\'eserve deux fibrations g\'en\'eriquement transverses, 
la suite $(\|f^{n*}\|)_n$ est born\'ee.
\end{itemize}

\bigskip

\noindent On a une situation analogue, et m\^eme plus pr\'ecise,
pour le groupe $\mathrm{Aut}[\mathbb{C}^2].$ \`A partir du th\'eor\`eme~\ref{jung}, Friedland et Milnor donnent une classification des automorphismes polynomiaux
de~$\mathbb{C}^2.$ Avant de l'\'enoncer on introduit le sous-ensemble
$\mathrm{H}$ de $\mathrm{Aut}[\mathbb{C}^2]$ constitu\'e des 
applications \textsl{de type H\'enon}, {\it i.e.} des automorphismes
de la forme 
\begin{align*}
& f=\phi h_1\ldots h_n\phi^{-1}, && \phi\in\mathrm{Aut}[\mathbb{C}^2],
&& h_i \text{ application de H\'enon g\'en\'eralis\'ee}.
\end{align*}

\begin{thm}[\cite{FM}]\label{FriedlandMilnor}
{\sl Soit $f$ un automorphisme polynomial de $\mathbb{C}^2.$ On a l'alternative
suivante
\smallskip
\begin{itemize}
\item ou bien $f$ est conjugu\'e \`a un \'el\'ement de $\mathtt{E};$

\item ou bien $f$ appartient \`a $\mathrm{H}.$
\end{itemize}}
\end{thm}

\noindent On peut reformuler le th\'eor\`eme \ref{FriedlandMilnor} en utilisant le 
premier degr\'e dynamique: si $f$ est un automorphisme polynomial de 
$\mathbb{C}^2,$ alors
\smallskip
\begin{itemize}
\item $f$ est conjugu\'e \`a un \'el\'ement de 
$\mathtt{E}$ si et seulement si $\lambda(f)=1;$

\item $f$ est de type H\'enon si et seulement
si $\lambda(f)>1.$
\end{itemize}
\smallskip

\noindent Soit $f$ un automorphisme polynomial du plan
complexe. \`A $f$ on peut associer le sous-arbre $\mathrm{Fix}(f)$ de 
$\mathcal{T}$ constitu\'e des sommets et des ar\^etes fix\'es
par l'action de $f.$ Le th\'eor\`eme \ref{FriedlandMilnor}
correspond aussi \`a l'alternative: $\mathrm{Fix}(f)\not=\emptyset$
ou $\mathrm{Fix}(f)=\emptyset;$ en effet le stabilisateur 
du sommet $f\mathtt{E}$ (resp. du sommet~$f\mathtt{A},$ resp. 
de l'ar\^ete $f\mathtt{S}$) est le groupe $f\mathtt{E}f^{-1}$ 
(resp. $f\mathtt{A}f^{-1},$ resp. $f\mathtt{S}f^{-1}$). Lorsque
$f$ est dans $\mathrm{H}$ il existe une 
unique g\'eod\'esique infinie, {\it i.e.} un sous-arbre 
isomorphe \`a $\mathbb{R},$ not\'ee $\mathrm{Geo}(f),$
telle que l'action de~$f$ restreinte \`a $\mathrm{Geo}(f)$ soit une 
translation.

\bigskip

\noindent Plus pr\'ecis\'ement on a la \og dichotomie\fg\hspace{1mm} suivante.
\medskip
\begin{itemize}
\item Ou bien $f$ appartient \`a $\mathtt{E}$ \`a conjugaison pr\`es,

\noindent le premier degr\'e dynamique de $f$ vaut $1,$ 

\noindent la transformation $f$ pr\'eserve une fibration rationnelle,

\noindent le centralisateur de $f$ est non d\'enombrable,

\noindent l'ensemble $\mathrm{Fix}(f)$ est non vide; 

\item ou bien $f$ est dans $\mathrm{H},$

\noindent le premier degr\'e de $f$ est strictement sup\'erieur \`a $1,$ 

\noindent l'\'el\'ement $f$ ne pr\'eserve pas de courbe rationnelle,

\noindent le centralisateur de $f$ est d\'enombrable,

\noindent l'ensemble $\mathrm{Fix}(f)$ est vide, 

\noindent l'application $f$ admet une infinit\'e de points p\'eriodiques hyperboliques.
\end{itemize}

\medskip

\noindent La situation est-elle la m\^eme pour le groupe de Cremona ?

\begin{thm}[\cite{Ca}] 
{\sl  Soit $f$ une transformation de Cremona entropique. Si $g$ est une 
transformation birationnelle de $\mathbb{P}^2(\mathbb{C})$ dans lui-m\^eme 
qui commute avec $f$ il existe deux entiers~$m$ dans $\mathbb{N}^*$ et $n$
dans $\mathbb{Z}$ tels que $g^m=f^n.$}
\end{thm}

\noindent N\'eanmoins contrairement au cas des automorphismes
polynomiaux il n'y a pas d'\'equivalence entre le fait d'\^etre 
entropique et d'avoir un centralisateur d\'enombrable: consid\'erons la 
famille de transformations birationnelles $(f_{\alpha,\beta})$ donn\'ee par
\begin{align*}
& \mathbb{P}^2(\mathbb{C})\dashrightarrow\mathbb{P}^2(\mathbb{C}), && (x:y:z)\mapsto((\alpha x+y)z:\beta y
(x+z):z(x+z)), && \alpha,\hspace{0.1cm}\beta\in\mathbb{C}^*,
\end{align*}
\noindent soit dans la carte affine $z=1$ 
\begin{align*}
& f_{\alpha,\beta}= \left(\frac{\alpha x+y}{x+1},\beta y\right). 
\end{align*}

\begin{thm}[\cite{De2}] 
{\sl  Les $f_{\alpha,\beta}$ sont des twists de de Jonqui\`eres. 

\noindent Supposons que $\alpha$ et $\beta$ soient de module $1$ et g\'en\'eriques.

\noindent Si $g$ commute \`a $f_{\alpha,\beta},$ alors $g$ est une puissance
de $f_{\alpha,\beta};$ en particulier le centralisateur de $f_{\alpha,\beta}$ est 
d\'enombrable.

\noindent Les \'el\'ements $f^2_{\alpha,\beta}$ poss\`edent deux points fixes
$p_1,$ $p_2$ et 
\smallskip
\begin{itemize}
\item il existe un voisinage $\mathcal{V}_1$ de $p_1$ sur lequel $f_{\alpha,\beta}$
est conjugu\'e \`a $(\alpha x,\beta y);$ en particulier l'adh\'erence de l'orbite
d'un point de $\mathcal{V}_1$ sous l'action de $f_{\alpha,\beta}$ est un tore de dimension $2;$ 

\item il existe un voisinage $\mathcal{V}_2$ de $p_2$ tel que la dynamique de 
$f_{\alpha,\beta}^2$ soit localement lin\'eaire sur $\mathcal{V}_2;$ l'adh\'erence
de l'orbite g\'en\'erique d'un point de $\mathcal{V}_2$ sous l'action de $f^2_{\alpha,\beta}$ est un cercle.
\end{itemize}}
\end{thm} 

\begin{rem} 
La
droite $z=0$ est contract\'ee par $f_{\alpha,\beta}$ sur le point $(0:1:0)$ qui
est lui-m\^eme \'eclat\'e sur $z=0;$ la transformation $f_{\alpha,
\beta}$ n'est pas alg\'ebriquement stable, c'est pour cette raison qu'on a 
consid\'er\'e $f_{\alpha,\beta}^2$ plut\^ot que $f_{\alpha,\beta}.$
\end{rem}

\noindent On va \'etablir un r\'esultat qui nous servira 
au paragraphe suivant.

\begin{lem}[\cite{De3}]\label{commut}
{\sl Soient $f$ et $g$ deux transformations birationnelles sur une surface~$\mathrm{S}.$ 
Supposons que $f$ et $g$ soient virtuellement isotopes \`a l'identit\'e. Si 
$f$ et $g$ commutent, il existe une surface $\tilde{\mathrm{S}},$
une transformation birationnel\-le~$\phi\colon \mathrm{S}\dashrightarrow
\tilde{\mathrm{S}}$ et un entier $k$ tels que 
\begin{align*}
& \phi f\phi^{-1},\hspace{1mm}\phi g\phi^{-1}\in\mathrm{Aut}(\tilde{\mathrm{S}})&& 
\text{et}&&\phi f^k\phi^{-1},\hspace{1mm}\phi g^k\phi^{-1}\in\mathrm{Aut}^0(\tilde{\mathrm{S}}).
\end{align*}}
\end{lem}

\begin{proof}[{\sl D\'emonstration}]
Par hypoth\`ese il existe une surface $\tilde{\mathrm{S}},$ une
transformation birationnel\-le~$\zeta\colon \mathrm{S}\dasharrow~\tilde{\mathrm{S}}$ et
un entier $n$ tels que $\zeta^{-1}f^n\zeta$ soit un automorphisme
de $\tilde{\mathrm{S}}$ isotope \`a l'identit\'e. Pla\c{c}ons-nous 
sur~$\tilde{\mathrm{S}};$ pour simplifier on note encore $f$
l'automorphisme $\zeta^{-1}f^n\zeta$ et $g$ la transformation
$\zeta^{-1}g\zeta.$

\noindent Commen\c{c}ons par montrer qu'il existe $\eta
\colon Y \dashrightarrow\tilde{\mathrm{S}}$ birationnel tel que
$\eta^{-1}f^ \ell\eta$ soit un automorphisme de $Y$ isotope \`a
l'identit\'e pour un certain $\ell$ et $\eta^{-1}g\eta$ soit
alg\'ebriquement stable. On note~$\nu(g)$ le nombre minimum
d'\'eclatements n\'ecessaires pour rendre $g$ alg\'ebriquement
stable. La preuve proc\`ede par induction sur $\nu(g).$

\noindent Si $\nu(g)$ est nul, alors $\eta=\mathrm{id}$ convient.

\noindent Supposons le lemme d\'emontr\'e pour les transformations
$f$ et $g$ satisfaisant $\nu(g)\leq j;$ consid\'e\-rons un couple~$(\tilde{f},
\tilde{g})$ v\'erifiant les hypoth\`eses de l'\'enonc\'e et 
$\nu(\tilde{g})=j+1.$ Puisque $\tilde{g}$ n'est pas alg\'ebriquement
stable, il existe une courbe $V$ dans $\mathrm{Exc}(\tilde{g})$ et
un entier $q$ tels que $\tilde{g}^q( V)$ soit un point
d'ind\'etermination $p$ de $\tilde{g}$. Comme $\tilde{f}$ et
$\tilde{g}$ commutent, $\tilde{f}^k$ fixe les composantes
irr\'eductibles de $\mathrm{Ind}(\tilde{g})$ pour un certain
entier $k.$ 
Consid\'erons $\kappa$ l'\'eclatement de $\widetilde{\mathrm{S}}$ au point $p;$ ce point
\'etant fix\'e par $\tilde{f}^k,$ d'une part $\kappa^{-1}
\tilde{f}^k\kappa$ est un automorphisme, d'autre part
$\nu(\kappa^{-1}\tilde{g} \kappa)=j.$ Alors, par hypoth\`ese de
r\'ecurrence, il existe $\eta\colon Y
\dashrightarrow\tilde{\mathrm{S}}$ et $\ell$ tels que
$\eta^{-1}\tilde{f}^\ell \eta$ soit un automorphisme isotope \`a
l'identit\'e et $\eta^{-1} \tilde{g}\eta$ soit alg\'ebriquement
stable. \medskip

\noindent Posons $\overline{f}=\eta^{-1}f^\ell\eta$ et
$\overline{g}=\eta^{-1}g\eta.$ La premi\`ere \'etape 
consiste \`a consid\'erer $\varepsilon_1$ la contraction d'une
courbe de $\mathrm{Exc}(\overline{g}^{-1})~;$ puisque les courbes
contract\'ees par $\overline{g}^{-1}$ sont, d'apr\`es
\cite{DiFa}, d'auto-intersection n\'egative et que
$\overline{f}$ est isotope \`a l'identit\'e, elles sont fix\'ees
par $\overline{f}$ donc par $\varepsilon_1\overline{f}
\varepsilon_1^{-1}.$ La $i$-\`eme \'etape consistant \`a
r\'ep\'eter la premi\`ere avec $\varepsilon_{i-1} \ldots
\varepsilon_1\overline{f}\varepsilon_1^{-1}\ldots
\varepsilon_{i-1}^{-1}$ et $\varepsilon_{i-1} \ldots
\varepsilon_1\overline{g}\varepsilon_1^{-1}\ldots
\varepsilon_{i-1}^{-1},$ on obtient le r\'esultat souhait\'e.
D'apr\`es \cite{DiFa} le proc\'ed\'e termine. Toujours par
l'argument de Diller et Favre, une puissance de
$\varepsilon^{-1}g\varepsilon$ est isotope \`a l'identit\'e.
\end{proof}

\bigskip

\subsection{Caract\'erisation des courbes invariantes par certaines
transformations birationnelles}

\noindent Au \S\ref{gpefini} on a \'evoqu\'e l'article \cite{DJS} sans en 
pr\'eciser les r\'esultats; le vocabulaire ayant \'et\'e introduit,  
donnons quelques-uns de leurs \'enonc\'es. 

\begin{thm}[\cite{DJS}]
{\sl Soient $S$ une surface projective complexe et $f$ une transformation
birationnelle sur $S$ alg\'ebriquement stable et entropique. Soit 
$\mathcal{C}$ une courbe connexe invariante par~$f.$

\smallskip

\noindent La courbe $\mathcal{C}$ est de genre $0$ ou $1.$

\smallskip

\noindent Supposons que $\mathcal{C}$ soit de genre $1.$ 
Quitte \`a contracter certaines courbes sur $\mathrm{S},$ il existe une 
$2$-for\-me~$\omega$ m\'eromorphe telle que 
\begin{itemize}
\item $f^*\omega=\alpha\omega,$ $\alpha\in\mathbb{C};$ 

\item et $-C$ soit le diviseur de $\omega.$ 
\end{itemize}
La constante $\alpha$ est d\'etermin\'ee par la 
courbe $\mathcal{C}$ et la restriction de $f$ \`a $\mathcal{C}.$ }
\end{thm}

\noindent Les auteurs s'int\'eressent aussi au nombre de 
composantes irr\'eductibles d'une courbe invariante par 
une transformation birationnelle sur une surface rationnelle
$S.$ Ils montrent que, sauf dans un cas bien particulier, 
ce nombre est born\'e par une quantit\'e qui ne 
d\'epend que de $\mathrm{S}.$ 

\begin{thm}[\cite{DJS}]
{\sl Soient $\mathrm{S}$ une surface rationnelle et $f$ une transformation 
birationnelle sur $\mathrm{S}$ alg\'ebriquement stable et entropique. Soit 
$\mathcal{C}$ une courbe de $\mathrm{S}$ invariante par~$f.$ 

\smallskip

\noindent Si l'une des composantes connexes de $\mathcal{C}$
est de genre $1$ le nombre de composantes irr\'eductibles de
$\mathcal{C}$ est born\'e par $\dim \mathrm{H}^{1,1}(\mathrm{S})+2.$

\smallskip

\noindent Si toutes les composantes connexes de $\mathcal{C}$
sont de genre $0,$ alors
\begin{itemize}
\item ou bien $\mathcal{C}$ a au plus $\dim \mathrm{H}^{1,1}
(\mathrm{S})+1$ composantes irr\'eductibles;

\item ou bien il existe une application holomorphe $\pi
\colon\mathrm{S}\to\mathbb{P}^1(\mathbb{C}),$ unique modulo automorphisme
de $\mathbb{P}^1(\mathbb{C}),$ telle que $\mathcal{C}$ contienne exactement
$k\geq 2$ fibres distinctes de $\pi$ et $\mathcal{C}$ 
compte au plus $\dim \mathrm{H}^{1,1}(\mathrm{S})+k-1$ composantes 
irr\'eductibles. 
\end{itemize}}
\end{thm}

\bigskip

\section{Sous-groupes de type fini, sous-groupes nilpotents}\label{typefini}

\bigskip

\noindent Rappelons une conjecture \'enonc\'ee dans le programme de Zimmer
(\cite{Zi1}): un r\'eseau d'un groupe de Lie 
r\'eel simple connexe $\mathrm{G}$ ne peut pas agir fid\`element sur une vari\'et\'e
compacte dont la dimension serait strictement inf\'erieure au rang r\'eel de 
$\mathrm{G}.$ Motiv\'e par cette conjecture on d\'ecrit, \`a l'aide du th\'eor\`eme \ref{jefcha}, les repr\'esentations de certains r\'eseaux dans le groupe de 
Cremona.

\begin{thm}[\cite{De3}]\label{plongement}
{\sl Soient $\Gamma$ un sous-groupe
d'indice fini de $\mathrm{SL}_3(\mathbb{Z})$ et $\iota$ un morphisme injectif de $\Gamma$
dans $\mathrm{Bir}(\mathbb{P}^2).$ Alors $\iota$ est, \`a conjugaison 
birationnelle pr\`es, le plongement canonique ou la contragr\'ediente ($u\mapsto\transp u^{-1}$).}
\end{thm}

\noindent On en d\'eduit une g\'en\'eralisation du th\'eor\`eme 
\ref{autaut}.

\begin{cor}[\cite{De4}]
{\sl Soit $\phi$ un endomorphisme 
non trivial de $\mathrm{Bir}(\mathbb{P}^2).$ Il existe $\kappa$ un plongement du 
corps $\mathbb{C}$ dans lui-m\^eme et $\psi$ une transformation birationnelle tels que 
\begin{align*}
& \phi(f)=\kappa(\psi f\psi^{-1}), &&\forall\hspace{1mm} f\in\mathrm{Bir}(\mathbb{P}^2).
\end{align*}}
\end{cor}

\noindent Une autre cons\'equence est qu'une large classe de groupes ne se 
plonge pas dans le groupe de Cremona; on a par exemple l'\'enonc\'e
suivant.

\begin{cor}[\cite{De3}]
{\sl Soit $\Gamma$ un sous-groupe
d'indice fini de $\mathrm{SL}_n(\mathbb{Z}).$ D\`es que $n$ est sup\'erieur ou \'egal
\`a $4,$ le groupe $\Gamma$ ne se plonge pas dans $\mathrm{Bir}(\mathbb{P}^2).$}
\end{cor}

\noindent Avant de signaler un r\'esultat d\^u \`a Cantat qui traite une plus large 
classe de groupes, rappelons qu'un groupe $\mathrm{G}$ a la \textsl{propri\'et\'e (T) 
de Kazhdan} si toute action continue de $\mathrm{G}$ sur un espace 
de Hilbert par d\'eplacement unitaire a un point fixe global. Par exemple
$\mathrm{SL}_n(\mathbb{Z})$ poss\`ede la propri\'et\'e~(T) de Kazhdan d\`es
que $n\geq 3.$

\begin{thm}[\cite{Ca}]
{\sl Soit $\Gamma$ un groupe infini non d\'enombrable de transformations
de Cremona. Si $\Gamma$ a la propri\'et\'e (T), il existe une application 
birationnelle qui conjugue $\Gamma$ \`a un sous-groupe de $\mathrm{Aut}
(\mathbb{P}^2(\mathbb{C})).$}
\end{thm}

\noindent Dans le m\^eme esprit on s'est int\'eress\'e aux sous-groupes 
nilpotents\footnote{Soit $\mathrm{G}$ un groupe; posons $\mathrm{G}^{(0)}=\mathrm{G}$
et $\mathrm{G}^{(k)}=[\mathrm{G},\mathrm{G}^{(k-1)}]=\langle aba^{-1}b^{-1}\,\vert
\, a\in\mathrm{G},\, b\in\mathrm{G}^{(k-1)}\rangle$ pour $k\geq 1.$ Le groupe $\mathrm{G}$
est nilpotent de longueur $k$ si $\mathrm{G}^{(k)}=\mathrm{id}.$} du groupe de Cremona.

\begin{thm}[\cite{De5}]
{\sl Soit $\mathrm{G}$ un groupe nilpotent; supposons que $\mathrm{G}$
ne soit pas ab\'elien \`a indice fini pr\`es. 
Soit $\iota$ un morphisme injectif de $\mathrm{G}$ dans $\mathrm{Bir}
(\mathbb{P}^2).$ Le groupe $\mathrm{G}$ v\'erifie une 
des propri\'et\'es suivantes
\smallskip
\begin{itemize}
\item $\mathrm{G}$ est de torsion;

\item le premier groupe d\'eriv\'e de $\mathrm{G}$ est ab\'elien \`a indice fini pr\`es.
\end{itemize} }
\end{thm}

\noindent On en d\'eduit des obstructions \`a ce que certains groupes
se plongent dans $\mathrm{Bir}(\mathbb{P}^2).$

\begin{cor}[\cite{De5}]
{\sl Soit $\mathrm{G}$ un groupe contenant un sous-groupe nilpotent, non 
virtuellement m\'etab\'elien et sans torsion. Il n'existe pas de repr\'esentation
fid\`ele de $\mathrm{G}$ dans le groupe de Cremona.}
\end{cor}

\bigskip

\subsection{Id\'ee de la d\'emonstration du th\'eor\`eme \ref{plongement}}
Pour simplifier, on va se restreindre au cas~$\Gamma=
\mathrm{SL}_3(\mathbb{Z}),$ {\it i.e.} on consid\`ere un plongement 
$\iota\colon\mathrm{SL}_3(\mathbb{Z})\hookrightarrow\mathrm{Bir}
(\mathbb{P}^2).$ 

\bigskip

\subsubsection{Quelques notions sur les groupes}

\bigskip

\noindent Notons $\delta_{ij}$ la matrice de Kronecker de taille $3\times
3$ et~$\mathrm{e}_{ij}=\mathrm{id}+\delta_{ij}.$

\begin{pro}[\cite{St}] \label{presensl3}
{\sl Le groupe $\mathrm{SL}_3(\mathbb{Z})$ a pour pr\'esentation 
\begin{align*}
\langle
\mathrm{e}_{ij},\hspace{1mm} i\not=j\hspace{1mm}\vert\hspace{1mm}
(\mathrm{e}_{12}\mathrm{e}_{21}^{-1}\mathrm{e}_{12})^4=\mathrm{id},\hspace{1mm} 
[\mathrm{e}_{ij},\mathrm{e}_{kl}]=\left\{\begin{array}{lll}
\mathrm{id} \text{ si }i\not=l \text{ et } j\not=k\\
\mathrm{e}_{il} \text{ si } i\not=l \text{ et } j=k\\
\mathrm{e}_{kj}^{-1} \text{ si } i=l \text{ et } j\not=k\end{array}\right.
\rangle .
\end{align*}}
\end{pro}

\noindent Les $\mathrm{e}_{ij}$ seront dits
\textsl{g\'en\'erateurs standards} de $\mathrm{SL}_3(
\mathbb{Z}).$

\begin{defi}
On appelle $k$-\textsl{groupe de Heisenberg} tout groupe
de pr\'esentation
\begin{align*}
&\mathcal{H}_k=\langle \mathrm{f},\hspace{1mm}\mathrm{g},
\hspace{1mm}\mathrm{h}\hspace{1mm}\vert\hspace{1mm}[\mathrm{f},
\mathrm{h}]=[\mathrm{g},\mathrm{h}]=
\mathrm{id},\hspace{1mm}[\mathrm{f},\mathrm{g}]=\mathrm{h}^k\rangle.
\end{align*}
\end{defi}

\noindent Par convention $\mathcal{H}_1=\mathcal{H}$ est un 
groupe de Heisenberg.

\noindent Le groupe $\mathrm{SL}_3(\mathbb{Z})$ contient de nombreux
groupes de Heisenberg~; par exemple le sous-grou\-pe~$\langle \mathrm{e}_{12},\hspace{1mm} \mathrm{e}_{13},\hspace{1mm} 
\mathrm{e}_{23}\rangle$ en
est un. Plus g\'en\'eralement le groupe engendr\'e par trois 
g\'en\'erateurs \og qui se suivent sur le dessin suivant\fg
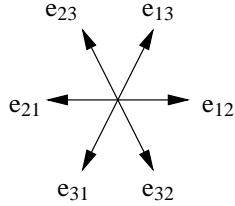
\begin{figure}[H]
\begin{center}
\input{sl.pstex_t}
\end{center}
\end{figure}
\noindent est un groupe de Heisenberg. Remarquons 
que sur la figure qui pr\'ec\`ede on a
\begin{itemize}
\item \og $\mathrm{e}_{12}+\mathrm{e}_{23}=\mathrm{e}_{13}$ \fg\hspace{1mm} ce qui traduit que
$[\mathrm{e}_{12},\mathrm{e}_{23}]=\mathrm{e}_{13};$

\item \og $\mathrm{e}_{12}+\mathrm{e}_{13}$ \fg\hspace{1mm} n'appartient pas au 
diagramme ce qui signifie que $[\mathrm{e}_{12},\mathrm{e}_{13}]=\mathrm{id}.$
\end{itemize}

\begin{defi}
Soient $\mathrm{G}$ un groupe de type
fini, $\{a_1$, $\ldots$, $a_n\}$ une partie g\'en\'eratrice de $\mathrm{G}$ 
et~$f$ un \'el\'ement de $\mathrm{G}.$ 

\noindent La \textsl{longueur} de $f,$ not\'ee
$\vert f\vert,$ est le plus petit entier $k$ pour lequel il existe une
suite $(s_1,\ldots,s_k)$ d'\'el\'ements de $\{a_1,\ldots,a_n,
a_1^{-1}, \ldots,a_n^{-1}\}$ telle que $f=s_1\ldots s_k.$ 

\noindent Un
\'el\'ement $f$ de $\mathrm{G}$ est \textsl{distordu} s'il est
d'ordre infini et si la quantit\'e $\displaystyle\lim_{k\to+\infty}\frac{
\vert f^k\vert}{k}$ est nulle. 
\end{defi}

\begin{eg}
Le g\'en\'erateur $\mathrm{h}$ d'un groupe de Heisenberg est distordu. 
\`A partir des \'egali\-t\'es~$[\mathrm{f},\mathrm{h}]=[\mathrm{g},\mathrm{h}]=
\mathrm{id}$ et $[\mathrm{f},\mathrm{g}]=\mathrm{h}$ on obtient 
\begin{align*}
& [\mathrm{f}^p,\mathrm{g}^q]=\mathrm{h}^{pq} &&\forall\hspace{1mm}
p,\hspace{1mm} q\in \mathbb{N}.
\end{align*}
\noindent Il en r\'esulte que $\vert\mathrm{h}^{p^2}\vert=\vert[\mathrm{f}^p,\mathrm{g}^p]\vert\leq 4p$
d'o\`u l'in\'egalit\'e $\frac{\vert\mathrm{h}^{p^2}\vert}{p^2}\leq \frac{4}{p}.$
\end{eg}

\medskip

\noindent En particulier les g\'en\'erateurs standards de $\mathrm{SL}_3
(\mathbb{Z})$ sont distordus.

\bigskip

\subsubsection{Dynamique de l'image d'un groupe de Heisenberg}

\bigskip

\begin{lem}[\cite{De3}]
{\sl Soient $f$ un \'el\'ement d'un groupe de type fini $\mathrm{G}$ et
$\phi$ un morphisme de $\mathrm{G}$ dans~$\mathrm{Bir}(\mathbb{P}^2)$. Si $f$ est
distordu, le premier degr\'e dynamique de $\phi(f)$ vaut $1$.}
\end{lem}

\begin{proof}[{\sl D\'emonstration}]
Soit $\{a_1$,
$\ldots$, $a_n\}$ une partie g\'en\'eratrice de $\mathrm{G}.$ Les
in\'egalit\'es 
\begin{align*}
&\lambda(\phi(f))^n\leq\deg \phi(f)^n\leq\max_i(\deg\phi(a_i))^{ \vert f^n\vert}
\end{align*}
\noindent conduisent \`a 
\begin{align*}
0\leq\log\lambda(\phi(f))\leq\frac{\vert f^n\vert}{n}\log(\max_i(\deg
\phi(a_i))).
\end{align*}
\noindent Si $f$ est distordu, la quantit\'e $\displaystyle\lim_{k
\to +\infty}\frac{\vert f^k\vert}{k}$ est nulle et le degr\'e dynamique de
$\phi(f)$ vaut un.
\end{proof}

\noindent D'apr\`es le th\'eor\`eme \ref{jefcha} on a 
\smallskip
\begin{itemize}
\item ou bien $\iota(\mathrm{e}_{ij})$ est un twist de de Jonqui\`eres
ou d'Halphen donc en particulier laisse une 
unique fibration (rationnelle ou elliptique) invariante;

\item ou bien $\iota(\mathrm{e}_{ij})$ est virtuellement isotope \`a l'identit\'e.
\end{itemize}
\smallskip

\noindent Les relations satisfaites par les $\mathrm{e}_{ij}$ assurent que si 
$\iota(\mathrm{e}_{i_0j_0})$ pr\'eserve une unique fibration alors tous les 
$\iota(\mathrm{e}_{ij})$ la pr\'eservent. D'o\`u l'alternative
\smallskip
\begin{itemize}
\item l'un des $\iota(\mathrm{e}_{ij})$  laisse une unique fibration 
(rationnelle ou elliptique) invariante;

\item tous les $\iota(\mathrm{e}_{ij})$ sont virtuellement isotopes \`a l'identit\'e. 
\end{itemize}

\bigskip

\subsubsection{Fibration invariante}

\bigskip

\noindent Rappelons une des propri\'et\'es satisfaites par les groupes
de type fini ayant la propri\'et\'e (T) qui nous sera utile.

\begin{lem}[\cite{De3}]\label{kaz} 
{\sl Soient $\Gamma$ un groupe de type fini ayant la propri\'et\'e (T)
et $\phi$ un morphisme de $\Gamma$ dans $\mathrm{PGL}_2(\mathbb{C}(y))$ (resp.
$\mathrm{PGL}_2(\mathbb{C})$). L'image de $\phi$ est finie.}
\end{lem}

\begin{pro}[\cite{De3}] 
{\sl Soit $\phi$ un morphisme de $\mathrm{SL}_3(\mathbb{Z})$ dans 
$\mathrm{Bir}(\mathbb{P}^2)$. Si 
l'un des~$\phi(\mathrm{e}_{ij})$ pr\'eserve une unique fibration, l'image 
de $\phi$ est finie.}
\end{pro}

\begin{proof}[{\sl D\'emonstration}]
Notons $\tilde{\mathrm{e}}_{ij}$ l'image de $\mathrm{e}_{ij}$ par $\phi;$
les diff\'erents g\'en\'erateurs
jouent un r\^ole identique, on peut donc supposer, sans
perdre de g\'en\'eralit\'e, que $\tilde{\mathrm{e}}_{12}$ pr\'eserve une
unique fibration~$\mathcal{F}.$

\noindent Les relations satisfaites par les $\mathrm{e}_{ij}$ (proposition \ref{presensl3}) 
entra\^inent que $\mathcal{F}$ est invariante par tous les~$\tilde{\mathrm{e}}_{ij}$. 
\noindent Ainsi, pour tout $\tilde{\mathrm{e}}_{ij},$ il existe
$\nu_{ij}$ dans~$\mathrm{PGL}_2(\mathbb{C})$ et $F\colon
\mathbb{P}^2(\mathbb{C})\to\mathrm{Aut}(\mathbb{P}^1
(\mathbb{C}))$ d\'efinissant $\mathcal{F}$ tels que $F\circ \tilde{\mathrm{e}}_{i
j}=\nu_{ij}\circ F.$ Soit $\varsigma$ le morphisme d\'efini
par
\begin{align*}
&\varsigma\colon\mathrm{SL}_3(\mathbb{Z})\to\mathrm{PGL}_2(\mathbb{C}),
&&\mathrm{e}_{ij}\mapsto \nu_{ij}.
\end{align*}
\noindent Puisque $\mathrm{SL}_3(\mathbb{Z})$ a la propri\'et\'e (T) de Kazhdan,
le lemme \ref{kaz} assure que $\mathrm{im}\,\varsigma$ est
finie, {\it i.e.} $\mathrm{ker}\,\varsigma$ est d'indice fini donc 
en particulier poss\`ede la propri\'et\'e (T). 

\noindent Si
$\mathcal{F}$ est rationnelle, on peut supposer que
l'image de $\phi$ est contenue dans le groupe de de 
Jonqui\`eres; par suite la restriction de $\phi$
\`a $\mathrm{ker}\,\varsigma$ est \og \`a valeurs\fg\hspace{1mm} dans~$\mathrm{PGL}_2(\mathbb{C}(y)),$ elle
ne peut donc \^etre injective. Dans ce cas $\phi(\mathrm{ker}\,\varsigma)$ est
fini ce qui implique que $\mathrm{im}\,\phi$ l'est aussi.

\noindent La fibration $\mathcal{F}$ ne peut \^etre
elliptique~; en effet le groupe des transformations birationnelles
qui pr\'eservent une fibration elliptique fibre \`a fibre est
m\'etab\'elien et un sous-groupe d'indice fini de~$\mathrm{SL}_3(\mathbb{Z})$
ne peut pas l'\^etre.
\end{proof}

\bigskip

\subsubsection{Factorisation dans un groupe d'automorphismes}

\bigskip

\noindent \'Etudions le cas o\`u tout $\iota(\mathrm{e}_{ij})$  
est virtuellement isotope \`a l'identit\'e. Avec des techniques similaires \`a celles utilis\'ees
pour \'etablir le lemme \ref{commut} on montre 
l'\'enonc\'e suivant.

\begin{pro}\label{pui}
{\sl Soit $\varsigma$ une repr\'esentation de $\mathcal{H}_k$ dans le groupe
de Cremona. Supposons que tout g\'en\'erateur standard de
$\varsigma(\mathcal{H}_k)$ soit virtuellement isotope \`a l'identit\'e.
Il existe une surface~$\tilde{\mathrm{S}},$
une transformation birationnelle $\phi\colon \mathrm{S}\dashrightarrow
\tilde{\mathrm{S}}$ et un entier $k$ tels que 
\begin{align*}
& \phi \varsigma(\ell)\phi^{-1}\in\mathrm{Aut}(\tilde{\mathrm{S}})&& \text{et} &&
\phi \varsigma(\ell)^k\phi^{-1}\in\mathrm{Aut}^0(\tilde{\mathrm{S}}), && \ell\in\{\mathrm{f},\,
\mathrm{g},\,\mathrm{h}\}.
\end{align*}}
\end{pro}

\noindent Ainsi les images de 
$\mathrm{e}_{12}^n,$ $\mathrm{e}_{13}^n$ et $\mathrm{e}_{2 3}^n$ par
$\iota$ sont, pour un certain $n,$ des automorphismes d'une m\^eme
surface $\mathrm{S};$ de plus quitte \`a changer $n$ on peut supposer que 
$\iota(\mathrm{e}_{12}^n),$ $\iota(\mathrm{e}_{13}^n)$ et $\iota(
\mathrm{e}_{23}^n)$ appartiennent \`a~$\mathrm{Aut}^0 (\mathrm{S})$ (proposition
\ref{pui}).

\noindent Comme un automorphisme isotope \`a l'identit\'e 
pr\'eserve les courbes d'auto-intersection n\'egative on peut supposer
\`a conjugaison birationnelle pr\`es et indice fini pr\`es que $\mathrm{S}$ est
minimale, {\it i.e.}~$\mathrm{S}=~\mathbb{P}^2(\mathbb{C})$ ou $\mathbb{P}^1(
\mathbb{C})\times \mathbb{P}^1(\mathbb{C})$ ou 
$\mathrm{F}_m.$ Rappelons que 
\begin{align*}
& \mathrm{Aut}(\mathbb{P}^1(\mathbb{C})\times\mathbb{P}^1(\mathbb{C}))
=(\mathrm{PGL}_2(\mathbb{C})\times\mathrm{PGL}_2(\mathbb{C}))\rtimes
(y,x), && \mathrm{Aut}(\mathbb{P}^2(\mathbb{C}))=\mathrm{PGL}_3
(\mathbb{C})
\end{align*}

\noindent et pour $m\geq 1$

\begin{small}
\begin{align*}
&\mathrm{Aut}(\mathrm{F}_m)=\left\{\left(\frac{\alpha x+P(y)}{(cy+
d)^m},\frac{ay+b}{cy+d}\right)\hspace{0.1cm}\Big\vert\hspace{0.1cm}\left[\begin{array}{cc}a
& b \\ c & d\end{array}\right]\in\mathrm{PGL}_2(\mathbb{C}),\hspace{0.1cm}\alpha
\in\mathbb{C}^*,\hspace{0.1cm} P\in\mathbb{C}[y],\hspace{0.1cm}\deg P\leq m\right\}.
\end{align*}
\end{small}

\noindent Il n'y a pas de plongement de $\langle \mathrm{e}_{12}^n,\hspace{1mm}
\mathrm{e}_{13}^n,\hspace{1mm} \mathrm{e}_{23}^n\rangle$ dans 
$\mathrm{Aut}(\mathbb{P}^1(\mathbb{C})\times \mathbb{P}^1(\mathbb{C}))$
donc $\mathrm{S}=\mathbb{P}^2(\mathbb{C})$ ou~$\mathrm{F}_m.$ 

\noindent Si $\mathrm{S}=\mathrm{F}_m$ on peut montrer que l'inclusion 
$\iota(\langle \mathrm{e}_{12}^n,\hspace{1mm}\mathrm{e}_{13}^n,\hspace{1mm}
\mathrm{e}_{23}^n\rangle)\subset\mathrm{Aut}(\mathrm{F}_m)$ entra\^ine 
que 
\begin{align*}
\iota(\langle \mathrm{e}_{ij}^n,\hspace{1mm} i\not=j\hspace{1mm}\rangle)
\subset\mathrm{Aut}[\mathbb{C}^2].
\end{align*}
\noindent Alors \cite{Ca-La} implique que $\iota(\langle \mathrm{e}_{ij}^n,\, i
\not=j\,\rangle)$ est n\'ecessairement un sous-groupe de $\mathrm{PGL}_3(\mathbb{C}).$

\noindent Enfin lorsque $\mathrm{S}=\mathbb{P}^2(\mathbb{C})$ on obtient que 
si l'image par $\iota$ de $\langle\mathrm{e}_{12}^n,\,
\mathrm{e}_{13}^n,\,\mathrm{e}_{23}^n\rangle$ est un sous-groupe
d'automorphismes de $\mathbb{P}^2(\mathbb{C}),$ alors
$\iota(\langle \mathrm{e}_{ij}^n,\, i\not=j\,\rangle)$
l'est aussi.

\bigskip

\subsection{Conclusion.}

\bigskip

\noindent D'apr\`es ce qui pr\'ec\`ede l'image de tout
g\'en\'erateur standard de $\mathrm{SL}_3(\mathbb{Z})$ est 
virtuellement isotope \`a l'identit\'e et $\iota(\mathrm{e}_{12}^n)$,
$\iota(\mathrm{e}_{13}^n)$ et $\iota(\mathrm{e}_{23}^n)$ sont, pour un certain
$n,$ conjugu\'es \`a des automorphismes d'une surface minimale
$\mathrm{S}$ avec $\mathrm{S}=\mathbb{P}^2(\mathbb{C})$ ou $\mathrm{S}=\mathrm{F}_m$, $m\geq 1.$ 
On montre finalement qu'\`a conjugaison pr\`es
$\iota(\langle\mathrm{e}_{ij}^n\rangle)$ est, pour un certain $n,$ 
un sous-groupe de~$\mathrm{PGL}_3(\mathbb{C}).$ La restriction de $\iota$ \`a $\langle\mathrm{e}_{ij}^n\rangle$ se
prolonge alors en un morphisme de groupe de Lie de 
$\mathrm{PGL}_3(\mathbb{C})$ dans lui-m\^eme (\cite{St})~; 
par simplicit\'e de 
$\mathrm{PGL}_3(\mathbb{C})$, ce prolongement est injectif et donc surjectif. Or
les automorphismes lisses
de~$\mathrm{PGL}_3(\mathbb{C})$ s'obtiennent \`a partir des automorphismes
int\'erieurs et de la contragr\'e\-diente (th\'eor\`eme \ref{dieudon})~; ainsi, \`a conjugaison
lin\'eaire pr\`es, $\iota_{\vert\langle\mathrm{e}_{ij}^n\rangle}$ 
co\"incide avec le plongement canonique ou la
contragr\'e\-diente. 

\noindent Soit $f$ un \'el\'ement de $\iota(\mathrm{SL}_3(\mathbb{Z}))\setminus \iota(
\langle\mathrm{e}_{ij}^n\rangle)$ qui contracte au moins une courbe~$\mathcal{
C}=\mathrm{Exc}(f)$. Le groupe $\langle\mathrm{e}_{ij}^n\rangle$ est distingu\'e dans
$\mathrm{SL}_3(\mathbb{Z});$ la courbe $\mathcal{C}$ est donc invariante par 
$\iota(\langle\mathrm{e}_{ij}^n\rangle)$ et par suite l'est aussi par l'adh\'erence 
de Zariski de $\iota(\langle\mathrm{e}_{ij}^n\rangle),$ qui n'est autre que 
$\mathrm{PGL}_3(\mathbb{C}),$ ce qui est impossible. Il en r\'esulte que $f$ appartient \`a $\mathrm{PGL}_3
(\mathbb{C})$ et que $\iota(\mathrm{SL}_3(\mathbb{Z}))$ est inclus dans $\mathrm{PGL}_3(\mathbb{C}).$

\section{Alternative de Tits}\label{alttits}

\noindent Le groupe lin\'eaire satisfait l'alternative de Tits. 

\begin{thm}[\cite{Ti}]\label{tits}
{\sl Soient $\Bbbk$ un corps de caract\'eristique nulle et $\Gamma$ un sous-groupe 
de type fini de $\mathrm{GL}_n(\Bbbk).$ Alors
\begin{itemize}
\item ou bien $\Gamma$ contient un groupe libre non ab\'elien;

\item ou bien $\Gamma$ contient un sous-groupe r\'esoluble\footnote{Soit 
$\mathrm{G}$ un groupe; posons $\mathrm{G}^{(0)}=\mathrm{G}$
et $\mathrm{G}^{(k)}=[\mathrm{G}^{(k-1)},\mathrm{G}^{(k-1)}]=\langle 
aba^{-1}b^{-1}\,\vert\, a,\, b\in \mathrm{G}^{(k-1)}\rangle$ pour $k\geq 1.$ Le groupe $\mathrm{G}$
est r\'esoluble s'il existe un entier $k$ tel que $\mathrm{G}^{(k)}=\{\mathrm{id}\}.$} d'indice fini.
\end{itemize}}
\end{thm}

\noindent Signalons que le groupe des diff\'eomorphismes d'une 
vari\'et\'e r\'eelle de dimension sup\'erieure ou \'egale \`a $1$ ne 
satisfait pas l'alternative de Tits (\emph{voir} \cite{Gh} et les r\'ef\'erences
qui s'y trouvent). Par contre le groupe des
automorphismes polynomiaux de $\mathbb{C}^2$ v\'erifie  
l'alternative de Tits~(\cite{La}); Lamy obtient cet \'enonc\'e \`a 
partir de la classification des sous-groupes de  $\mathrm{Aut}[\mathbb{C}^2],$
classification \'etablie en utilisant l'action de ce groupe sur $\mathcal{T}.$

\begin{thm}[\cite{La}]\label{autsteph}
{\sl Soit $\mathrm{G}$ un sous-groupe de $\mathrm{Aut}[\mathbb{C}^2].$ Une 
et une seule des \'eventualit\'es suivantes est r\'ealis\'ee.
\begin{itemize}
\item Chaque \'el\'ement de $\mathrm{G}$ est \'el\'ementaire, on a 
alors l'alternative
\begin{itemize}
\item $\mathrm{G}$ est conjugu\'e \`a un sous-groupe de $\mathtt{E}$
ou $\mathtt{A};$

\item $\mathrm{G}$ est ab\'elien, $\mathrm{G}$ s'\'ecrit $\bigcup_{i\in\mathbb{N}}\mathrm{G}_i$
o\`u $\mathrm{G}_i\subset\mathrm{G}_{i+1}$ et 
chaque $\mathrm{G}_i$ est conjugu\'e \`a un groupe cyclique
fini du type $\langle(\alpha x,\beta y)\rangle$ avec $\alpha,$ $\beta$
racines de l'unit\'e du m\^eme ordre. Chaque \'el\'ement
de $\mathrm{G}$ admet un unique point fixe\footnote{en tant qu'automorphisme polynomial de $\mathbb{C}^2$} et ce 
point fixe est le m\^eme pour tous les \'el\'ement de $\mathrm{G}.$ 
Enfin l'action de $\mathrm{G}$ fixe un bout de l'arbre $\mathcal{T}.$
\end{itemize}

\item $\mathrm{G}$ contient des \'el\'ements de type H\'enon et ceux-ci
admettent tous la m\^eme g\'eod\'esique auquel cas $\mathrm{G}$ est 
r\'esoluble et contient un sous-groupe d'indice fini isomorphe \`a $\mathbb{Z}.$

\item $\mathrm{G}$ poss\`ede deux \'el\'ements de type H\'enon avec
des g\'eod\'esiques distinctes, $\mathrm{G}$ contient alors un groupe 
libre \`a deux g\'en\'erateurs.
\end{itemize}}
\end{thm}

\noindent En s'appuyant, entre autres, sur le th\'eor\`eme \ref{jefcha}, Cantat 
caract\'erise les sous-groupes de type fini de $\mathrm{Bir}(\mathbb{P}^2);$
Favre a reformul\'e, dans \cite{Fa}, cette classification de la fa\c{c}on 
suivante. 

\begin{thm}[\cite{Ca}]\label{clastypefini}
{\sl Soit $\mathrm{G}$ un sous-groupe de type fini du groupe de Cremona. 
Une et une seule des possibilit\'es suivantes est r\'ealis\'ee.
\begin{itemize}
\item Tous les \'el\'ements de $\mathrm{G}$ sont virtuellement isotopes 
\`a l'identit\'e, alors
\begin{itemize}
\item ou bien $\mathrm{G}$ est, \`a indice fini pr\`es et conjugaison birationnelle pr\`es, 
contenu dans la composante neutre de $\mathrm{Aut}(\mathrm{S})$ o\`u $\mathrm{S}$ 
d\'esigne une surface rationnelle minimale;

\item ou bien $\mathrm{G}$ pr\'eserve une fibration rationnelle.
\end{itemize}

\item $\mathrm{G}$ contient un twist de Jonqui\`eres ou 
d'Halphen et aucun 
\'el\'ement entropique auquel cas $\mathrm{G}$ pr\'eserve une fibration
rationnelle ou elliptique.

\item $\mathrm{G}$ contient deux \'el\'ements entropiques $f$ et $g$ 
tels que $\langle f,\hspace{0.1cm} g\rangle$ soit libre.

\item $\mathrm{G}$ contient un \'el\'ement entropique et pour tout 
couple $(f,g)$ d'\'el\'ements entropiques, $\langle f,\hspace{0.1cm} g
\rangle$ n'est pas libre, alors
\begin{align*}
& 1\longrightarrow\ker\rho\longrightarrow\mathrm{G}\stackrel{\rho}{\longrightarrow}
\mathbb{Z}\longrightarrow 1
\end{align*}
\noindent et $\ker\rho$ est constitu\'e uniquement d'\'el\'ements virtuellement isotopes
\`a l'identit\'e.
\end{itemize}}
\end{thm}

\noindent Une des cons\'equences est l'\'enonc\'e suivant.

\begin{thm}[\cite{Ca}] {\sl Le groupe $\mathrm{Bir}(\mathbb{P}^2)$ 
satisfait l'alternative de Tits.} 
\end{thm}

\noindent Un des ingr\'edients communs aux d\'emonstrations des th\'eor\`emes
\ref{tits}, \ref{autsteph} \ref{clastypefini} est le \og lemme du ping-pong\fg,
crit\`ere utilis\'e de nombreuses fois par 
Klein pour construire des produits libres (la formulation
donn\'ee ci-apr\`es est toutefois plus r\'ecente). 

\begin{lem}[\og Lemme du ping-pong\fg]\label{pingpong}
{\sl Soit $\mathrm{G}$ un groupe agissant sur un ensemble
$\mathrm{X}.$ Consid\'erons $\Gamma_1$ et $\Gamma_2$ deux sous-groupes de $\mathrm{G}$ et
posons $\Gamma=\langle\Gamma_1,\Gamma_2\rangle.$ On suppose que
\begin{itemize}
\item $\Gamma_1$ (resp. $\Gamma_2$) compte au moins $3$ (resp. $2$) \'el\'ements,

\item il existe $X_1$ et $X_2$ deux parties non vides de $X$
telles que
\begin{align*}
&\mathrm{X}_2 \nsubseteq\mathrm{X}_1;
&&\forall\, \alpha \in \Gamma_1 \setminus \{\mathrm{id}\}, \hspace{3mm}
\alpha(\mathrm{X}_2) \subset\mathrm{X}_1;
&&\forall\, \beta \in \Gamma_2 \setminus \{\mathrm{id}\}, \hspace{3mm}
\beta(\mathrm{X}_1) \subset\mathrm{X}_2.
\end{align*}
\end{itemize}
\noindent Alors $\Gamma$ est isomorphe au produit libre $\Gamma_1 \ast \Gamma_2$ de $\Gamma_1$ et
$\Gamma_2.$}
\end{lem}

\begin{proof}[{\sl D\'emonstration}]
Consid\'erons le morphisme surjectif $\Gamma_1 \ast \Gamma_2
\twoheadrightarrow \Gamma.$ Soit $m$ un mot r\'eduit, non vide,
compos\'e de lettres dans l'union disjointe de $\Gamma_1 \setminus
\{\mathrm{id}\} \cup \Gamma_2 \setminus \{\mathrm{id}\}$
\begin{align*}
&m=(\alpha_1)\beta_1\alpha_2\beta_2\ldots \beta_k (\alpha_{k+1}).
\end{align*}

\noindent Commen\c{c}ons par supposer que $m$ soit de la forme
\begin{align*}
\alpha_1\beta_1\alpha_2\beta_2\ldots \beta_k \alpha_{k+1}.
\end{align*}
Comme $\alpha_i(\mathrm{X}_2)\subset \mathrm{X}_1$ et $\beta_i(
\mathrm{X}_1)\subset\mathrm{X}_2$, on
obtient $m(\mathrm{X}_2)\subset\mathrm{X}_1$ donc $m\not=\mathrm{id}.$

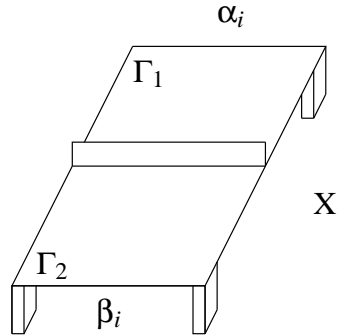
\begin{figure}[H]
\begin{center}
\input{pingpong2.pstex_t}
\end{center}
\end{figure}

\noindent Supposons maintenant que $m$ soit du type
$\beta_1\alpha_2\beta_2\ldots \beta_k.$ Soit $\alpha \in \Gamma_1
\setminus \{\mathrm{id}\},$ l'\'el\'ement $m$ est trivial si et seulement si
$\alpha m\alpha^{-1}=\mathrm{id};$ or d'apr\`es le cas
pr\'ec\'edent $\alpha m\alpha^{-1}\not=\mathrm{id}$ donc
$m\not=\mathrm{id}.$

\noindent Supposons d\'esormais que 
$m$ s'\'ecrive $\alpha_1\beta_1\alpha_2\beta_2\ldots \beta_k.$
Soit $\alpha \in \Gamma_1
\setminus\{\mathrm{id},\alpha_1\},$ le mot $m$ est trivial
si et seulement si
$\alpha^{-1}m\alpha=\mathrm{id}$ ce qui est
exclu d'apr\`es ce qui pr\'ec\`ede.

\noindent Finalement consid\'erons l'\'eventualit\'e o\`u $m=\beta_1 \alpha_2
\beta_2\ldots \alpha_k$ auquel cas $m$ est trivial si et
seulement si $m^{-1}=\alpha_k^{-1}\ldots\beta_1^{-1}$ l'est 
ce qui est impossible d'apr\`es ce qu'on vient de voir.
\end{proof}

\begin{eg}\label{egpingpong}
Les matrices $\left[%
\begin{array}{cc}
  1 & 2 \\
  0 & 1 \\
\end{array}%
\right]$ et $\left[%
\begin{array}{cc}
  1 & 0 \\
  2 & 1 \\
\end{array}%
\right]$ engendrent un sous-groupe libre de
rang~$2$ dans $\mathrm{SL}_2(\mathbb{Z})$. En effet, posons
\begin{align*}
&\Gamma_1=\left\{\left[%
\begin{array}{cc}
  1 & 2 \\
  0 & 1 \\
\end{array}%
\right]^n\,\vert \, n \in \mathbb{Z} \right\},
&&\Gamma_2=\left\{\left[%
\begin{array}{cc}
  1 & 0 \\
  2 & 1 \\
\end{array}%
\right]^n\,\vert\, n\in \mathbb{Z} \right\},
\end{align*}
\begin{align*}
&\mathrm{X}_1=\left\{(x,y)\in\mathbb{R}^2 \,\vert\, \vert x\vert>\vert y\vert\right\}
&&
\text{et} &&\mathrm{X}_2=\left\{(x,y)\in\mathbb{R}^2
\,\vert\, \vert x\vert<\vert y\vert\right\}.
\end{align*}
\noindent Consid\'erons $\gamma$ un \'el\'ement de $\Gamma_1\setminus \{\mathrm{id}\}$
et $(x,y)$ un \'el\'ement de $\mathrm{X}_2,$ on remarque alors que $\gamma(x,y)$
est de la forme $(x+my,y),$ avec $\vert m\vert \geq
2;$ par suite $\gamma(x,y)$ appartient \`a $\mathrm{X}_1$. De m\^eme si $\gamma$ appartient
\`a $\Gamma_2\setminus \{\mathrm{id}\}$ et $(x,y)$ \`a $\mathrm{X}_1,$ l'image de $(x,y)$
par $\gamma$ est dans $\mathrm{X}_2.$ D'apr\`es le lemme \ref{pingpong} on a
\begin{align*}
\langle\left[%
\begin{array}{cc}
  1 & 2 \\
  0 & 1 \\
\end{array}%
\right] ,\,\left[%
\begin{array}{cc}
  1 & 0 \\
  2 & 1 \\
\end{array}%
\right]\rangle \simeq \mathrm{F}_2= \mathbb{Z}\ast
\mathbb{Z}=\Gamma_1\ast\Gamma_2.
\end{align*}
\smallskip

\noindent De la m\^eme mani\`ere on obtient que 
\begin{align*}
&\left[%
\begin{array}{cc}
  1 & k \\
  0 & 1 \\
\end{array}%
\right] && \text{et} &&\left[%
\begin{array}{cc}
  1 & 0 \\
  k & 1 \\
\end{array}%
\right] 
\end{align*}
\noindent engendrent un groupe libre de rang $2$ dans
$\mathrm{SL}_2(\mathbb{Z})$ pour tout $k \geq 2.$ Par contre ce n'est pas 
le cas pour~$k=~1,$ les matrices 
\begin{align*}
&\left[%
\begin{array}{cc}
  1 & 1 \\
  0 & 1 \\
\end{array}%
\right]  && \text{et} &&\left[%
\begin{array}{cc}
  1 & 0 \\
  1 & 1 \\
\end{array}%
\right] 
\end{align*}
\noindent engendrent $\mathrm{SL}_2(\mathbb{Z})$.
\end{eg}

\begin{eg}
Deux matrices prises au hasard dans
$\mathrm{SL}_\nu(\mathbb{C}),$ o\`u $\nu$ d\'esigne un entier sup\'erieur ou \'egal \`a  $2,$ engendrent un groupe libre
isomorphe \`a $\mathrm{F}_2.$ En effet, soit $\mathcal{F}$ le sous-ensemble 
de~$\mathrm{SL}_\nu(\mathbb{C}) \times \mathrm{SL}_\nu(\mathbb{C})$ d\'efini par
\begin{align*}
\mathcal{F}=\{(A,B) \in \mathrm{SL}_\nu(\mathbb{C}) \times
\mathrm{SL}_\nu(\mathbb{C}) \,\vert\, \langle A,\,B\rangle \simeq \mathrm{F}_2\}.
\end{align*}
\noindent Le compl\'ementaire $\complement\mathcal{F}$ de
$\mathcal{F}$ est une r\'eunion d\'enombrable de sous-vari\'et\'es
alg\'ebriques: \`a tout mot $m=a^{n_1} b^{\ell_1}\ldots
a^{n_k}b^{\ell_k}$ r\'eduit non trivial dans $\mathrm{F}_2,$ on peut associer la
sous-vari\'et\'e $V_m$ d'\'equa\-tion~$A^{n_1} B^{\ell_1}\ldots
A^{n_k}B^{\ell_k}=\mathrm{id}.$ Le sous-ensemble $\mathcal{F}$
n'est pas vide: il contient 
\begin{align*}
(A,B)=\left(\left[%
\begin{array}{cccc}
  1 & 2 &  & \mbox{\LARGE {0}} \\
  0 & 1 &  &  \\
   &  & \ddots &  \\
  \mbox{\LARGE {0}} &  &  & 1 \\
\end{array}%
\right],\left[%
\begin{array}{cccc}
  1 & 2 &  & \mbox{\LARGE {0}} \\
  0 & 1 &  &  \\
    &   & \ddots &  \\
\mbox{\LARGE {0}} &  &  & 1 \\
\end{array}%
\right]\right).
\end{align*}
\noindent Ainsi $\complement\mathcal{F}$ est une union d\'enombrable
d'ensembles alg\'ebriques de codimension $\geq 1$ strictement
incluse dans $\mathrm{SL}_\nu(\mathbb{C}) \times \mathrm{SL}_\nu(\mathbb{C});$ 
la mesure de Lebesgue de $\complement\mathcal{F}$ est donc nulle. Autrement dit deux matrices
prises au hasard dans $\mathrm{SL}_\nu(\mathbb{C})$ engendrent un groupe
libre.
\end{eg}

\noindent La strat\'egie de Cantat pour \'etablir le th\'eor\`eme \ref{clastypefini} est la suivante. 
On peut \'etudier certaines propri\'et\'es d'une transformation birationnelle
en consid\'erant son action sur la cohomologie $\mathrm{H}^2(\mathrm{X},\mathbb{R})$
o\`u $\mathrm{X}$ d\'esigne un mod\`ele birationnel ad\'equat du plan projectif complexe.
Le choix d'un tel mod\`ele n'\'etant pas unique, Manin a introduit, dans \cite{Ma2},
l'espace $\mathbb{H}^\infty$ de toutes les classes de cohomologie de tous 
les mod\`eles birationnels de $\mathbb{P}^2(\mathbb{C}).$ Le groupe de 
Cremona agit isom\'etriquement sur cet espace de dimension infinie,
$\mathrm{Bir}(\mathbb{P}^2)$ se plonge donc dans le groupe des
isom\'etries de $\mathbb{H}^\infty.$ L'espace $\mathbb{H}^\infty$
est hyperbolique au sens de Gromov; Cantat utilise alors,
entre autres, un \og argument de ping-pong\fg.

\section{Automorphismes sur les surfaces rationnelles}\label{autrat}

\noindent Une surface compacte projective admettant un automorphisme d'entropie positive est ou bien un tore, ou bien une surface
K$3$\footnote{Une {\sl surface K$3$} est une 
surface $\mathrm{S}$ complexe, compacte, simplement connexe, \`a fibr\'e 
canonique trivial. En particulier il existe une $2$-forme holomorphe 
$\omega$ sur $\mathrm{S}$ qui ne s'annule pas; $\omega$ est unique \`a 
multiplication pr\`es par un scalaire non nul. Toute surface quartique lisse 
dans $\mathbb{P}^3(\mathbb{C})$ est une surface K$3.$
Tout rev\^etement double de $\mathbb{P}^2(\mathbb{C})$ ramifi\'e le long de 
courbes sextiques lisses est une surface K$3.$}, ou bien une surface 
rationnelle (\cite{Ca3}). Pour l'\'etude des automorphismes sur les surfaces K$3$
on renvoie \`a \cite{Ca4}. Il semble que les surfaces rationnelles poss\'edant
des automorphismes d'entropie positive soient relativement rares: les
premi\`eres familles infinies de tels automorphismes sont connues depuis peu (\cite{BK3, BK1, McM}).
N\'eanmoins c'est sur les surfaces rationnelles que les automorphismes 
d'entropie positive sont les plus abondants au sens o\`u elles poss\`edent
des familles de dimension arbitrairement grande (\cite{BK2}). Les travaux
cit\'es se concentrent sur les surfaces obtenues en \'eclatant $\mathbb{P}^2
(\mathbb{C});$ ceci est justifi\'e par le th\'eor\`eme suivant de Nagata: soient $\mathrm{S}$
une surface rationnelle et $f$ un automorphisme sur $\mathrm{S}$ tel que 
$f_*$ soit d'ordre infini. Il existe une suite finie d'\'eclatements $\pi_{j+1}\colon
\mathrm{S}_{j+1}\to \mathrm{S}_j$ telle que 
\begin{align*}
& \mathrm{S}_1=\mathbb{P}^2(\mathbb{C}), && \mathrm{S}_{N+1}=
\mathrm{S}, && \pi_{j+1} \text{ \'eclatement de }p_j\in\mathrm{S}_j. 
\end{align*}
N\'eanmoins une surface obtenue \`a partir de $\mathbb{P}^2(\mathbb{C})$ apr\`es
une suite g\'en\'erique d'\'eclatements n'a pas d'automorphisme non trivial 
(\cite{GrH, Ko}). La classification des automorphismes sur les surfaces rationnelles
est un probl\`eme ouvert.

\noindent Rappelons d'une part que pour une surface rationnelle $\mathrm{S}$ les groupes $\mathrm{H}^2(\mathrm{S},\mathbb{Z}),$
$\mathrm{H}^{1,1}(\mathrm{S})\cap\mathrm{H}^2(\mathrm{S},\mathbb{Z})$ et $\mathrm{Pic}(\mathrm{S})$ co\"incident et d'autre part que $\mathrm{Pic}(\mathbb{P}^2(\mathbb{C}))\simeq\mathbb{Z}.$
Soient $\pi\colon\mathrm{S}\to\mathbb{P}^2(\mathbb{C})$ l'\'eclatement d'un point $p$ 
de $\mathbb{P}^2(\mathbb{C})$ et $E$ le diviseur exceptionnel. Les fonctions rationnelles
sur $\mathrm{S}$ sont par d\'efinition les tir\'es en arri\`ere par~$\pi$ des
fonctions rationnelles sur $\mathbb{P}^2(\mathbb{C}).$ Par suite $\mathrm{Pic}
(\mathrm{S})$ est engendr\'e par~$E$ et le tir\'e en arri\`ere d'une droite 
$\mathrm{L}$ de $\mathbb{P}^2(\mathbb{C}).$ Plus g\'en\'eralement on a l'\'enonc\'e suivant.

\begin{thm}
{\sl Soient $p_1,$ $\ldots,$ $p_N$ des points distincts de $\mathbb{P}^2(\mathbb{C}).$ Notons
$\pi\colon\mathrm{S}\to\mathbb{P}^2(\mathbb{C})$ la suite d'\'eclatements
des $p_i.$ Si $E_j=\pi^{-1}p_j$ d\'esignent les fibres exceptionnelles et $\mathrm{L}$ une droite
g\'en\'erique de $\mathbb{P}^2(\mathbb{C}),$ alors
$\mathrm{Pic}(\mathrm{S})=\mathbb{Z}E_1\oplus\ldots\oplus\mathbb{Z}E_N\oplus\mathbb{Z}\pi^*\mathrm{L}.$}
\end{thm}

\begin{eg}\label{tau}
Consid\'erons l'involution $\tau$ introduite dans les Exemples \ref{egsptdind}
\begin{align*}
& \tau\colon\mathbb{P}^2(\mathbb{C})\dashrightarrow\mathbb{P}^2(
\mathbb{C}), && (x:y:z)\mapsto(x^2:xy:y^2-xz)
\end{align*}
\noindent Rappelons qu'elle a un seul point d'ind\'etermination $(0:0:1)$ et
une seule courbe contract\'ee, la droite $x=0$ que nous noterons $\Theta.$

\noindent Soit $\pi_1\colon\mathcal{X}_1\to\mathbb{P}^2(\mathbb{C})$ l'\'eclatement de $(0:0:1).$ On d\'esigne par $E_1$ le diviseur exceptionnel et on se place dans les coordonn\'ees locales $(\xi,s)$ avec $\pi_1(\xi,s)=(\xi s:s:1).$ On remarque que
\begin{align*}
&\Theta=\{\xi=0\}, && E_1=\{s=0\}, && \pi_1^{-1}(x:y:z)=\left(\frac{x}{y},\frac{y}{z}\right).
\end{align*}

\noindent La transformation $\tau$ est donn\'ee sur $\mathcal{X}_1$ par
\begin{align*}
&\tau_{\mathcal{X}_1}=\pi_1^{-1}\tau\pi_1\colon\mathcal{X}_1\to\mathcal{X}_1, &&
(\xi,s)\mapsto\left(\xi,\frac{\xi s}{s-\xi}\right).
\end{align*}

\noindent La restriction de $\tau_{\mathcal{X}_1}$ \`a $E_1$ est l'identit\'e; en particulier $E_1$ n'est pas exceptionnel.

\noindent Au voisinage de $\Theta$ on a
$$\pi_1^{-1}\tau(x:y:z)=\left(\frac{x}{y},\frac{xy}{y^2-xz}\right);$$
ainsi $\Theta$ est contract\'ee par $\tau_{\mathcal{X}_1}$ sur le point $(0,0)=E_1\cap \Theta.$

\noindent Notons $\pi_2\colon\mathcal{X}_2\to\mathcal{X}_1$ l'\'eclatement de $(0,0).$ On d\'esigne par $E_2$ le diviseur exceptionnel. Choisissons un syst\`eme de coordonn\'ees local $(u,v)$ dans lequel 
\begin{align*}
&\pi_2\colon\mathcal{X}_2\to\mathcal{X}_1, && (u,v)\mapsto(u,uv).
\end{align*}

\noindent On remarque que $E_2=\{u=0\}$ et $\tau_{\mathcal{X}_2}(u,v)=\left(u,\frac{v}{v-1}\right).$ Par suite $E_2$ est invariant par $\tau_{\mathcal{X}_2}$ et au voisinage de $\Theta$ on a
$$\pi_2^{-1}\pi_1^{-1}\tau(x:y:z)=\left(\frac{x}{y},\frac{y^2}{y^2-xz}\right);$$
ainsi $\Theta$ est contract\'ee sur le point $(0,1)\in E_2.$

\noindent Finalement consid\'erons l'\'eclatement $\pi_3\colon\mathcal{X}_3\to\mathcal{X}_2$ du point $(0,1).$ Utilisons les coordonn\'ees locales $(\eta,\mu)$ pour lesquelles $\pi_3(\eta,\mu)=(\eta,\eta\mu+1).$ Le diviseur exceptionnel est $E_3=\{\eta=0\}$ et au voisinage de $\Theta$ on a 
$$\pi_3^{-1}\pi_2^{-1}\pi_1^{-1}\tau=\left(\frac{x}{y},\frac{yz}{y^2-xz}\right);$$
en particulier l'image de $\Theta$ par $\tau_{\mathcal{X}_3}$ est $E_3$ et $\tau$ \'etant une involution $E_3$ est envoy\'e sur $\Theta.$ Il en r\'esulte que ni $E_3,$ ni $\Theta$ n'est exceptionnel; la transformation $\tau_{\mathcal{X}_3}$ est un automorphisme.

\bigskip

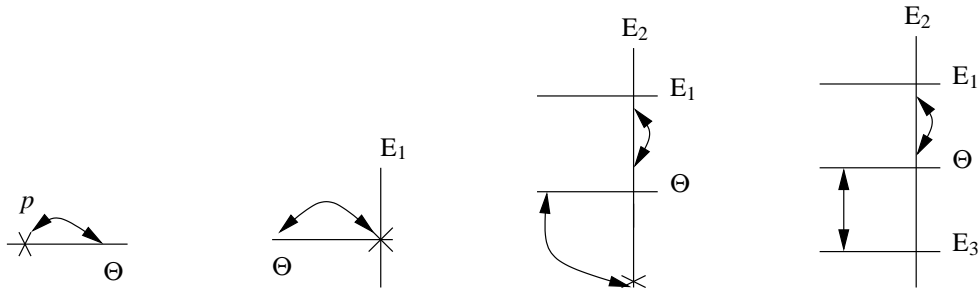
\begin{tabular}{cccc}
\hspace{0.7cm} \input{aut.pstex_t}\hspace{0.7cm} &\hspace{0.7cm} \input{aut2.pstex_t}\hspace{0.7cm}  &\hspace{0.7cm} \input{aut3.pstex_t} \hspace{0.7cm} &\hspace{0.7cm} \input{aut4.pstex_t}\hspace{0.7cm} 
\end{tabular}

\bigskip 

\noindent On va d\'eterminer l'action de $\tau^*_{\mathcal{X}_3}$ sur $\mathrm{Pic}(\mathcal{X}_3)$ dont une base est, d'apr\`es le th\'eor\`eme pr\'ec\'edent, $\{\mathrm{L},\,E_1,\,E_2,\,E_3\}.$ On a vu d\'ej\`a 
vu que
\begin{align*}
&E_1\to E_1, && E_2\to E_2, && E_3\to \Theta=\{x=0\}.
\end{align*}

\noindent On remarque que $\mathrm{L}=\Theta$ est un \'el\'ement de $\mathrm{Pic}(\mathbb{P}^2(\mathbb{C})).$ Apr\`es le premier \'eclatement on a $\mathrm{L}=\Theta+E_1\in\mathrm{Pic}(\mathcal{X}_1).$ Le centre du second \'eclatement est $\Theta\cap E_1\in\mathcal{X}_1$ ainsi, par tir\'e en arri\`ere, on obtient deux copies de $E_2$ d'o\`u $\mathrm{L}=\Theta+E_1+2E_2\in\mathrm{Pic}(\mathcal{X}_2).$ Enfin on \'eclate un point de~$E_2\setminus(\Theta\cup E_1);$ par suite $$\mathrm{L}=\Theta+E_1+2E_2+2E_3\in\mathrm{Pic}(\mathcal{X}_3).$$

\noindent Maintenant consid\'erons $\mathrm{L}=\{a_0x+a_1y+a_2z=0\};$ on constate que $\tau^*\mathrm{L}=\{a_0x^2+a_1xy+a_2(-xz+y^2)=0\}$ qui est de degr\'e $2.$ Il en r\'esulte que $\tau^*\mathrm{L}=2\mathrm{L}+\sum m_jE_j.$ D\'eterminons les $m_j.$ Pour $m_1$ on tire en arri\`ere $\ell=a_0x+a_1y+a_2z$ par $\tau\pi_1;$ on trouve $a_0\xi^2s^2+a_1\xi s^2+a_2(-\xi s+s^2)$ qui s'annule \`a l'ordre $1$ sur $E_1=\{s=0\}$ d'o\`u $m_1=1.$ Un calcul montre que $\ell\tau\pi_1\pi_2$ (resp. $\ell\tau\pi_1\pi_2\pi_3$) s'annule \`a l'ordre $2$ (resp. $3$) sur  $E_2=\{u=0\}$ (resp. $E_3$); autrement dit $m_2=2$ et $m_3=3.$ 

\noindent Il en r\'esulte que
\begin{align*}
&\mathrm{L}\to\mathrm{L}-E_1-2E_2-3E_3, && E_1\to E_1, && E_2\to E_2, 
&& E_3\to \Theta=\mathrm{L}-E_1-2E_2-2E_3;
\end{align*}
autrement dit $\tau^*_{\mathcal{X}_3}$ est donn\'e par $\left[\begin{array}{cccc}
2 & 0 & 0 & 1\\
-1 & 1 & 0 & -1 \\
-2 & 0 & 1 & -2\\
-3 & 0 & 0 & -2
\end{array}
\right].$
\end{eg}

\bigskip

\noindent Soient $\mathrm{S}$ une surface obtenue en \'eclatant $\mathbb{P}^2
(\mathbb{C})$ en un nombre fini de points et $f$ un automorphisme sur $\mathrm{S};$ on a
\begin{align*}
\mathrm{H}^*(\mathrm{S};\mathbb{C})=\mathrm{H}^0(\mathrm{S};\mathbb{C})
\oplus\mathrm{H}^{1,1}(\mathrm{S};\mathbb{C})\oplus\mathrm{H}^4
(\mathrm{S};\mathbb{C}).
\end{align*}

\noindent L'application $f^*$ agit sur chacun des facteurs.
 La dimension de 
$\mathrm{H}^0$ est le nombre de composantes connexes, ici $1,$ et 
$f^*_{\vert\mathrm{H}^0(\mathrm{S};\mathbb{C})}=\mathrm{id}.$ 
On a $\mathrm{H}^4(\mathrm{S};\mathbb{C})\simeq\mathbb{C}$ et~$f^*_{\vert
\mathrm{H}^4(\mathrm{S};\mathbb{C})}$ est la multiplication par
$\deg f=~1.$  Ainsi si $\mathrm{S}$ est obtenue en \'eclatant $\mathbb{P}^2
(\mathbb{C})$ un nombre fini de fois et si $f$ est un automorphisme de $\mathrm{S},$ alors 
\begin{align*}
\mathrm{Per}_n=2+\text{ trace }(f^{*n}_{\vert\mathrm{H}^{1,1}})
\end{align*}

\noindent o\`u $\mathrm{Per}_n$ d\'esigne le nombre de points p\'eriodiques
de p\'eriode $n$ compt\'es avec multiplicit\'e.

\bigskip

\noindent Rappelons qu'une m\'etrique  k\"{a}hl\'erienne sur
une vari\'et\'e $M$ est une m\'etrique hermitienne $h$ sur le fibr\'e
tangent v\'erifiant: si localement $h$ s'\'ecrit $\sum h_{ij}\mathrm{d}z_i\wedge \mathrm{d}
\overline{z_j},$ la forme de K\"{a}hler associ\'ee $\omega=\sum
h_{ij}\mathrm{d}z_i\wedge \mathrm{d}\overline{z_j}$ est ferm\'ee. Une {\sl vari\'et\'e 
k\"{a}hl\'erienne} est une vari\'et\'e munie 
d'une m\'etrique k\"{a}hl\'erienne. 

\begin{egs}\hspace{1mm}
\begin{itemize}
\item Puisque toute sous-vari\'et\'e d'une 
vari\'et\'e k\"{a}hl\'erienne est k\"{a}hl\'erienne et que la m\'etrique de 
Fubini-Study sur $\mathbb{P}^n(\mathbb{C})$ est k\"{a}hl\'erienne toute vari\'et\'e 
alg\'ebrique projective est k\"{a}hl\'erienne. 

\item Soit $\Lambda$ un r\'eseau dans $\mathbb{C}^n;$ le tore $\mathbb{C}^n/\Lambda$
est une vari\'et\'e k\"{a}hl\'erienne.

\item Toute surface de Riemann est une vari\'et\'e k\"{a}hl\'erienne.
\end{itemize}
\end{egs}

\begin{thm}[\cite{DiFa}]
{\sl Soit $f$ un automorphisme sur une surface k\"{a}hl\'erienne tel 
que~$\lambda(f)>1.$ Alors $\lambda(f)$ est une valeur propre de $f^*$ avec
multiplicit\'e $1$ et c'est l'unique valeur propre de module strictement
sup\'erieur \`a $1.$

\noindent Si $\eta$ est une valeur propre de $f^*$ soit $\eta$
vaut $\lambda(f)$ ou $1/\lambda(f),$ soit $\vert\eta\vert=1.$

\noindent De plus $\lambda(f)$ est un nombre de Salem\footnote{Un entier
alg\'ebrique r\'eel $\alpha > 1$, est un nombre de Salem
si ses conjugu\'es sont $\alpha,$ $\alpha^{-1}$ et $\mu_i$ avec 
$\vert\mu_i\vert=1.$}.}
\end{thm}

\noindent \`A partir d'une transformation birationnelle $f$ sur le plan projectif complexe; peut-on trouver une suite finie d'\'eclatements $\pi\colon\mathcal{X}\to\mathbb{P}^2(\mathbb{C})$ telle que l'application induite $f_\mathcal{X}$  soit un automorphisme ? La r\'eponse est oui si $f$ est d'ordre fini, mais qu'en est-il si $f$ est d'ordre infini ? Si~$f$ n'est pas un automorphisme il existe une courbe $\mathcal{C}$ contract\'ee sur un point $p_1;$ la premi\`ere \'etape consiste \`a \'eclater le point $p_1$ via $\pi_1\colon\mathcal{X}_1\to\mathbb{P}^2(\mathbb{C}).$  Il se peut que tout se passe pour le mieux et que l'application induite $f_{\mathcal{X}_1}$ ne contracte pas $\mathcal{C}$ sur un point mais l'envoie sur $E_1;$ mais si $p_1$ n'est pas d'ind\'etermination alors $f_{\mathcal{X}_1}$ contracte $E_1$ sur $p_2=f(p_1)$ ... Ainsi cette proc\'edure qui consiste \`a \'eclater un nombre fini de points ne porte pas ses fruits si $f$ est alg\'ebriquement stable. 

\noindent Consid\'erons par exemple la famille de transformations birationnelles \'etudi\'ee dans \cite{BK1} donn\'ee par
\begin{align*}
& f\colon\mathbb{P}^2(\mathbb{C})\dashrightarrow\mathbb{P}^2(\mathbb{C}), && (x:y:z)\mapsto(x(bx+y):z(bx+y):x(ax+z)).
\end{align*}

\noindent La transformation $f$ contracte trois droites: $\Theta_0=\{x=0\},$ $\Theta_\beta=\{y+bx=0\}$ et $\Theta_\gamma=\{ax+z=~0\}$ sur respectivement $p_1=(0:1:0),$ $p_2=(0:0:1)$ et $q=(1:-a:0).$ Les points d'ind\'etermination sont $p_2,$ $p_1$ et $m=(1:-b:-a)$ qui sont respectivement \'eclat\'es sur $\Theta_0,$ $\Theta_B$ et $\Theta_C.$ Soit $\pi\colon \mathrm{Y}\to\mathbb{P}^2(\mathbb{C})$ l'\'eclatement de $p_1$ et $p_2.$ On peut v\'erifier que ni $\Theta_0,$ ni $\Theta_\beta$ sont exceptionnels pour $f_\mathrm{Y};$ plus pr\'ecis\'ement on a 
\begin{equation}\label{action}
\Theta_\beta\to E_2\to\Theta_0\to E_1\to\Theta_B=\{z=0\}.
\end{equation}
De plus $\Theta_\gamma$ est la seule courbe exceptionnelle et $m$ le seul point d'ind\'etermination de $f_\mathrm{Y}.$ D\'eterminons $f_{\mathrm{Y}*}.$ Une base de $\mathrm{Pic}(\mathrm{Y})$ est $\{\mathrm{L},\,E_1,\,E_2\}.$ La droite $\Theta_B$ contient un seul centre d'\'eclatement,~$p_1,$ donc $\Theta_B=\mathrm{L}-E_1.$  De m\^eme $\Theta_0$ contenant $p_1$ et $p_2$ on a $\Theta_0=\mathrm{L}-E_1-E_2.$ Une droite g\'en\'erique intersecte $\Theta_0,$ $\Theta_\beta$ et $\Theta_\gamma;$ donc l'image de $\mathrm{L}$ par $f$ est une conique passant par l'image de chacune de ces droites \`a savoir $p_1,$ $p_2$ et $q.$ Par suite $f_{\mathrm{Y}*}\mathrm{L}=2\mathrm{L}-E_1-E_2.$ \`A l'aide de (\ref{action}) on obtient que $f_{\mathrm{Y}*}$ est donn\'e par la matrice
$M=\left[\begin{array}{ccc}
2 & 1 & 1\\
-1 & -1 & -1\\
-1 & 0 & -1
\end{array}\right]$ dont le polyn\^ome caract\'eristique est $t^3-t-1.$ Soient $\mathrm{D}$ une droite et $\{\mathrm{D}\}$ sa classe dans $\mathrm{Pic}(\mathrm{Y}).$ Si $\mathrm{D}\cup f_\mathrm{Y}\mathrm{D}$ ne contient pas~$m,$ alors $f^2_\mathrm{Y}\mathrm{D}$ est donn\'e par $M^2\left[\begin{array}{ccc} 1 \\ 0\\0\end{array} \right]=\left[\begin{array}{ccc} 2 \\ 0\\-1\end{array} \right]=2\mathrm{L}-E_2,$ autrement dit $f_\mathrm{Y}^2\mathrm{D}$ est une conique qui intersecte $E_2$ mais pas $E_1.$ De m\^eme si $m$ n'appartient pas \`a $\mathrm{D}\cup f_\mathrm{Y}\mathrm{D}\cup f^2_\mathrm{Y}\mathrm{D}$ alors $f_\mathrm{Y}^3\mathrm{D}$ est une cubique intersectant $E_1$ et $E_2$ avec multiplicit\'e $1.$ Si $m$ n'appartient pas \`a~$\mathrm{D}\cup\ldots\cup f_\mathrm{Y}^{n-1}\mathrm{D},$ les it\'er\'es de $f_\mathrm{Y}$ sont holomorphes au voisinage de $\mathrm{D}$ et $(f_\mathrm{Y}^*)^n\mathrm{L}=\{f_\mathrm{Y}^n\mathrm{D}\}.$ Les param\`etres $a$ et~$b$ sont g\'en\'eriques si $m$ n'appartient pas \`a $\displaystyle\bigcup_{j=0}^\infty f_\mathrm{Y}^j\mathrm{D}.$ Pour $a,$ $b$ g\'en\'eriques on a $(f_\mathrm{Y}^*)^n=(f_\mathrm{Y}^n)^*$ et $\lambda(f_\mathrm{Y})$ est la plus grande racine en module du polyn\^ome $t^3-t-1.$ Soit $\mathcal{V}_n$ le sous-ensemble de $\mathbb{C}^2$ d\'efini par
\begin{align*}
\mathcal{V}_n=\{(a,b)\in\mathbb{C}^2\,\vert\, f_\mathrm{Y}^j(q)\not=m, \, 0\leq j<n,\, f_\mathrm{Y}^nq=m\}.
\end{align*}
Remarquons que si le couple $(a,b)$ n'appartient pas \`a $\mathcal{V}_n$ alors $\lambda(f_\mathrm{Y})$ est la plus grande racine en valeur absolue de $t^3-t-1$ dont on remarque qu'elle n'est pas un nombre de Salem. Il s'en suit que $f_\mathrm{Y}$ n'est pas un automorphisme modulo une suite finie d'\'eclatements. R\'eciproquement supposons que $(a,b)$ soit dans $\mathcal{V}_n;$ notons $\mathrm{Z}$ la surface obtenue en \'eclatant les points $q,$ $f_\mathrm{Y}(q),$ $\ldots,$ $f_\mathrm{Y}^n(q)=m.$ On peut montrer que la transformation induite $f_\mathrm{Z}$ est alors un automorphisme. Ainsi $f$ est un automorphisme sur un mod\`ele birationnel \`a $\mathbb{P}^2(\mathbb{C})$ si et seulement si $(a,b)$ appartient \`a $\mathcal{V}_n.$ En \'etudiant $f_{\mathrm{Z}*}$ on obtient que si $(a,b)$ appartient \`a $\mathcal{V}_n$ le polyn\^ome caract\'eristique de la matrice de $f_{\mathrm{Z}*}$ est $x^{n+1}(x^3-x-1)+x^3+x^2-1.$ Le premier degr\'e dynamique de $f$ est la plus grande racine en valeur absolue de ce polyn\^ome; il est strictement sup\'erieur \`a $1$ d\`es que $n\geq 7.$ Dans l'esprit de \cite{DJS} Bedford et Kim \'etudient les courbes invariantes par $f$ ainsi que la dynamique de $f$ (\emph{voir} \cite{BK1}).

\section{Le groupe de Cremona est-il simple ?}

\noindent Rappelons qu'un groupe est simple s'il ne poss\`ede pas
de sous-groupe distingu\'e non trivial. Commen\c{c}ons par remarquer que $\mathrm{Aut}[\mathbb{C}^2]$
n'est pas simple; soit $\phi$ le morphisme de~$\mathrm{Aut}[
\mathbb{C}^2]$ dans~$\mathbb{C}^*$ qui \`a $f$ associe $\det\mathrm{jac}\,
f.$ Le noyau de $\phi$ est un sous-groupe 
normal propre de $\mathrm{Aut}[\mathbb{C}^2].$ Dans les ann\'ees $1970$ 
Danilov \'etablit que $\ker\,\phi$
n'est pas simple (\cite{Da}). \`A l'aide de r\'esultats de Schupp
sur la th\'eorie des petites simplifications (\cite{Sc}) il montre que le sous-groupe 
normal\footnote{\hspace{1mm}Soient $\mathrm{G}$ un groupe et $f$ un 
\'el\'ement de $\mathrm{G};$ le sous-groupe normal 
engendr\'e par $f$ dans $\mathrm{G}$ est~$\langle hfh^{-1}\hspace{1mm}\vert\hspace{1mm} h\in~\mathrm{G}\rangle.$} engendr\'e par 
l'automorphisme
\begin{align*}
&(ea)^{13}, && a=(y,-x), && e=(x,y+3x^5-5x^4)
\end{align*}
\noindent est strictement contenu dans $\mathrm{Aut}[\mathbb{C}^2].$

\medskip

\noindent Plus r\'ecemment Furter et Lamy ont donn\'e un \'enonc\'e un 
peu plus pr\'ecis que celui de Danilov. Avant de l'\'enoncer introduisons une 
longueur $\ell(.)$ pour les \'el\'ements de $\mathrm{Aut}[\mathbb{C}^2].$
\begin{itemize}
\item Si $f$ appartient \`a $\mathtt{A}\cap\mathtt{E},$ alors $\ell(f)=0;$ 

\item sinon $\ell(f)$ est l'entier minimal $n$ pour lequel $f$ s'\'ecrit $g_1\ldots g_n$
avec $g_i$ dans $\mathtt{A}$ ou $\mathtt{E}.$
\end{itemize}
L'\'el\'ement exhib\'e par Danilov est de longueur $26.$

\noindent On note que $\ell(.)$ est invariant par conjugaison int\'erieure, on 
peut donc supposer~$f$ de longueur minimale dans sa
classe de conjugaison.

\begin{thm}[\cite{FL}]
{\sl Soit $f$ un \'el\'ement de $\mathrm{Aut}[\mathbb{C}^2]$ de d\'eterminant jacobien
$1.$ Supposons que~$f$ soit de longueur minimale dans sa classe de conjugaison.
\begin{itemize}
\item Si $f$ est non trivial et $\ell(f)\leq 8,$ le sous-groupe normal engendr\'e
par $f$ co\"incide avec le groupe des automorphismes polynomiaux du plan de 
d\'eterminant jacobien $1;$

\item si $f$ est g\'en\'erique\footnote{\hspace{1mm}On renvoie \`a \cite{FL} pour 
plus de d\'etail.} et $\ell(f)\geq 14,$ le sous-groupe normal engendr\'e par $f$ est 
strictement contenu dans le sous-groupe de $\mathrm{Aut}
[\mathbb{C}^2]$ form\'e des \'el\'ements de d\'eterminant jacobien $1.$
\end{itemize} }
\end{thm}

\bigskip

\noindent Le groupe de Cremona est-il simple ? Cette question a \'et\'e abord\'ee
dans \cite{Gi3, CD} mais reste ouverte.

\bigskip

\subsection*{Remerciements} 
Je tiens \`a remercier Zindine Djadli pour m'avoir demand\'e avec enthousiasme d'\'ecrire ce texte et de m'avoir permis de le faire sans contrainte. Merci \`a Dominique Cerveau pour nos conversations math\'ematiques aussi fr\'equentes qu'enrichissantes et pour ses nombreuses remarques. Je remercie Serge Cantat pour ses commentaires et pour m'avoir signal\'e comment am\'eliorer (\cite{De}, lemme~2.1) pour obtenir le lemme \ref{amel} et J\'er\'emy Blanc pour ses multiples remarques et suggestions. Merci \`a Julien Grivaux pour m'avoir aid\'ee \`a \og\'eclaircir\fg\, un paragraphe. Enfin merci \`a Igor Dolgachev pour m'avoir signal\'e les r\'ef\'erences \cite{Gi3, Gi4}.

\vspace*{8mm}

\bibliographystyle{alpha}
\bibliography{biblio}
\nocite{*}

\end{document}

%% file: arbre.pstex_t
\begin{picture}(0,0)%
\includegraphics{arbre.pstex}%
\end{picture}%
\setlength{\unitlength}{3947sp}%
\begingroup\makeatletter\ifx\SetFigFont\undefined%
\gdef\SetFigFont#1#2#3#4#5{%
  \reset@font\fontsize{#1}{#2pt}%
  \fontfamily{#3}\fontseries{#4}\fontshape{#5}%
  \selectfont}%
\fi\endgroup%
\begin{picture}(3474,2089)(-1136,-194)
\put(-974,1739){\makebox(0,0)[lb]{\smash{{\SetFigFont{12}{14.4}{\familydefault}{\mddefault}{\updefault}\begin{small}$ea\mathtt{E}$ \end{small}}}}}
\put(1951,1739){\makebox(0,0)[lb]{\smash{{\SetFigFont{12}{14.4}{\familydefault}{\mddefault}{\updefault}\begin{small}$ae\mathtt{A}$ \end{small}}}}}
\put(-899,464){\makebox(0,0)[lb]{\smash{{\SetFigFont{12}{14.4}{\familydefault}{\mddefault}{\updefault}\begin{small}$\widetilde{e}a\mathtt{E}$ \end{small}}}}}
\put(-974,1064){\makebox(0,0)[lb]{\smash{{\SetFigFont{12}{14.4}{\familydefault}{\mddefault}{\updefault}\begin{small}$e\widetilde{a}\mathtt{E}$ \end{small}}}}}
\put(-524,164){\makebox(0,0)[lb]{\smash{{\SetFigFont{12}{14.4}{\familydefault}{\mddefault}{\updefault}\begin{small}$\widetilde{e}\mathtt{A}$ \end{small}}}}}
\put(  1,914){\makebox(0,0)[lb]{\smash{{\SetFigFont{12}{14.4}{\familydefault}{\mddefault}{\updefault}\begin{small}$\mathrm{id}\mathtt{E}$ \end{small}}}}}
\put(-524,1439){\makebox(0,0)[lb]{\smash{{\SetFigFont{12}{14.4}{\familydefault}{\mddefault}{\updefault}\begin{small}$e\mathtt{A}$ \end{small}}}}}
\put(1051,914){\makebox(0,0)[lb]{\smash{{\SetFigFont{12}{14.4}{\familydefault}{\mddefault}{\updefault}\begin{small}$\mathrm{id}\mathtt{A}$ \end{small}}}}}
\put(1576,164){\makebox(0,0)[lb]{\smash{{\SetFigFont{12}{14.4}{\familydefault}{\mddefault}{\updefault}\begin{small}$\widetilde{a}\mathtt{E}$ \end{small}}}}}
\put(-974,-136){\makebox(0,0)[lb]{\smash{{\SetFigFont{12}{14.4}{\familydefault}{\mddefault}{\updefault}\begin{small}$\widetilde{e}\widetilde{a}\mathtt{E}$ \end{small}}}}}
\put(1951,-136){\makebox(0,0)[lb]{\smash{{\SetFigFont{12}{14.4}{\familydefault}{\mddefault}{\updefault}\begin{small}$\widetilde{a}\widetilde{e}\mathtt{A}$ \end{small}}}}}
\put(1876,464){\makebox(0,0)[lb]{\smash{{\SetFigFont{12}{14.4}{\familydefault}{\mddefault}{\updefault}\begin{small}$\widetilde{a}e\mathtt{A}$ \end{small}}}}}
\put(1951,1064){\makebox(0,0)[lb]{\smash{{\SetFigFont{12}{14.4}{\familydefault}{\mddefault}{\updefault}\begin{small}$a\widetilde{e}\mathtt{A}$ \end{small}}}}}
\put(1576,1439){\makebox(0,0)[lb]{\smash{{\SetFigFont{12}{14.4}{\familydefault}{\mddefault}{\updefault}\begin{small}$a\mathtt{E}$ \end{small}}}}}
\end{picture}%

%% file: decomp.pstex_t
\begin{picture}(0,0)%
\includegraphics{decomp.pstex}%
\end{picture}%
\setlength{\unitlength}{3947sp}%
\begingroup\makeatletter\ifx\SetFigFont\undefined%
\gdef\SetFigFont#1#2#3#4#5{%
  \reset@font\fontsize{#1}{#2pt}%
  \fontfamily{#3}\fontseries{#4}\fontshape{#5}%
  \selectfont}%
\fi\endgroup%
\begin{picture}(6687,3214)(1426,-12344)
\put(7351,-11461){\makebox(0,0)[lb]{\smash{{\SetFigFont{12}{14.4}{\familydefault}{\mddefault}{\updefault}\begin{small}$\mathbb{P}^2(\mathbb{C})$ \end{small}}}}}
\put(7426,-11011){\makebox(0,0)[lb]{\smash{{\SetFigFont{12}{14.4}{\familydefault}{\mddefault}{\updefault}\begin{small}$\widetilde{L}_{AB}$ \end{small}}}}}
\put(7801,-10336){\makebox(0,0)[lb]{\smash{{\SetFigFont{12}{14.4}{\familydefault}{\mddefault}{\updefault}\begin{small}$E_B$ \end{small}}}}}
\put(7051,-10336){\makebox(0,0)[lb]{\smash{{\SetFigFont{12}{14.4}{\familydefault}{\mddefault}{\updefault}\begin{small}$E_A$ \end{small}}}}}
\put(7351,-9586){\makebox(0,0)[lb]{\smash{{\SetFigFont{12}{14.4}{\familydefault}{\mddefault}{\updefault}\begin{small}$E_C$ \end{small}}}}}
\put(6676,-9586){\makebox(0,0)[lb]{\smash{{\SetFigFont{12}{14.4}{\familydefault}{\mddefault}{\updefault}\begin{small}$\widetilde{L}_{AC}$ \end{small}}}}}
\put(8101,-9586){\makebox(0,0)[lb]{\smash{{\SetFigFont{12}{14.4}{\familydefault}{\mddefault}{\updefault}\begin{small}$\widetilde{L}_{BC}$ \end{small}}}}}
\put(4951,-9886){\makebox(0,0)[lb]{\smash{{\SetFigFont{12}{14.4}{\familydefault}{\mddefault}{\updefault}\begin{small}$\widetilde{L}_{BC}$ \end{small}}}}}
\put(4651,-10786){\makebox(0,0)[lb]{\smash{{\SetFigFont{12}{14.4}{\familydefault}{\mddefault}{\updefault}\begin{small}$\widetilde{L}_{AC}$ \end{small}}}}}
\put(5176,-11011){\makebox(0,0)[lb]{\smash{{\SetFigFont{12}{14.4}{\familydefault}{\mddefault}{\updefault}\begin{small}$E_C$ \end{small}}}}}
\put(4051,-10936){\makebox(0,0)[lb]{\smash{{\SetFigFont{12}{14.4}{\familydefault}{\mddefault}{\updefault}\begin{small}$E_A$ \end{small}}}}}
\put(4651,-9286){\makebox(0,0)[lb]{\smash{{\SetFigFont{12}{14.4}{\familydefault}{\mddefault}{\updefault}\begin{small}$E_B$ \end{small}}}}}
\put(4351,-9511){\makebox(0,0)[lb]{\smash{{\SetFigFont{12}{14.4}{\familydefault}{\mddefault}{\updefault}\begin{small}$\widetilde{L}_{AB}$ \end{small}}}}}
\put(2026,-11461){\makebox(0,0)[lb]{\smash{{\SetFigFont{12}{14.4}{\familydefault}{\mddefault}{\updefault}\begin{small}$\mathbb{P}^2(\mathbb{C})$\end{small} }}}}
\put(2176,-9436){\makebox(0,0)[lb]{\smash{{\SetFigFont{12}{14.4}{\familydefault}{\mddefault}{\updefault}\begin{small}$C$ \end{small}}}}}
\put(1576,-10261){\makebox(0,0)[lb]{\smash{{\SetFigFont{12}{14.4}{\familydefault}{\mddefault}{\updefault}\begin{small}$L_{AC}$ \end{small}}}}}
\put(2026,-10861){\makebox(0,0)[lb]{\smash{{\SetFigFont{12}{14.4}{\familydefault}{\mddefault}{\updefault}\begin{small}$L_{AB}$ \end{small}}}}}
\put(1426,-10861){\makebox(0,0)[lb]{\smash{{\SetFigFont{12}{14.4}{\familydefault}{\mddefault}{\updefault}\begin{small}$A$ \end{small}}}}}
\put(2776,-10861){\makebox(0,0)[lb]{\smash{{\SetFigFont{12}{14.4}{\familydefault}{\mddefault}{\updefault}\begin{small}$B$ \end{small}}}}}
\put(4801,-12286){\makebox(0,0)[lb]{\smash{{\SetFigFont{12}{14.4}{\familydefault}{\mddefault}{\updefault}\begin{small}$\sigma$ \end{small}}}}}
\put(2476,-10261){\makebox(0,0)[lb]{\smash{{\SetFigFont{12}{14.4}{\familydefault}{\mddefault}{\updefault}\begin{small}$L_{BC}$ \end{small}}}}}
\put(3451,-10036){\makebox(0,0)[lb]{\smash{{\SetFigFont{12}{14.4}{\familydefault}{\mddefault}{\updefault}\begin{small}$\pi_1$ \end{small}}}}}
\put(6076,-10036){\makebox(0,0)[lb]{\smash{{\SetFigFont{12}{14.4}{\familydefault}{\mddefault}{\updefault}\begin{small}$\pi_2$ \end{small}}}}}
\end{picture}%

%% file: zariski.pstex_t
\begin{picture}(0,0)%
\includegraphics{zariski.pstex}%
\end{picture}%
\setlength{\unitlength}{3947sp}%
\begingroup\makeatletter\ifx\SetFigFont\undefined%
\gdef\SetFigFont#1#2#3#4#5{%
  \reset@font\fontsize{#1}{#2pt}%
  \fontfamily{#3}\fontseries{#4}\fontshape{#5}%
  \selectfont}%
\fi\endgroup%
\begin{picture}(5025,3364)(1726,-15944)
\put(2551,-15886){\makebox(0,0)[lb]{\smash{{\SetFigFont{12}{14.4}{\familydefault}{\mddefault}{\updefault}\begin{small}$\phi^ {-1}(\mathcal{C}_2)$ \end{small}}}}}
\put(4201,-15436){\makebox(0,0)[lb]{\smash{{\SetFigFont{12}{14.4}{\familydefault}{\mddefault}{\updefault}\begin{small}$f$ \end{small}}}}}
\put(5401,-14161){\makebox(0,0)[lb]{\smash{{\SetFigFont{12}{14.4}{\familydefault}{\mddefault}{\updefault}\begin{small}$\pi$ \end{small}}}}}
\put(3076,-14161){\makebox(0,0)[lb]{\smash{{\SetFigFont{12}{14.4}{\familydefault}{\mddefault}{\updefault}\begin{small}$\phi$ \end{small}}}}}
\put(4801,-13336){\makebox(0,0)[lb]{\smash{{\SetFigFont{12}{14.4}{\familydefault}{\mddefault}{\updefault}\begin{small}$q_2$ \end{small}}}}}
\put(5026,-13111){\makebox(0,0)[lb]{\smash{{\SetFigFont{12}{14.4}{\familydefault}{\mddefault}{\updefault}\begin{small}$E$\end{small} }}}}
\put(3901,-12736){\makebox(0,0)[lb]{\smash{{\SetFigFont{12}{14.4}{\familydefault}{\mddefault}{\updefault}\begin{small}$\mathcal{C}_1$ \end{small}}}}}
\put(4876,-12736){\makebox(0,0)[lb]{\smash{{\SetFigFont{12}{14.4}{\familydefault}{\mddefault}{\updefault}\begin{small}$\mathcal{C}_2$ \end{small}}}}}
\put(6001,-15811){\makebox(0,0)[lb]{\smash{{\SetFigFont{12}{14.4}{\familydefault}{\mddefault}{\updefault}\begin{small}$\pi(\mathcal{C}_2)$ \end{small}}}}}
\put(2851,-15361){\makebox(0,0)[lb]{\smash{{\SetFigFont{12}{14.4}{\familydefault}{\mddefault}{\updefault}\begin{small}$m$ \end{small}}}}}
\put(1726,-15361){\makebox(0,0)[lb]{\smash{{\SetFigFont{12}{14.4}{\familydefault}{\mddefault}{\updefault}\begin{small}$\phi^ {-1}(\mathcal{C}_1)$ \end{small}}}}}
\put(6001,-15061){\makebox(0,0)[lb]{\smash{{\SetFigFont{12}{14.4}{\familydefault}{\mddefault}{\updefault}\begin{small}$p=f(m)$ \end{small}}}}}
\put(5251,-15586){\makebox(0,0)[lb]{\smash{{\SetFigFont{12}{14.4}{\familydefault}{\mddefault}{\updefault}\begin{small}$\pi(\mathcal{C}_1)$ \end{small}}}}}
\put(3826,-13261){\makebox(0,0)[lb]{\smash{{\SetFigFont{12}{14.4}{\familydefault}{\mddefault}{\updefault}\begin{small}$q_1$\end{small} }}}}
\put(1726,-14836){\makebox(0,0)[lb]{\smash{{\SetFigFont{12}{14.4}{\familydefault}{\mddefault}{\updefault}\begin{small}$\mathrm{S}$ \end{small}}}}}
\put(6751,-14761){\makebox(0,0)[lb]{\smash{{\SetFigFont{12}{14.4}{\familydefault}{\mddefault}{\updefault}\begin{small}$\mathrm{S}'$ \end{small}}}}}
\put(3076,-12961){\makebox(0,0)[lb]{\smash{{\SetFigFont{12}{14.4}{\familydefault}{\mddefault}{\updefault}\begin{small}$\widetilde{\mathrm{S}}$ \end{small}}}}}
\end{picture}%

%% file: el1.pstex_t
\begin{picture}(0,0)%
\includegraphics{el1.pstex}%
\end{picture}%
\setlength{\unitlength}{3947sp}%
\begingroup\makeatletter\ifx\SetFigFont\undefined%
\gdef\SetFigFont#1#2#3#4#5{%
  \reset@font\fontsize{#1}{#2pt}%
  \fontfamily{#3}\fontseries{#4}\fontshape{#5}%
  \selectfont}%
\fi\endgroup%
\begin{picture}(5049,1714)(2239,-11744)
\put(7126,-11311){\makebox(0,0)[lb]{\smash{{\SetFigFont{12}{14.4}{\familydefault}{\mddefault}{\updefault}\begin{small}$s_{n+1}$ \end{small}}}}}
\put(7276,-10861){\makebox(0,0)[lb]{\smash{{\SetFigFont{12}{14.4}{\familydefault}{\mddefault}{\updefault}\begin{small}$-(n+1)$\end{small} }}}}
\put(6826,-10336){\makebox(0,0)[lb]{\smash{{\SetFigFont{12}{14.4}{\familydefault}{\mddefault}{\updefault}\begin{small}$0$\end{small} }}}}
\put(6001,-10636){\makebox(0,0)[lb]{\smash{{\SetFigFont{12}{14.4}{\familydefault}{\mddefault}{\updefault}\begin{small}$\pi_2$ \end{small}}}}}
\put(5251,-11011){\makebox(0,0)[lb]{\smash{{\SetFigFont{12}{14.4}{\familydefault}{\mddefault}{\updefault}\begin{small}$-(n+1)$ \end{small}}}}}
\put(5101,-11461){\makebox(0,0)[lb]{\smash{{\SetFigFont{12}{14.4}{\familydefault}{\mddefault}{\updefault}\begin{small}$\widetilde{s_n}$ \end{small}}}}}
\put(5101,-10411){\makebox(0,0)[lb]{\smash{{\SetFigFont{12}{14.4}{\familydefault}{\mddefault}{\updefault}\begin{small}$-1$ \end{small}}}}}
\put(4351,-10186){\makebox(0,0)[lb]{\smash{{\SetFigFont{12}{14.4}{\familydefault}{\mddefault}{\updefault}\begin{small}$-1$ \end{small}}}}}
\put(4276,-10561){\makebox(0,0)[lb]{\smash{{\SetFigFont{12}{14.4}{\familydefault}{\mddefault}{\updefault}\begin{small}$\widetilde{f}$ \end{small}}}}}
\put(3676,-10636){\makebox(0,0)[lb]{\smash{{\SetFigFont{12}{14.4}{\familydefault}{\mddefault}{\updefault}\begin{small}$\pi_1$ \end{small}}}}}
\put(2851,-11311){\makebox(0,0)[lb]{\smash{{\SetFigFont{12}{14.4}{\familydefault}{\mddefault}{\updefault}\begin{small}$s_n$ \end{small}}}}}
\put(3001,-10861){\makebox(0,0)[lb]{\smash{{\SetFigFont{12}{14.4}{\familydefault}{\mddefault}{\updefault}\begin{small}$-n$ \end{small}}}}}
\put(2776,-10636){\makebox(0,0)[lb]{\smash{{\SetFigFont{12}{14.4}{\familydefault}{\mddefault}{\updefault}\begin{small}$p$ \end{small}}}}}
\put(2551,-10336){\makebox(0,0)[lb]{\smash{{\SetFigFont{12}{14.4}{\familydefault}{\mddefault}{\updefault}\begin{small}$0$ \end{small}}}}}
\put(2326,-10636){\makebox(0,0)[lb]{\smash{{\SetFigFont{12}{14.4}{\familydefault}{\mddefault}{\updefault}\begin{small}$f$ \end{small}}}}}
\put(6676,-11686){\makebox(0,0)[lb]{\smash{{\SetFigFont{12}{14.4}{\familydefault}{\mddefault}{\updefault}\begin{small}$\mathrm{F}_{n+1}$\end{small} }}}}
\put(2476,-11686){\makebox(0,0)[lb]{\smash{{\SetFigFont{12}{14.4}{\familydefault}{\mddefault}{\updefault}\begin{small}$\mathrm{F}_n$\end{small} }}}}
\end{picture}%

%% file: el2.pstex_t
\begin{picture}(0,0)%
\includegraphics{el2.pstex}%
\end{picture}%
\setlength{\unitlength}{3947sp}%
\begingroup\makeatletter\ifx\SetFigFont\undefined%
\gdef\SetFigFont#1#2#3#4#5{%
  \reset@font\fontsize{#1}{#2pt}%
  \fontfamily{#3}\fontseries{#4}\fontshape{#5}%
  \selectfont}%
\fi\endgroup%
\begin{picture}(5530,2020)(2176,-11744)
\put(7126,-11311){\makebox(0,0)[lb]{\smash{{\SetFigFont{12}{14.4}{\familydefault}{\mddefault}{\updefault}\begin{small}$s_n$ \end{small}}}}}
\put(5251,-11386){\makebox(0,0)[lb]{\smash{{\SetFigFont{12}{14.4}{\familydefault}{\mddefault}{\updefault}\begin{small}$\widetilde{s_{n+1}}$ \end{small}}}}}
\put(5476,-11011){\makebox(0,0)[lb]{\smash{{\SetFigFont{12}{14.4}{\familydefault}{\mddefault}{\updefault}\begin{small}$-(n+1)$ \end{small}}}}}
\put(7276,-10861){\makebox(0,0)[lb]{\smash{{\SetFigFont{12}{14.4}{\familydefault}{\mddefault}{\updefault}$-n$ }}}}
\put(6076,-10636){\makebox(0,0)[lb]{\smash{{\SetFigFont{12}{14.4}{\familydefault}{\mddefault}{\updefault}\begin{small}$\pi_2$\end{small} }}}}
\put(4726,-10711){\makebox(0,0)[lb]{\smash{{\SetFigFont{12}{14.4}{\familydefault}{\mddefault}{\updefault}\begin{small}$\widetilde{f}$\end{small} }}}}
\put(6826,-10336){\makebox(0,0)[lb]{\smash{{\SetFigFont{12}{14.4}{\familydefault}{\mddefault}{\updefault}\begin{small}$0$ \end{small}}}}}
\put(4951,-10336){\makebox(0,0)[lb]{\smash{{\SetFigFont{12}{14.4}{\familydefault}{\mddefault}{\updefault}\begin{small}$-1$ \end{small}}}}}
\put(4426,-10111){\makebox(0,0)[lb]{\smash{{\SetFigFont{12}{14.4}{\familydefault}{\mddefault}{\updefault}\begin{small}$-1$ \end{small}}}}}
\put(3826,-10636){\makebox(0,0)[lb]{\smash{{\SetFigFont{12}{14.4}{\familydefault}{\mddefault}{\updefault}\begin{small}$\pi_1$ \end{small}}}}}
\put(3001,-11011){\makebox(0,0)[lb]{\smash{{\SetFigFont{12}{14.4}{\familydefault}{\mddefault}{\updefault}\begin{small}$-(n+1)$ \end{small}}}}}
\put(2851,-11311){\makebox(0,0)[lb]{\smash{{\SetFigFont{12}{14.4}{\familydefault}{\mddefault}{\updefault}\begin{small}$s_{n+1}$ \end{small}}}}}
\put(2701,-10636){\makebox(0,0)[lb]{\smash{{\SetFigFont{12}{14.4}{\familydefault}{\mddefault}{\updefault}\begin{small}$p$ \end{small}}}}}
\put(2551,-10336){\makebox(0,0)[lb]{\smash{{\SetFigFont{12}{14.4}{\familydefault}{\mddefault}{\updefault}\begin{small}$0$ \end{small}}}}}
\put(2176,-10636){\makebox(0,0)[lb]{\smash{{\SetFigFont{12}{14.4}{\familydefault}{\mddefault}{\updefault}\begin{small}$f$ \end{small}}}}}
\put(6751,-11686){\makebox(0,0)[lb]{\smash{{\SetFigFont{12}{14.4}{\familydefault}{\mddefault}{\updefault}\begin{small}$\mathrm{F}_n$ \end{small}}}}}
\put(2401,-11686){\makebox(0,0)[lb]{\smash{{\SetFigFont{12}{14.4}{\familydefault}{\mddefault}{\updefault}\begin{small}$\mathrm{F}_{n+1}$\end{small} }}}}
\end{picture}%

%% file: as.pstex_t
\begin{picture}(0,0)%
\includegraphics{as.pstex}%
\end{picture}%
\setlength{\unitlength}{3947sp}%
\begingroup\makeatletter\ifx\SetFigFont\undefined%
\gdef\SetFigFont#1#2#3#4#5{%
  \reset@font\fontsize{#1}{#2pt}%
  \fontfamily{#3}\fontseries{#4}\fontshape{#5}%
  \selectfont}%
\fi\endgroup%
\begin{picture}(5337,895)(1051,-8444)
\put(5776,-7711){\makebox(0,0)[lb]{\smash{{\SetFigFont{12}{14.4}{\familydefault}{\mddefault}{\updefault}\begin{small}$f$ \end{small}}}}}
\put(4951,-7711){\makebox(0,0)[lb]{\smash{{\SetFigFont{12}{14.4}{\familydefault}{\mddefault}{\updefault}\begin{small}$f$ \end{small}}}}}
\put(3226,-7711){\makebox(0,0)[lb]{\smash{{\SetFigFont{12}{14.4}{\familydefault}{\mddefault}{\updefault}\begin{small}$f$ \end{small}}}}}
\put(2401,-7711){\makebox(0,0)[lb]{\smash{{\SetFigFont{12}{14.4}{\familydefault}{\mddefault}{\updefault}\begin{small}$f$ \end{small}}}}}
\put(1576,-7711){\makebox(0,0)[lb]{\smash{{\SetFigFont{12}{14.4}{\familydefault}{\mddefault}{\updefault}\begin{small}$f$ \end{small}}}}}
\put(1051,-8386){\makebox(0,0)[lb]{\smash{{\SetFigFont{12}{14.4}{\familydefault}{\mddefault}{\updefault}\begin{small}$\mathcal{C}$ \end{small}}}}}
\end{picture}%

%% file: as4.pstex_t
\begin{picture}(0,0)%
\includegraphics{as4.pstex}%
\end{picture}%
\setlength{\unitlength}{3947sp}%
\begingroup\makeatletter\ifx\SetFigFont\undefined%
\gdef\SetFigFont#1#2#3#4#5{%
  \reset@font\fontsize{#1}{#2pt}%
  \fontfamily{#3}\fontseries{#4}\fontshape{#5}%
  \selectfont}%
\fi\endgroup%
\begin{picture}(5337,895)(1051,-8444)
\put(1051,-8386){\makebox(0,0)[lb]{\smash{{\SetFigFont{12}{14.4}{\familydefault}{\mddefault}{\updefault}\begin{small}$\mathcal{C}$\end{small}}}}}
\put(1951,-8086){\makebox(0,0)[lb]{\smash{{\SetFigFont{12}{14.4}{\familydefault}{\mddefault}{\updefault}\begin{small}$p_1$ \end{small}}}}}
\put(2776,-8086){\makebox(0,0)[lb]{\smash{{\SetFigFont{12}{14.4}{\familydefault}{\mddefault}{\updefault}\begin{small}$p_2$ \end{small}}}}}
\put(4501,-8086){\makebox(0,0)[lb]{\smash{{\SetFigFont{12}{14.4}{\familydefault}{\mddefault}{\updefault}\begin{small}$p_{k-1}$ \end{small}}}}}
\put(5326,-8086){\makebox(0,0)[lb]{\smash{{\SetFigFont{12}{14.4}{\familydefault}{\mddefault}{\updefault}\begin{small}$p_k$ \end{small}}}}}
\end{picture}%

%% file: sl.pstex_t
\begin{picture}(0,0)%
\includegraphics{sl.pstex}%
\end{picture}%
\setlength{\unitlength}{3947sp}%
\begingroup\makeatletter\ifx\SetFigFont\undefined%
\gdef\SetFigFont#1#2#3#4#5{%
  \reset@font\fontsize{#1}{#2pt}%
  \fontfamily{#3}\fontseries{#4}\fontshape{#5}%
  \selectfont}%
\fi\endgroup%
\begin{picture}(1200,1339)(1126,-1619)
\put(1951,-436){\makebox(0,0)[lb]{\smash{{\SetFigFont{12}{14.4}{\familydefault}{\mddefault}{\updefault}\begin{small}$\mathrm{e}_{13}$\end{small}}}}}
\put(1351,-436){\makebox(0,0)[lb]{\smash{{\SetFigFont{12}{14.4}{\familydefault}{\mddefault}{\updefault}\begin{small}$\mathrm{e}_{23}$\end{small}}}}}
\put(1951,-1561){\makebox(0,0)[lb]{\smash{{\SetFigFont{12}{14.4}{\familydefault}{\mddefault}{\updefault}\begin{small}$\mathrm{e}_{32}$\end{small}}}}}
\put(2326,-1036){\makebox(0,0)[lb]{\smash{{\SetFigFont{12}{14.4}{\familydefault}{\mddefault}{\updefault}\begin{small}$\mathrm{e}_{12}$\end{small}}}}}
\put(1126,-1036){\makebox(0,0)[lb]{\smash{{\SetFigFont{12}{14.4}{\familydefault}{\mddefault}{\updefault}\begin{small}$\mathrm{e}_{21}$\end{small}}}}}
\put(1426,-1561){\makebox(0,0)[lb]{\smash{{\SetFigFont{12}{14.4}{\familydefault}{\mddefault}{\updefault}\begin{small}$\mathrm{e}_{31}$\end{small}}}}}
\end{picture}%

%% file: pingpong2.pstex_t
\begin{picture}(0,0)%
\includegraphics{pingpong2.pstex}%
\end{picture}%
\setlength{\unitlength}{3947sp}%
\begingroup\makeatletter\ifx\SetFigFont\undefined%
\gdef\SetFigFont#1#2#3#4#5{%
  \reset@font\fontsize{#1}{#2pt}%
  \fontfamily{#3}\fontseries{#4}\fontshape{#5}%
  \selectfont}%
\fi\endgroup%
\begin{picture}(1974,2133)(2389,-7573)
\put(4276,-6811){\makebox(0,0)[lb]{\smash{{\SetFigFont{12}{14.4}{\familydefault}{\mddefault}{\updefault}$\mathrm{X}$}}}}
\put(2551,-7186){\makebox(0,0)[lb]{\smash{{\SetFigFont{12}{14.4}{\familydefault}{\mddefault}{\updefault}$\Gamma_2$}}}}
\put(3676,-5611){\makebox(0,0)[lb]{\smash{{\SetFigFont{12}{14.4}{\familydefault}{\mddefault}{\updefault}$\alpha_i$}}}}
\put(3151,-5986){\makebox(0,0)[lb]{\smash{{\SetFigFont{12}{14.4}{\familydefault}{\mddefault}{\updefault}$\Gamma_1$}}}}
\put(2926,-7486){\makebox(0,0)[lb]{\smash{{\SetFigFont{12}{14.4}{\familydefault}{\mddefault}{\updefault}$\beta_i$}}}}
\end{picture}%

%% file: aut.pstex_t
\begin{picture}(0,0)%
\includegraphics{aut.pstex}%
\end{picture}%
\setlength{\unitlength}{3947sp}%
\begingroup\makeatletter\ifx\SetFigFont\undefined%
\gdef\SetFigFont#1#2#3#4#5{%
  \reset@font\fontsize{#1}{#2pt}%
  \fontfamily{#3}\fontseries{#4}\fontshape{#5}%
  \selectfont}%
\fi\endgroup%
\begin{picture}(774,664)(2239,-13619)
\put(2326,-13111){\makebox(0,0)[lb]{\smash{{\SetFigFont{12}{14.4}{\familydefault}{\mddefault}{\updefault}\begin{small}$p$ \end{small}}}}}
\put(2851,-13561){\makebox(0,0)[lb]{\smash{{\SetFigFont{12}{14.4}{\familydefault}{\mddefault}{\updefault}\begin{small}$\Theta$\end{small} }}}}
\end{picture}%

%% file: aut2.pstex_t
\begin{picture}(0,0)%
\includegraphics{aut2.pstex}%
\end{picture}%
\setlength{\unitlength}{3947sp}%
\begingroup\makeatletter\ifx\SetFigFont\undefined%
\gdef\SetFigFont#1#2#3#4#5{%
  \reset@font\fontsize{#1}{#2pt}%
  \fontfamily{#3}\fontseries{#4}\fontshape{#5}%
  \selectfont}%
\fi\endgroup%
\begin{picture}(774,993)(4489,-13648)
\put(5176,-12811){\makebox(0,0)[lb]{\smash{{\SetFigFont{12}{14.4}{\familydefault}{\mddefault}{\updefault}\begin{small}$\mathrm{E}_1$\end{small} }}}}
\put(4501,-13561){\makebox(0,0)[lb]{\smash{{\SetFigFont{12}{14.4}{\familydefault}{\mddefault}{\updefault}\begin{small}$\Theta$\end{small} }}}}
\end{picture}%

%% file: aut3.pstex_t
\begin{picture}(0,0)%
\includegraphics{aut3.pstex}%
\end{picture}%
\setlength{\unitlength}{3947sp}%
\begingroup\makeatletter\ifx\SetFigFont\undefined%
\gdef\SetFigFont#1#2#3#4#5{%
  \reset@font\fontsize{#1}{#2pt}%
  \fontfamily{#3}\fontseries{#4}\fontshape{#5}%
  \selectfont}%
\fi\endgroup%
\begin{picture}(837,1743)(7189,-14023)
\put(7726,-12436){\makebox(0,0)[lb]{\smash{{\SetFigFont{12}{14.4}{\familydefault}{\mddefault}{\updefault}\begin{small}$\mathrm{E}_2$\end{small} }}}}
\put(8026,-12811){\makebox(0,0)[lb]{\smash{{\SetFigFont{12}{14.4}{\familydefault}{\mddefault}{\updefault}\begin{small}$\mathrm{E}_1$\end{small} }}}}
\put(8026,-13411){\makebox(0,0)[lb]{\smash{{\SetFigFont{12}{14.4}{\familydefault}{\mddefault}{\updefault}\begin{small}$\Theta$\end{small} }}}}
\end{picture}%

%% file: aut4.pstex_t
\begin{picture}(0,0)%
\includegraphics{aut4.pstex}%
\end{picture}%
\setlength{\unitlength}{3947sp}%
\begingroup\makeatletter\ifx\SetFigFont\undefined%
\gdef\SetFigFont#1#2#3#4#5{%
  \reset@font\fontsize{#1}{#2pt}%
  \fontfamily{#3}\fontseries{#4}\fontshape{#5}%
  \selectfont}%
\fi\endgroup%
\begin{picture}(837,1818)(9814,-14098)
\put(10351,-12436){\makebox(0,0)[lb]{\smash{{\SetFigFont{12}{14.4}{\familydefault}{\mddefault}{\updefault}\begin{small}$\mathrm{E}_2$\end{small} }}}}
\put(10651,-12811){\makebox(0,0)[lb]{\smash{{\SetFigFont{12}{14.4}{\familydefault}{\mddefault}{\updefault}\begin{small}$\mathrm{E}_1$\end{small} }}}}
\put(10651,-13861){\makebox(0,0)[lb]{\smash{{\SetFigFont{12}{14.4}{\familydefault}{\mddefault}{\updefault}\begin{small}$\mathrm{E}_3$\end{small} }}}}
\put(10651,-13336){\makebox(0,0)[lb]{\smash{{\SetFigFont{12}{14.4}{\familydefault}{\mddefault}{\updefault}\begin{small}$\Theta$\end{small} }}}}
\end{picture}%